\documentclass[sort,11pt]{elsarticle}

\usepackage{amssymb}
\usepackage{amsmath}

\usepackage{lineno}

\usepackage{amsthm}
\usepackage{thmtools}
\newtheorem{theorem}{Theorem}
\newtheorem{remark}{Remark}
\newtheorem{assumption}{Assumption}
\newtheorem{lemma}[theorem]{Lemma}
\newtheorem{corollary}[theorem]{Corollary}

\usepackage{hyperref}
\usepackage{booktabs}
\usepackage{graphicx}
\usepackage{subcaption}
\usepackage{adjustbox}
\usepackage{color}

\usepackage{multirow} %

\usepackage{cleveref}
\usepackage{algorithm}
\usepackage{algorithmicx}
\usepackage{algpseudocode}
\algrenewcommand\algorithmicrequire{\textbf{Input:}}
\algrenewcommand\algorithmicensure{\textbf{Output:}}

\crefname{assumption}{Assumption}{Assumptions}
\Crefname{assumption}{Assumption}{Assumptions}
\crefname{lemma}{Lemma}{Lemmas}
\Crefname{lemma}{Lemma}{Lemmas}
\crefname{corollary}{Corollary}{Corollaries}
\Crefname{corollary}{Corollary}{Corollaries}
\crefname{appendix}{}{}
\Crefname{appendix}{}{}
\crefname{equation}{}{}
\Crefname{equation}{}{}

\newcommand{\R}{\mathbb{R}}
\newcommand{\di}{\,\mathrm{d}}

\newcommand{\X}{\mathcal{X}}

\newcommand{\M}{\mathcal{M}}

\newcommand{\mO}{\mathcal{O}}

\newcommand{\E}[1]{\mathbb{E}\left[#1\right]}

\newcommand{\Etx}[1]{\mathbb{E}_{t}^{x}\left[#1\right]}

\newcommand{\Etnx}[1]{\mathbb{E}_{t_n}^{x}\left[#1\right]}

\newcommand{\br}[1]{\left(#1\right)}
\newcommand{\fbr}[1]{(#1)}
\newcommand{\abs}[1]{\left\vert#1\right\vert}

\newcommand{\bbr}[1]{\left\{#1\right\}}

\newcommand{\vbr}[2]{#1(#2 #1)}
\newcommand{\rbr}[1]{\left[#1\right]}
\newcommand{\Pra}[1]{\mathbb{P}\left(#1\right)}

\newcommand{\argmin}{\mathop{\arg\min}}

\newcommand{\B}[1]{\boldsymbol{#1}}

\newcommand{\patho}{CodeAndFig0812/Subsection_4.1/CDE_d1e3---m122/outputs/}
\newcommand{\pathsocmnm}{CodeAndFig0812/Subsection_4.2/QLP-1_-2a_d1e4---m177/outputs/}

\newcommand{\pathsocmnl}{CodeAndFig0812/Subsection_4.2/QLP-2bv_d1e4---m179/outputs/}
\newcommand{\pathq}{CodeAndFig0812/Subsection_4.2/QLP-1_2a_d1e5---m183/outputs/}

\newcommand{\patht}{CodeAndFig0812/Subsection_4.2/QLP-2bv_d1e5---m187/outputs/}

\newcommand{\pathj}{CodeAndFig0812/Subsection_4.3/HJB1a---m92/outputs/}
\newcommand{\pathl}{CodeAndFig0812/Subsection_4.3/HJB1b---m94/outputs/}
\newcommand{\paths}{CodeAndFig0812/Subsection_4.3/HJB2---Code-is-in-subfolder-result---m165and169/}
\newcommand{\pathp}{CodeAndFig0812/Subsection_4.4/CDE_d1---m126/outputs/}

\journal{---}

\begin{document}

\begin{frontmatter}



	\title{Deep random difference method for high-dimensional quasilinear parabolic partial differential equations\tnoteref{t1}}
 \tnotetext[t1]{
    The work of the first author is supported by the Clements Chair in Applied Mathematics at SMU.
    The work of the second author is supported by the Postdoctoral Fellowship Program of CPSF (under grant GZB20250705) and China Postdoctoral Science 
    Foundation (under grant 2025M783150). 
    The work of the third author is supported by the NSF of China (under grant 12288201), the NSFC-RGC joint project (under grant 12461160275) and the Youth Innovation Promotion Association (CAS).
    }

\author[author1]{Wei Cai} 
\author[author2]{Shuixin Fang} 
\author[author2]{Tao Zhou\corref{cor}} 

\cortext[cor]{Corresponding author. }
\ead{tzhou@lsec.cc.ac.cn}

\affiliation[author1]{organization={Department of Mathematics, Southern Methodist University},
addressline={Clements Hall, 3100 Dyer St}, 
city={Dallas},
postcode={TX 75205}, 
state={Texas},
country={USA}}

\affiliation[author2]{organization={Institute of Computational Mathematics and Scientific/Engineering Computing,
Academy of Mathematics and Systems Science, Chinese Academy of Sciences},
addressline={No.55 Zhong Guan Cun East Road}, 
postcode={100190}, 
state={Beijing},
country={China}}
            
\begin{abstract}
    Solving high-dimensional parabolic partial differential equations (PDEs) with deep learning methods is often computationally and memory intensive, primarily due to the need for automatic differentiation (AD) to compute large Hessian matrices in the PDE. In this work, we propose a deep random difference method (DRDM) that addresses these issues by approximating the convection-diffusion operator using only first-order differences and the solution by deep neural networks, thus avoiding Hessian and other derivative computations. The DRDM is implemented within a Galerkin framework to reduce sampling variance, and the solution space is explored using stochastic differential equations (SDEs) to capture the dynamics of the convection-diffusion operator. The approach is then extended to solve Hamilton-Jacobi-Bellman (HJB) equations, which recovers existing martingale deep learning methods for PDEs [{\it SIAM J. Sci. Comput.}, 47 (2025), pp. C795-C819], without using stochastic calculus. The proposed method offers two main advantages: it avoids the need to compute derivatives in PDEs and enables parallel computation of the loss function in both time and space. Moreover, a rigorous error estimate is proven for the quasi-linear parabolic equation, showing first-order accuracy in $h$, the time step used in the discretization of the SDE paths by the Euler-Maruyama scheme. 
    Numerical experiments demonstrate that the method can efficiently and accurately solve quasilinear parabolic PDEs and HJB equations in dimensions up to $10^5$ and $10^4$, respectively.
\end{abstract}




\end{frontmatter}



\section{Introduction}

In this work, we focus on deep learning methods for solving high-dimensional quasilinear parabolic partial differential equations (PDEs) of the form
\begin{equation}\label{eq_pde}
    \mathcal{D} v(t, x) = f\big(t, x, v(t, x)\big), \quad (t, x) \in [0, T] \times \R^d
\end{equation}
with a terminal condition $v(T, x) = g(x)$, $x \in \R^d$. 
Here $v: [0, T] \times \R^d \to \R$ is the unknown solution, and $\mathcal{D}$ is the nonlinear second-order evolutionary differential operator defined as
\begin{equation}\label{eq_defD}
    \mathcal{D} := \partial_t + \mu^\top\vbr{\big}{t, x, v(t, x)} \partial_x + \frac{1}{2} \mathrm{Tr}\rbr{\sigma \sigma^{\top}\vbr{\big}{t, x, v(t, x)} \partial_{xx}}
\end{equation}
where $\mu$ and $\sigma$ are given functions valued in $\R^d$ and $\R^{d \times q}$, respectively, and $\mathrm{Tr}\rbr{\,\cdot\,}$ denotes the trace of a matrix.
The PDE \eqref{eq_pde} encompasses a wide range of important models, including the Allen-Cahn equation, Black-Scholes equation, Kolmogorov backward equation, Burgers' equation, and various reaction-diffusion equations.
Beyond~\eqref{eq_pde}, our proposed method is further extended to Hamilton-Jacobi-Bellman (HJB) equations, which are often high-dimensional and play a central role in stochastic optimal control problems.
Those high-dimensional PDEs have applications in a wide range of areas, including finance, physics, and engineering. 
Traditional numerical methods for high-dimensional PDEs often suffer from the curse of dimensionality (CoD), where the computational cost increases exponentially with the dimensionality.  
Recent advances in deep learning have demonstrated significant potential for overcoming the CoD, yielding various deep learning-based numerical methods for high-dimensional PDEs.
Currently, the popular deep learning methods for solving PDEs can be broadly classified into two categories: 
\begin{itemize}
    \item \textbf{Direct approach}: neural networks are trained to learn the solutions of \eqref{eq_pde} directly, as exemplified by the Physics-Informed Neural Networks (PINNs) \cite{Raissi2019Physics,hu2024sdgd,Hu2024Hutchinson,He2023Learning,wang2022is,Gao2023Failure}, Deep Galerkin Method (DGM) \cite{Sirignano2018DGM,Al2022Extensions}, Deep Ritz \cite{E2018deep}, Weak Adversarial Networks (WANs) \cite{Zang2020Weak,Chen2023Friedrichs}, etc. 
    
    \item \textbf{SDE-based approach}: the PDE~\eqref{eq_pde} is reformulated into a stochastic differential equation (SDE) and then solved by deep learning. 
    Relevant works include DeepBSDE \cite{weinan2017deep,han2018solving,han2018convergence}, forward-backward stochastic neural networks~\cite{raissi2018forwardbackward,Zhang2022FBSDE}, deep splitting method \cite{beck2019deep}, deep backward schemes \cite{hure2020deep,germain2022Approximation},
    Actor-critic method \cite{Zhou2021Actor}, 
    deep neural networks algorithms \cite{Hure2021Deep,Bachouch2022Deep},
    diffusion loss \cite{Nusken2023Interpolating}, diffusion Monte Carlo-like approach \cite{Han2020Solving}, martingale neural networks \cite{cai2024socmartnet,cai2023deepmartnet,cai2024martnet}, etc.
\end{itemize}

Direct approaches, such as PINN and DGM, train a neural network $\hat{v}$ to approximate the solution $v$ of the PDE \eqref{eq_pde} by minimizing the loss function  
\begin{equation}\label{eq_loss}
    L_{\mathrm{pinn}}(\hat{v}) := \int_0^T \int_{\R^d} \abs{\mathcal{D} \hat{v}(t, x) - f(t, x, \hat{v}(t, x))}^2 p(t, x) \di x \di t,
\end{equation}
where $p: [0, T] \times \R^d \to \R$ is a weight function.
For clarity of exposition, we omit the penalty term for the terminal condition, or equivalently, assume that $v$ satisfies the terminal condition a priori.
For implementation, the PINN loss in \eqref{eq_loss} is typically minimized using stochastic gradient descent, where the integrals in \eqref{eq_loss} are approximated by mini-batch samples $(t_i, x_i) \in [0, T] \times \R^d$, with $p(t, x)$ serving as the sampling distribution.

A key advantage of \eqref{eq_loss} is that the computations for different sampled points $(t_i, x_i)$ are independent, making the method highly suitable for parallelization. 
However, these methods depend on automatic differentiation (AD) to evaluate the differential operator $\mathcal{D} \hat{v}$. 
The back-propagation method to compute the $d \times d$ Hessian matrix is both computationally and memory intensive, especially in high dimensions, creating a bottleneck which has motivated a series of recent studies. 
For instance, \cite{He2023Learning} proposes a Gaussian-smoothed neural network, enabling explicit computation of derivatives via Stein's identity. 
A series of works \cite{hu2024sdgd,Hu2024Hutchinson,shi2024stochastic} address this challenge by introducing techniques such as stochastic dimension gradient descent, Hutchinson trace estimation, and stochastic Taylor derivative estimators in the framework of PINNs.
In addition, \cite{Li2024computational,Li2024DOF} propose computing high-order derivatives via forward propagation, which significantly reduces computational cost and memory usage without introducing statistical errors.

SDE-based approaches are grounded in the nonlinear Feynman-Kac formula, i.e., the solution $v$ to the PDE \eqref{eq_pde} can be represented as 
\begin{equation}\label{eq_feynman}
    Y_t = v(t, X_t), \quad Z_t = \partial_x v(t, X_t) \sigma(t, X_t, Y_t), \quad t \in [0, T],
\end{equation}
where $t \mapsto (X_t, Y_t, Z_t)$ are stochastic processes given by the forward-backward SDEs (FBSDEs)  
on a probability space $(\Omega, \mathcal{F}, \mathbb{P})$ endowed with the natural filtration $\{\mathcal{F}_t\}_{t \in [0, T]}$ generated by a $q$-dimensional standard Brownian motion $B$:
\begin{equation}\label{eq_bsde}
\left\{\begin{aligned}
    X_t &= x_0 + \int_0^t \mu(s, X_s, Y_s) \di s + \int_0^t \sigma(s, X_s, Y_s) \di B_s,\\
    Y_t &= g(X_T) - \int_t^T f(s, X_s, Y_s) \di s - \int_t^T Z_s \di B_s,
\end{aligned}
\right.\quad t \in [0, T],
\end{equation}
where $\int_0^t \cdot \di B_s$ denotes the It\^o integral with respect to $B_s$.
To clarify the connection with \eqref{eq_pde}, the FBSDEs~\eqref{eq_bsde} take a special form where the generator $f$ is independent of the process $Z$.

By reformulating the problem into the FBSDEs~\eqref{eq_bsde},
SDE-based approaches avoid the need to compute $\mathcal{D} \hat{v}$ directly. 
However, this advantage comes at the cost of reduced parallelism 
along the time dimension compared to direct approaches. 
Specifically, 
(i) For quasilinear parabolic PDEs \eqref{eq_pde}, many methods such as those in \cite{weinan2017deep,han2018solving,han2018convergence,raissi2018forwardbackward,Zhang2022FBSDE,Nusken2023Interpolating,Han2020Solving,Zhou2021Actor} rely on simulating sample paths that depend on the unknown process $Y_t = v(t, X_t)$ being learned. 
Consequently, these sample paths need to be updated iteratively during training, and such updates are inherently sequential in time, resulting in increased computational cost.
(ii) The methods in \cite{beck2019deep,hure2020deep,germain2022Approximation,Hure2021Deep,Bachouch2022Deep} propose to solve the problems backward in time, which can 
reduce training difficulty, but limits the time parallelizability of the algorithm as well.

More recently, a line of hybrid approaches has been developed to combine the strengths of both direct and SDE-based approaches.
In \cite{Han2020derivative}, PDE derivatives are approximated via the Feynman--Kac formula, and Brownian walks are used to explore the domain. This retains the parallelism of direct approaches while avoiding the computation of Hessians.
By exploring the connection between the residuals of the FBSDEs~\eqref{eq_bsde} and the PDE~\eqref{eq_pde}, 
the shotgun method \cite{zhang2025shotgun} uses a coarse time step for spatial exploration and finer time steps for enforcing the FBSDE residual, thereby improving computational efficiency. 
The SOC-MartNet method \cite{cai2024socmartnet} and its derivative-free variant \cite{cai2024martnet} recast the PDE~\eqref{eq_pde} as a martingale problem, which enables path sampling to be conducted offline, thereby allowing for highly parallelized online training. 

The aforementioned SDE-based and hybrid approaches typically rely on It\^o calculus, which may create a challenge for those without a solid background in stochastic analysis. 
Thus, it is the objective of this paper to develop the martingale methods proposed in \cite{cai2024socmartnet,cai2023deepmartnet,cai2024martnet} from a new perspective that is closer to the PINN framework and requires minimal stochastic calculus.
In addition, this paper provides further insights into error estimates, sampling methods, 
and the connections between PINN methods and SDE-based approaches.

Specifically, our main contributions are as follows:
\begin{enumerate}
    \item We propose a random difference method (RDM) to approximate the operator $\mathcal{D}$ in \eqref{eq_defD}, thereby recasting the PDE~\eqref{eq_pde} into an RDM formulation. 
    This formulation bridges the gap between PINNs and some SDE-based approaches. 
    In particular, this formulation can be incorporated into the PINN framework to eliminate the need for AD to compute large Hessian matrices. 
    On the other hand, it is equivalent to the martingale condition underlying the martingale methods in \cite{cai2024socmartnet,cai2023deepmartnet,cai2024martnet}, 
    and is also equivalent to the implicit Euler scheme for the FBSDEs \eqref{eq_bsde} as proposed in \cite{Zhang2004numerical,zhao2006new}.
    
    \item We provide rigorous error estimates for the RDM formulation. 
    This result justifies the first-order convergence observed in the martingale methods in \cite{cai2024socmartnet,cai2023deepmartnet,cai2024martnet}. 
    This convergence rate is higher than that of conventional SDE-based approaches, which typically achieve at most half-order accuracy, as analyzed in \cite{han2018convergence,hure2020deep,germain2022Approximation} and verified numerically in \cite{Zhang2022FBSDE}.
    
    \item We further establish a connection between the error and the residual of the neural network for the PDE~\eqref{eq_pde} in the case of linear drift-diffusion equations. 
    This analysis reveals that an appropriate spatial sampling weight is given by the solution to the adjoint equation, namely, the Fokker-Planck equation.
    As a result, this provides a theoretical justification for SDE-based approaches, where the process generated by the operator $\mathcal{D}$ serves as the sample generator to explore the solution space. 
    We extend this sampling method to the quasilinear PDE~\eqref{eq_pde}, and numerical experiments demonstrate its effectiveness, particularly for convection-dominated problems, even in high dimensions. 
    
    \item The RDM approximation introduces an additional expectation inside the loss function. 
    Directly estimating this expectation via Monte Carlo methods is inefficient, resulting in a convergence rate slower than the standard $1/\sqrt{n}$ rate, where $n$ is the number of samples.
    To address this issue, we reformulate the RDM using a Galerkin method as a minimax problem. 
    This reformulation enables efficient sampling with standard convergence rates, and can be implemented with adversarial learning, an approach that has proven more suitable for HJB equations in high dimension \cite{wang2022is}.
\end{enumerate}

With the above contributions, a deep random difference method (DRDM) is developed to solve the PDE \eqref{eq_pde}.
The DRDM has several advantages. 
First, although it is equivalent to the derivative-free martingale method in \cite{cai2024martnet}, its derivation is directly based on PDE theory rather than It\^o's stochastic calculus. 
The RDM formulation requires only Taylor expansions and basic properties of random variables (though the fundamental It\^o formula in stochastic calculus also relies on a Taylor expansion to the second order and the introduction of It\^o integral). 
The sampling method is grounded on the Fokker-Planck equation. 
Second, it combines the strengths of both direct and SDE-based approaches for PDEs: it avoids the need for computing derivatives in PDEs and enables parallel computation of the loss function in both time and space.
This significantly improves computational efficiency.
Numerical experiments demonstrate that the proposed method can efficiently and accurately solve quasilinear parabolic PDEs and HJB equations in dimensions up to $10^5$ and $10^4$, respectively.

The remainder of this paper is organized as follows. 
\Cref{sec_meth} introduces the RDM for differential operators, its application to the PDE \eqref{eq_pde}, and its connections to SDE-based approaches. 
\Cref{sec_rdo_pde} presents the deep learning implementation of the DRDM for solving \eqref{eq_pde} and extends the approach to HJB equations. 
The appendix contains theoretical error estimates, a variance analysis for mini-batch sampling, detailed parameter settings for the numerical experiments, and an extension of the RDM to quadratic gradient terms.

\section{Random difference method}\label{sec_meth}

In this section, we first present the random difference approximation to the differential operator defined in \eqref{eq_defD}. 
This technique forms a central component of our RDM for solving PDEs.
We then discuss the spatial sampling method for the RDM, and clarify its connections to both the martingale methods proposed in \cite{cai2024socmartnet,cai2023deepmartnet,cai2024martnet} and the implicit Euler schemes \cite{Zhang2004numerical,zhao2006new} for FBSDEs \eqref{eq_bsde}.

\subsection{Basic expansion}\label{sec_basic}

We begin by considering the Laplacian operator $\Delta = \sum_{i=1}^q \partial_{z_i}^2$ acting on a function $z \mapsto F(z) \in \R$ for $z \in \R^q$, 
where all the fourth-order partial derivatives of $F$ are assumed to be bounded on $\R^q$.
Let $\xi = \br{\xi_1, \xi_2, \cdots, \xi_q}^{\top}$ be a $\R^q$-valued random variable, and let $h > 0$ be a small step size. 
Applying the Taylor expansion of $F(z)$ at $z = 0$, substituting $z = \sqrt{h} \xi$, and taking the expectation with respect to $\xi$, we obtain
\begin{equation}\label{eq_taylor}
    \E{F(\sqrt{h} \xi)} = F(0) + \sum_{k=1}^3 \frac{h^{k/2}}{k!}\E{\vbr{\big}{\xi^{\top} \partial_z}^k F(0)} + \frac{h^2}{4!}\E{\vbr{\big}{\xi^{\top} \partial_z}^4 F(c\sqrt{h} \xi)},
\end{equation}
where $\xi^{\top} \partial_z := \sum_{i=1}^q \xi_i \partial_{z_i}$, and $c \in [0, 1]$ depending on $\sqrt{h} \xi$.
Assume the components $\xi_1, \xi_2, \cdots, \xi_q$ of $\xi$ satisfy:
\begin{equation}\label{eq_condi_xi}
    \E{\xi_i} = 0, \quad \E{\xi_i \xi_j} = \delta_{ij}, \quad \E{\xi_i \xi_j \xi_k} = 0, \quad \E{\abs{\xi_i \xi_j \xi_k \xi_l}} < \infty
\end{equation}
for $i, j, k, l = 1, 2, \cdots, q$, 
where $\delta_{ij}$ denotes the Kronecker delta. 
Under the conditions in \eqref{eq_condi_xi}, the first- and third-order derivative terms in \eqref{eq_taylor} vanish, yielding
\begin{equation}\label{eq_rdo_expa}
    \E{F(\sqrt{h} \xi)} = F(0) + \frac{h}{2} \sum_{i=1}^q \partial_{z_i}^2 F(0) + \mO\br{h^2}. 
\end{equation}
The above expansion serves as the basis for the RDM to be introduced in the next subsection.

\begin{remark}
    Typical examples of $\xi$ that satisfy \eqref{eq_condi_xi} include: 
    \begin{itemize}
        \item $\xi \sim \mathrm{N}(0, I_q)$, where $I_q$ is the $q$-dimensional identity matrix;
        \item $\xi$ has mutually independent components, each satisfying
        \begin{equation}
            \Pra{\xi_i = c} = \Pra{\xi_i = -c} = \frac{1}{2 c^2}, \quad \Pra{\xi_i = 0} = 1 - \frac{1}{c^2}
        \end{equation}  
        for some constant $c \geq 1$.  
    \end{itemize}
\end{remark}

\subsection{RDM for convection-diffusion operator}\label{sec_condiff}

Next, we consider the convection-diffusion operator $\mathcal{D}$ as defined in \eqref{eq_defD}.
For a fixed point $(t, x) \in [0, T] \times \R^d$, define the function
\begin{equation}\label{eq_defVsz}
    V(s, z) := v(t + s, x + \mu s + \sigma z), \quad (s, z) \in [0, T - t] \times \R^q, 
\end{equation}
where, for clarity, we suppress the dependence of $V$, $\mu$, and $\sigma$ on $(t, x, v)$.
Applying the expansion \eqref{eq_rdo_expa} with $F$ replaced by $z \mapsto V(h, z)$, we have
\begin{equation}\label{eq_expandV}
    \E{V(h, \sqrt{h}\xi)} = V(h, 0) + \frac{h}{2} \sum_{i=1}^q \frac{\partial^2 V}{\partial z_i^2} (h, 0) + \mO\br{h^2}. 
\end{equation}
By the definition in \eqref{eq_defVsz}, the Taylor expansion and the chain rule imply
\begin{align}
    V(h, 0) &= V(0, 0) + \br{\partial_t + \mu^{\top} \partial_x} v(t, x) h + \mO(h^2), \label{eq1_expandV}\\
    \sum_{i=1}^q \frac{\partial^2 V}{\partial z_i^2} (h, 0) &= \sum_{i=1}^q \frac{\partial^2 V}{\partial z_i^2}(0, 0) + \mO(h) = \mathrm{Tr}\rbr{\sigma \sigma^{\top} \partial_{xx} v(t, x)} + \mO(h). \label{eq2_expandV}
\end{align}
Inserting \eqref{eq1_expandV} and \eqref{eq2_expandV} into \eqref{eq_expandV}, we obtain 
\begin{equation}\label{eq_expandV2}
    \E{V(h, \sqrt{h}\xi)} = V(0, 0) + h \vbr{\big}{\partial_t + \mu^{\top} \partial_x + \frac{1}{2} \mathrm{Tr}\rbr{\sigma \sigma^{\top} \partial_{xx}}} v(t, x) + \mO(h^2). 
\end{equation}
By substituting \eqref{eq_defD} and \eqref{eq_defVsz} into the above equation, we arrive at the random difference approximation for $\mathcal{D}$ as 
\begin{equation}\label{eq_rdo}
    \mathcal{D}_h v(t, x) := \E{\frac{v(t + h, x + \mu h + \sigma \sqrt{h} \xi) - v(t, x)}{h}} = \mathcal{D}v(t, x) + \mO(h). 
\end{equation}
A rigorous derivation of \eqref{eq_rdo}, along with an explicit upper bound for the remainder term $\mO(h)$, is provided in \Cref{coro_consist} in \cref{sec_theory}.
An extension of the random difference approximation for the quadratic gradient term $\abs{\partial_x v}^2$ is presented in \cref{sec_rdm_quad}. 

A notable feature of \eqref{eq_rdo} is that it enables the approximation of second-order derivatives using only first-order differences, thereby avoiding explicit computation of the Hessian.

\begin{remark}[Regularity condition]\label{rmk_reg}
The approximation \eqref{eq_rdo} can also be derived from It\^o's formula \cite[Theorem 2.3.3]{Zhang2017Backward}.
This idea was exploited in \cite{Han2020derivative,zhang2025shotgun} to develop deep learning methods for PDEs.
Although the It\^o formula itself only requires $v \in C^{1,2}$, i.e., once continuously differentiable in $t$ and twice in $x$, 
we still need a stronger regularity condition $v \in C^{2,4}$ to ensure the remainder in \eqref{eq_rdo} achieves the order $\mO(h)$.
This requirement is mild for parabolic equations:
under some strong parabolicity conditions on $\sigma$ \cite[p.~409, (8.2)]{Taylor2023PDEIII},
the solution $v(t, x)$ of \eqref{eq_pde} is infinitely differentiable for $(t, x) \in [0, T) \times \R^d$; 
see \cite[p.~412, Proposition 8.3]{Taylor2023PDEIII}.
\end{remark}

\begin{remark}[Comparison to randomized smoothing PINNs (RS-PINNs)]\label{rmk_rspinn}
    RS-PINNs \cite{He2023Learning,Hu2025Bias} share the same motivation as the RDM in avoiding the use of AD for PDE derivatives, but employ a different approach.
    Specifically, RS-PINNs use randomized smoothing neural networks to approximate the PDE solution, i.e., $v(x) \approx v_{\theta}(x) := \E{\psi_{\theta}(x + \sqrt{h}\xi)}$ for $x \in \R^d$, where $\xi$ is a $d$-dimensional standard Gaussian random vector and $\psi_{\theta}$ is a neural network parameterized by $\theta$. 
    This smoothing allows the computation of derivatives via Stein's identity \cite{Stein1981Estimation}, namely,
    \begin{equation}\label{eq_stein}
        \partial_x v_{\theta}(x) = \E{\frac{\xi}{\sqrt{h}} \psi_{\theta}(x + \sqrt{h} \xi)}, \;\; \partial_{xx} v_{\theta}(x) = \E{\frac{\xi \xi^{\top} - I_d}{h} \psi_{\theta}(x + \sqrt{h} \xi)}.    
    \end{equation}
    Compared to the RDM in \eqref{eq_rdo}, RS-PINNs are more like a central difference method, while the RDM is an upwind one similar to the material derivative in fluid dynamics. 
    This distinction makes the RDM more suitable for convection-dominated problems, as demonstrated by our numerical experiments in \cref{sec_cde,sec_hjb,sec_ablation}.
\end{remark}

\subsection{RDM for solving quasilinear parabolic PDEs}\label{sec_rdm}

By approximating $\mathcal{D}$ with $\mathcal{D}_h$ defined in \eqref{eq_rdo}, the PDE in \eqref{eq_pde} can be reformulated as
\begin{equation}\label{eq_ER0}
    \E{R(t, x, \xi; v)} = \mO(h), \quad (t, x) \in [0, T-h] \times \R^d,
\end{equation}
where $R(t, x, \xi; v)$ is the residual function defined by
\begin{align}
    R(t, x, \xi; v) &:= \frac{v(t + h, x + \xi_h) - v(t, x)}{h} - f\vbr{\big}{t, x, v(t, x)}, \label{eq_defR}\\
    \xi_h &:= \mu\vbr{\big}{t, x, v(t, x)} h + \sigma\vbr{\big}{t, x, v(t, x)} \sqrt{h} \xi. \label{eq_def_xih}
\end{align}
Furthermore, analogous to the PINN loss defined in \eqref{eq_loss}, 
the RDM loss for the formulation~\eqref{eq_ER0} is given by
\begin{equation}\label{eq_Lrdm}
    L_{\mathrm{rdm}}(\hat{v}) := \int_0^{T-h} \int_{\R^d} \big|\E{R(t, x, \xi; \hat{v})}\big|^2 p(t, x) \di x \di t, 
\end{equation}
where the expectation is taken with respect to the random variable $\xi$; $p: [0, T] \times \R^d \to \R$ is a weight function specified later in \cref{sec_pathsamp}; $\hat{v}: [0, T] \times \R^d \to \R$ is a neural network approximating the solution $v$ of \eqref{eq_pde}.
In \Cref{thm_convergence} in \cref{sec_theory}, we rigorously prove that the solution minimizing a discrete version of the RDM loss in \eqref{eq_Lrdm} achieves a first-order convergence rate in the time step size $h$.

The RDM loss in \eqref{eq_Lrdm} avoids explicit computation of spatial partial derivatives, including the Hessian; 
however, an expectation in \eqref{eq_Lrdm} is needed for each $(t, x)$, which may lead to inefficient mini-batch estimation. 
Actually, suppose $\hat{L}_{\mathrm{rdm}}(\hat{v})$ is an unbiased mini-batch estimator of the RDM loss in \eqref{eq_Lrdm}, using $M$ sampled points $(t_i, x_i)$ to approximate the integral over $[0, T-h] \times \R^d$, and $2 K$ samples of $\xi$ to estimate the expectation in \eqref{eq_Lrdm} at each $(t_i, x_i)$; see \eqref{eq_Lrdm_mini} in \cref{sec_variance}. 
Thus $\hat{L}_{\mathrm{rdm}}(\hat{v})$ requires at least $2 M K$ evaluations of the residual $R(t, x, \xi; \hat{v})$. 
By \eqref{eq_VarLrdm_final} in \cref{sec_variance}, such an unbiased mini-batch estimator has a variance of order
\begin{equation}\label{eq_varLrdm}
    \mathrm{Var}\rbr{\hat{L}_{\mathrm{rdm}}(\hat{v})} = \mO\br{\frac{1}{M} \vbr{\Big}{1 + \frac{1}{K^2}}}, \quad M, K \to \infty,
\end{equation}
which is slower than the standard Monte Carlo convergence rate of $1/(MK)$. 
In \cref{sec_rdo_pde}, we address this issue by reformulating \eqref{eq_Lrdm} within a Galerkin framework, which reduces the variance of the mini-batch estimator to order $\mathcal{O}\vbr{\big}{1/(MK)}$.

\subsection{Spatial sampling rule based on error propagation analysis}\label{sec_pathsamp}

Now, we analyze how the residual of the PDE \eqref{eq_pde} propagates and affects the error of the numerical solution. 
This analysis provides guidance for designing the weight function $p(t, x)$ in \eqref{eq_Lrdm}.

\subsubsection{Error equation}

Let $\hat{v}$ be an approximation to the solution $v$ of \eqref{eq_pde}.
We define the error function $\epsilon$ and the residual function $\hat{r}$ as 
\begin{equation}\label{eq_def_rtx}
    \epsilon(t, x) := \hat{v}(t, x) - v(t, x), \quad \hat{r}(t, x) := (\partial_t + \mathcal{L}) \hat{v}(t, x) - f(t, x, v)
\end{equation}
for $(t, x) \in [0, T] \times \R^d$, where $\partial_t + \mathcal{L}$ is a decomposition of the operator $\mathcal{D}$ defined in \eqref{eq_defD}:
\begin{equation}\label{eq_decompD}
    \mathcal{D} = \partial_t + \mathcal{L}, \quad \mathcal{L} := \mu^{\top}(t, x, v) \partial_x + \frac{1}{2} \mathrm{Tr}[\sigma \sigma^{\top} (t, x, v) \partial_{xx}].
\end{equation} 
We emphasize that $\mathcal{L}$ and $f$ in \eqref{eq_def_rtx} are evaluated at the exact solution $v$, rather than at $\hat{v}$. 
Thus, $\hat{r}(t, x)$ is a linearized residual, enabling us to establish an explicit relationship between $\epsilon$ and $\hat{r}$.

Combining \eqref{eq_pde} and \eqref{eq_def_rtx}, we can obtain the error equation
\begin{equation}\label{eq_error_eq}
    (\partial_t + \mathcal{L}) \epsilon(t, x) = \hat{r}(t, x), \quad (t, x) \in [0, T) \times \R^d. 
\end{equation}
Let $p(t, x)$ be a weight function defined on $[0, T] \times \R^d$.
Taking the inner product of $\epsilon(t,x)$ with $p(t, x)$ over $x \in \R^d$ and then differentiating with respect to $t$, we obtain the integral form of the error equation:
\begin{equation}\label{eq_dt_vp}
\begin{aligned}
    \partial_t \int_{\R^d} \epsilon p \di x &= \int_{\R^d} (\partial_t \epsilon) p \di x + \int_{\R^d} \epsilon (\partial_t p) \di x \\
    &= \int_{\R^d} \hat{r} p \di x - \int_{\R^d} (\mathcal{L} \epsilon) p \di x + \int_{\R^d} \epsilon \partial_t p \di x, 
\end{aligned}
\end{equation}
where the second line follows from \eqref{eq_error_eq}. 

\subsubsection{Error propagation from residuals}

Now, we shall simplify the right-hand side of \eqref{eq_dt_vp} to obtain an ordinary differential equation with an explicit solution for $\int \epsilon p \,\di x$. 
Let $\mathcal{L}^*$ denote the $L^2$-adjoint of $\mathcal{L}$ on $\R^d$, so that
\begin{equation}\label{eq_Ladj}
    \int_{\R^d} (\mathcal{L} \epsilon) p \di x = \int_{\R^d} \epsilon (\mathcal{L}^* p) \di x. 
\end{equation}
Inserting \eqref{eq_Ladj} into \eqref{eq_dt_vp} yields
\begin{equation}\label{eq_error_eq_lin}
    \partial_t \int_{\R^d} \epsilon p \di x = \int_{\R^d} \hat{r} p \di x + \int_{\R^d} \epsilon (\partial_t p - \mathcal{L}^* p) \di x. 
\end{equation}
In particular, let $p$ satisfy the adjoint equation:
\begin{equation}\label{eq_adjpde}
\left\{\begin{aligned}
    &\fbr{\partial_t - \mathcal{L}^*} p = 0, \quad (t, x) \in (0, T] \times \R^d,\\
    &p(0, x) = \delta_{x_0}(x), \quad x \in \R^d,
\end{aligned}\right.
\end{equation}
where $\delta_{x_0}$ denotes the Dirac delta function centered at a spatial point $x_0 \in \R^d$ of interest.
By substituting \eqref{eq_adjpde} into the right-hand side of \eqref{eq_dt_vp}, the second integral term is eliminated. 
Integrating both sides with respect to $t$ over $[0, T]$, we obtain the error representation at $(0, x_0)$:
\begin{equation}\label{eq_v0x0}
\begin{aligned}
    \epsilon(0, x_0) &= \int_{\R^d} (\epsilon p)(0, x) \di x \\
    &= \int_{\R^d} \epsilon(T, x) p(T, x) \di x - \int_0^T \int_{\R^d} \hat{r}(s, x) p(s, x) \di x \di s. 
\end{aligned}
\end{equation}

The representation \eqref{eq_v0x0} shows that the error at $(0, x_0)$ is determined by the terminal error, weighted by $p(T, x)$, and the accumulation of the residual $\hat{r}$ over time and space, weighted by $p(s, x)$. 
Therefore, $p(s, x)$ quantifies the sensitivity of the error at $(0, x_0)$ to the residual at $(s, x)$. 
A natural strategy is to sample $(t, x) \in [0, T] \times \R^d$ from a distribution whose density is proportional to $p$.
However, doing so directly is costly, as it entails solving the adjoint problem~\eqref{eq_adjpde} to obtain $p$.
The next subsection avoids this expense via an SDE-based sampling approach.

\subsubsection{Sampling method}\label{sec_pathsamp2}

The adjoint equation \eqref{eq_adjpde} coincides with the Fokker-Planck equation, whose solution $x \mapsto p(t, x)$ represents the PDF of $X_t$ governed by the forward SDE in \eqref{eq_bsde} \cite[Proposition 3.3]{pavliotis2014stochastic}.
Therefore, spatial samples can be efficiently generated by simulating the forward SDE in \eqref{eq_bsde} using the weak Euler-Maruyama scheme \cite[section 14.1]{Kloeden1992Numerical}. 
In practical implementation, the unknown solution $v$ in \eqref{eq_bsde} is replaced by a pilot approximation $\hat{v}$ obtained during neural network training.

Specifically, we introduce a time partition on $[0, T]$:
\begin{equation}\label{eq_timepart}
    t_n = n h, \quad h = T/N, \quad n = 0, 1, \cdots, N.
\end{equation}
The spatial samples $\bbr{X_n^m}_{m=1}^M \subset \R^d$ at time $t_n$ are generated by
\begin{equation}\label{eq_pathX}
    X_{n+1}^m = X_n^m + \mu\br{t_n, X_n^m, \hat{v}(t_n, X_n^m)}\,h + \sigma\br{t_n, X_n^m, \hat{v}(t_n, X_n^m)} \,\sqrt{h}\, \xi_{n+1}^{m},
\end{equation}
for $n = 0, 1, \cdots, N-1$, where for any $n$ and $m$, the random variables $\xi_n^m$ are i.i.d. samples of $\xi$ appearing in the random difference approximation in \eqref{eq_rdo}. 

\begin{remark}[Connection to SDE-based approaches]
    The sampling method~\eqref{eq_pathX} is widely adopted in SDE-based approaches for exploring the space $\R^d$; 
    see, e.g., \cite{weinan2017deep,han2018solving,raissi2018forwardbackward,Zhang2022FBSDE,beck2019deep,hure2020deep,germain2022Approximation,Zhou2021Actor,Hure2021Deep,Bachouch2022Deep,Nusken2023Interpolating}.
    Our analysis offers a theoretical justification for this sampling method from the perspective of error propagation.
\end{remark}

\subsection{Connection to martingale methods}\label{sec_con_mart}

The formulation \eqref{eq_ER0} is equivalent to a discrete version of the martingale formulation proposed in \cite{cai2024socmartnet,cai2023deepmartnet,cai2024martnet}. 
To see this connection, we recall the martingale condition for the function $v(t, x)$ being a solution to the PDE \eqref{eq_pde}:
\begin{equation}\label{eq_mart_cond}
    \Etx{\M_{t+h} - \M_{t}} = 0, \quad (t, x) \in [0, T-h] \times \R^d,
\end{equation}
where $\Etx{\,\cdot\,} := \E{\,\cdot\, \big\vert X_t = x}$,
and the process $\M$ (which should be a martingale) is defined by \cite{cai2024socmartnet,cai2023deepmartnet,cai2024martnet}
\begin{equation}\label{eq_defMt}
    \M_t := v(t, X_t) + \int_0^t f(s, X_s, v(s, X_s)) \di s,
\end{equation}
with $X_s$ given by the forward SDE in \eqref{eq_bsde}.

To discretize \eqref{eq_mart_cond}, we first substitute $\M_{t+h} - \M_{t}$ using \eqref{eq_defMt}, which yields
\begin{equation}\label{eq_mart_cond2}
    \Etx{v(t+h, X_{t+h})} - v(t, x) + \int_t^{t+h} \Etx{f(s, X_s, v(s, X_s))} \di s = 0
\end{equation}
for $(t, x) \in [0, T-h] \times \R^d$, where we have used the fact that $\Etx{v(t, X_t)} = v(t, x)$. 
Then the conditional expectation and the integral in \eqref{eq_mart_cond2} can be approximated by the weak Euler-Maruyama scheme \cite[section 14.1]{Kloeden1992Numerical} and the left-rectangle rule, respectively, yielding
\begin{align}
    &\Etx{v(t+h, X_{t+h})} = \E{v(t+h, x + \xi_h)} + \mO(h^2), \label{eq_EulerApprox}\\
    \int_t^{t+h} &\Etx{f(s, X_s, v(s, X_s))} \di s = h f(t, x, v(t, x)) + \mO(h^2), \label{eq_left_rect}
\end{align}
where $\xi_h$ is given by \eqref{eq_def_xih}, and the remainder terms hold under sufficient regularity conditions on $f$ and $v$.
Substituting these approximations into \eqref{eq_mart_cond2} and omitting the remainder terms, we arrive at the discretized martingale condition: for $(t, x) \in [0, T-h] \times \R^d$,
\begin{equation}\label{eq0_mart_cond3}
    \E{v(t+h, x + \xi_h)} - v(t, x) - h f(t, x, v(t, x)) = \mO(h^2). 
\end{equation}

Using the residual function defined in \eqref{eq_defR}, we can rewrite the left-hand side of \eqref{eq0_mart_cond3}. Consequently, \eqref{eq0_mart_cond3} becomes
\begin{equation}\label{eq_mart_cond3}
    h\,\E{R(t, x, \xi; v)} = \mO(h^2), \quad (t, x) \in [0, T-h] \times \R^d. 
\end{equation}
This clearly shows that the RDM formulation \eqref{eq_ER0} is equivalent to the discretized martingale condition \eqref{eq_mart_cond}.
In \cref{sec_rdo_pde}, we will show that the RDM formulation \eqref{eq_ER0} can be used to derive the derivative-free martingale method proposed in \cite{cai2024martnet}.

\subsection{Connection to the Euler scheme for FBSDEs}\label{sec_fbsde}

The RDM formulation~\eqref{eq_ER0} is also equivalent to the implicit Euler scheme for the FBSDEs~\eqref{eq_bsde}, which is a special case of general FBSDEs where the generator $f$ does not depend on the process $Z$. 
The explicit Euler scheme for general FBSDEs was originally proposed in \cite{Zhang2004numerical} and further generalized in \cite{zhao2006new,Zhao2012A}. 
To clarify the connection, we briefly outline the derivation of the aforementioned Euler scheme, following the approach proposed in \cite[sections 3.1, 3.2.1]{zhao2006new}.

By the well-posedness of FBSDEs, the backward SDE in \eqref{eq_bsde} can be restricted to the time grid defined in \eqref{eq_timepart}, and written as
\begin{equation*}
    Y_{t_n} = Y_{t_{n+1}} - \int_{t_n}^{t_{n+1}} f(s, X_s, Y_s) \di s - \int_{t_n}^{t_{n+1}} Z_s \di B_s
\end{equation*}
for $n = N-1, N-2, \cdots, 0$. 
Taking $\Etnx{\,\cdot\,} := \E{\,\cdot\, \big\vert X_{t_n} = x}$ on both sides of the above equation, the It\^o integral vanishes, yielding
\begin{equation*}
    \Etnx{Y_{t_n}} = \Etnx{Y_{t_{n+1}}} - \int_{t_n}^{t_{n+1}} \Etnx{f(s, X_s, Y_s)} \di s, \quad x \in \R^d. 
\end{equation*}
Substituting $Y$ using the first Feynman-Kac formula in \eqref{eq_feynman}, the above equation can be rewritten as
\begin{equation*}
    v(t_n, x) = \Etnx{v(t_{n+1}, X_{t_{n+1}})} - \int_{t_n}^{t_{n+1}} \Etnx{f(s, X_s, v(s, X_s))} \di s. 
\end{equation*}
Applying the approximations in \eqref{eq_EulerApprox} and \eqref{eq_left_rect} to the above equation, we obtain the implicit Euler scheme for the FBSDEs~\eqref{eq_bsde}, which reads
\begin{equation}\label{eq_EulerFBSDE}
    v(t_n, x) = \E{v(t_{n+1}, x + \xi_h)} - h f(t_n, x, v(t_n, x)), \quad x \in \R^d,
\end{equation}
where the index $n$ runs from $N-1$ down to $0$, and $\xi_h$ is defined in \eqref{eq_def_xih}.
This scheme is a special case of the $\theta$-schemes proposed in \cite{zhao2006new}; see \Cref{rmk_thetasch} for more details.

The scheme~\eqref{eq_EulerFBSDE} can be recovered from the RDM formulation~\eqref{eq_ER0} by substituting \eqref{eq_defR} into \eqref{eq_ER0} and rearranging terms to move $v(t_n, x)$ to the left-hand side.
Thus the RDM formulation~\eqref{eq_ER0} is equivalent to the implicit Euler scheme~\eqref{eq_EulerFBSDE} for the FBSDEs~\eqref{eq_bsde}.

\begin{remark}\label{rmk_thetasch}
    In some literature \cite{Chassagneux2014Runge,Chassagneux2014Linear} on numerical methods for FBSDEs, the equation \eqref{eq_EulerFBSDE} is referred to as the $Y$-part of a numerical scheme for FBSDEs, as it gives the relationship between $Y_{t_n}$ and $X_{t_n}$. 
    In more general settings, where the generator $f$ in \eqref{eq_bsde} also depends on the process $Z$, 
    \eqref{eq_EulerFBSDE} is modified as follows:
    \begin{equation*}
        v(t_n, x) = \E{v(t_{n+1}, x + \xi_h)} - h f\big(t_n, x, v(t_n, x), z(t_n, x)\big), \quad x \in \R^d. 
    \end{equation*}
    Here, $z(t_n, x)$ denotes a numerical approximation to $\Etnx{Z_{t_n}}$. 
    In this case, a second equation, referred to as the $Z$-part of the scheme, is required to determine the value of $z(t_n, x)$. 
    A simple example of such a $Z$-part is given as a special case of \cite[(3.12)]{zhao2006new}:
    \begin{equation*}
        z(t_n, x) = \frac{1}{h} \E{v(t_{n+1}, x + \xi_h) \sqrt{h} \xi^{\top}}, \quad x \in \R^d
    \end{equation*}
    with $(\xi_h, \xi)$ constrained by \eqref{eq_def_xih}. 
    For more details about numerical schemes for FBSDEs, see \cite{zhao2010stable,Zhao2012A,zhao2014new,Yang2020unified,Wang2022Sinc} and the references therein.
\end{remark}

In summary, the RDM formulation~\eqref{eq_ER0} provides a unified perspective that recovers both the discrete martingale formulation proposed in \cite{cai2024socmartnet,cai2023deepmartnet,cai2024martnet}, 
and the implicit Euler scheme \cite{Zhang2004numerical,zhao2006new} for the specific FBSDEs~\eqref{eq_bsde} with the generator $f$ independent of $Z$. 
Notably, this connection is established without relying on stochastic calculus.

\section{Deep neural network with RDM (DRDM)}\label{sec_rdo_pde}

In this section, we introduce the deep learning implementation for learning the solution $v$ from the RDM formulation~\eqref{eq_ER0}, and subsequently extend the method to HJB equations.
Most of the content in this section runs parallel to the development in \cite{cai2024martnet} for the martingale method based on \eqref{eq_mart_cond}.

A straightforward approach is to train neural networks by minimizing the RDM loss in \eqref{eq_Lrdm} using stochastic gradient descent (SGD) to estimate the integral and the expectation. 
However, this approach has two main drawbacks: (i) as discussed in \cref{sec_rdm}, the RDM loss in \eqref{eq_Lrdm} is inefficient for mini-batch estimation due to its relatively high variance; (ii) the RDM loss is essentially an $L^2$ loss over the time-space domain, which may not guarantee convergence for high-dimensional HJB equations, as observed in \cite[Theorem 4.4]{wang2022is}.

To address these issues, we introduce a Galerkin method that reformulates the RDM loss into a minimax problem, which is then solved via adversarial learning. 
This approach enables more efficient mini-batch estimation with reduced variance, and avoids using $L^2$ loss for high-dimensional HJB equations.

\subsection{RDM with a Galerkin method and reduced variance}\label{sec_galerkin}

We first reformulate~\eqref{eq_ER0} into a weak formulation:
\begin{equation}\label{eq_weakform}
    \int_0^{T-h} \int_{\R^d} \rho(t, x) \E{R(t, x, \xi; v)} p(t, x) \di x \di t = \mO(h), \quad \forall \rho \in \mathcal{T},
\end{equation}
where $\mathcal{T}$ denotes the set of test functions, typically defined as:
\begin{equation}\label{eq_defrho}
    \mathcal{T} := \bbr{\rho: [0, T] \times \R^d \to \R^r \;\Big\vert\; \rho \text{ is smooth and bounded}},
\end{equation}
and $p(t, x)$ is a weight function defined on $[0, T] \times \R^d$.
It is worth noting that the test function $\rho$ is chosen to be vector-valued, which enhances the stability of the adversarial training process described in the next subsection.

For notational convenience, at each time $t \in [0, T]$, we introduce a random variable $X_t$, which is independent of $\xi$ in \eqref{eq_weakform} and has a PDF given by $x \mapsto p(t, x)$. 
Then,
\begin{equation}\label{eq_jointpdf}
    P_{(X_t, \xi)}(x, z) = p(t, x) P_{\xi}(z), \quad (t, x, z) \in [0, T] \times \R^d \times \R^q,
\end{equation}
where $P_{(X_t, \xi)}$ is the joint PDF of $(X_t, \xi)$, and $P_{\xi}$ is the marginal PDF of $\xi$.
Using \eqref{eq_jointpdf}, the left side of \eqref{eq_weakform} can be rewritten as
\begin{align}
    &\int_0^{T-h} \int_{\R^d} \rho(t, x) \E{R(t, x, \xi; v)} p(t, x) \di x \di t \notag\\
    =\;& \int_0^{T-h} \int_{\R^d} \int_{\R^q} \rho(t, x) R(t, x, z; v) P_{\xi}(z)  p(t, x) \di z \di x \di t \notag\\
    =\;& \int_0^{T-h} \int_{\R^d} \int_{\R^q} \rho(t, x) R(t, x, z; v) P_{(X_t, \xi)}(x, z) \di z \di x \di t \notag\\
    =\;& \int_0^{T-h} \E{\rho(t, X_t) R(t, X_t, \xi; v)} \di t, \label{eq_intERdt} 
\end{align}
where the expectation in the last line is taken with respect to $(X_t, \xi)$.
Substituting \eqref{eq_intERdt} into \eqref{eq_weakform}, we obtain the following minimax problem for approximating $v$ by a neural network $\hat v$:
\begin{equation}\label{eq_minmax}
    \min_{\hat{v} \in \mathcal{V}} \max_{\rho \in \mathcal{T}} \abs{L(\hat{v}, \rho)}^2, \quad L(\hat{v}, \rho) := \int_0^{T-h} \E{\rho(t, X_t) R(t, X_t, \xi; \hat{v})} \di t,
\end{equation}
where $\mathcal{V}$ represents the set of candidate functions $v$ satisfying the terminal condition, i.e., 
\begin{equation}\label{eq1_defV}
    \mathcal{V} := \bbr{\hat{v}: [0, T] \times \R^d \to \R \;\;\big\vert\; \hat{v} \in C^{2, 4}, \; \hat{v}(T, x) = g(x), \;\; x \in \R^d}
\end{equation}
with the regularity condition ensuring the validity of the RDM method in \eqref{eq_rdo}; see \Cref{rmk_reg}.
We note that the minimax loss $L(\hat{v}, \rho)$ is $\R^r$-valued, owing to the vector-valued test function $\rho$ defined in~\eqref{eq_defrho}.
The minimax formulation~\eqref{eq_minmax} is equivalent to the derivative-free martingale method proposed in \cite{cai2024martnet}.

Compared to the original loss function in \eqref{eq_Lrdm}, the minimax loss in \cref{eq_minmax} enables more efficient mini-batch estimation, as it involves only a single expectation over the whole time-space domain. 
In particular, if this expectation is approximated using $2MK$ samples of $(t, X_t, \xi)$ (see \eqref{eq_hatL} in \cref{sec_variance}), the resulting unbiased mini-batch estimator $|\hat{L}(\hat{v}, \rho)|^2$ requires $2MK$ evaluations of $R(t, x, \xi; \hat{v})$, matching the computational cost of the mini-batch RDM loss in \eqref{eq_Lrdm_mini}.
However, the variance of $|\hat{L}(\hat{v}, \rho)|^2$ is of order
\begin{equation}\label{eq_VarhatL2}
    \mathrm{Var}\rbr{\big|\hat{L}(\hat{v}, \rho)\big|^2} = \mO\br{\frac{1}{M K}}, 
\end{equation}
which is superior to that in \eqref{eq_varLrdm}. For a detailed derivation of these variances, see \cref{sec_variance}.

\subsection{Mini-batch estimation of the loss function}\label{sec_mini_batch}

Let $(X_n^m, \xi_n^m)$ for $n = 0, 1, \cdots, N$ and $m = 1, 2, \cdots, M$ be i.i.d. samples of $(X_n, \xi)$ generated by \eqref{eq_condi_xi} and \eqref{eq_pathX}.  
The loss function in \eqref{eq_minmax} can be estimated in an unbiased manner using its mini-batch version given by
\begin{align}
    \abs{L(\hat{v}, \rho)}^2 &\approx L^{\top}(\hat{v}, \rho; \mathbb{A}_1) L(\hat{v}, \rho; \mathbb{A}_2),\label{eq_approxG2}\\
    L(\hat{v}, \rho; \mathbb{A}_i) &:= \frac{T}{\abs{\mathbb{A}_i}} \sum_{(n, m) \in \mathbb{A}_i} \rho(t_n, X_n^m) R(t_n, X_n^m, \xi_n^m; \hat{v}), \quad i = 1, 2, \label{eq_defL}
\end{align}
where $\abs{\mathbb{A}_i}$ denotes the number of elements in the index set $\mathbb{A}_i$, and $h = T/N$ is the time step size.
The index sets $\mathbb{A}_1$ and $\mathbb{A}_2$ are constructed by 
\begin{equation}\label{eq_defA1A2}
\begin{aligned}
    &\mathbb{A}_i = \mathbb{N}_i \times \mathbb{M}_i, \quad \mathbb{N}_i \subset \{0, 1, \cdots, N-1\}, \\ 
    &\mathbb{M}_i \subset \{1, 2, \cdots, M\}, \quad \mathbb{M}_1 \cap \mathbb{M}_2 = \emptyset
\end{aligned}
\end{equation}
for $i = 1, 2$. 
Here, inspired by \cite{Guo2022Monte,hu2024sdgd}, the disjointness between $\mathbb{M}_1$ and $\mathbb{M}_2$ is used to ensure the mini-batch estimation in \eqref{eq_approxG2} is unbiased. 
The sets $\mathbb{N}_1$ and $\mathbb{N}_2$ can be chosen randomly, without replacement, from $\{0, 1, \cdots, N-1\}$, or simply set to $\{0, 1, \cdots, N-1\}$ if $N$ is not too large. 
For clarity, we summarize the notation used in \eqref{eq_defL} in \Cref{tab_notations}.

\begin{table}[t]
    \centering
    \caption{Notations used in the mini-batch loss function~\eqref{eq_defL}.}\label{tab_notations}
    \resizebox{\textwidth}{!}{
    \begin{tabular}{l|l|l|l}
    \hline
    Notation & Description & Notation & Description \\
        \hline
        $\mathbb{N}_1, \mathbb{N}_2$ & Index set of time steps; see \eqref{eq_defA1A2}                   & $\mathbb{M}_1, \mathbb{M}_2$                     & Index set of sample paths; see \eqref{eq_defA1A2}                      \\
        $\mathbb{A}_i$           & Mini-batch index set $\mathbb{N}_i \times \mathbb{M}_i$, $i = 1, 2$ & $\abs{\mathbb{A}_i}$                     & Number of elements in $\mathbb{A}_i$, $i = 1, 2$           \\
        $h$                      & Time step size                                & $N$                                              & Number of time steps                          \\
        $M$                  & Number of sample paths retained in memory & $X_n^m$                                          & $m$-th sample path at $t = t_n$                     \\
        $\xi_n^m$                & $m$-th sample of $\xi$ at time $t_n$          & $R(t_n, X_n^m, \xi_n^m; v)$                      & Residual at $(t_n, X_n^m)$, defined in \eqref{eq_defR}                     \\
        \hline
    \end{tabular}
    }
\end{table}

\begin{remark}\label{rmk_offlinepath}
    The simulation of the sample paths $X_n^m$ in \eqref{eq_pathX} is sequential in time, which can be computationally expensive.
    However, since these samples are only used as spatial sampling points, it is not necessary to update them at every iteration.
    Instead, they can be generated offline and updated less frequently; for example, by updating the sample paths every $I_0$ iterations, where $I_0$ is a tunable hyperparameter.
\end{remark}

\subsection{Adversarial learning}

For numerical implementation, we represent the approximator $\hat{v}$ and the test function $\rho$ in \eqref{eq_minmax} using neural networks $v_{\theta}$ and $\rho_{\eta}$, parameterized by $\theta$ and $\eta$, respectively.
To enforce the terminal condition $v(T, x) = g(x)$, $x \in \R^d$, the neural network $v_{\theta}$ can be constructed as
\begin{equation}\label{eq_netv}
    v_{\theta}(t, x) = \left\{
    \begin{aligned}
        \phi_{\theta}(t, x), \quad &0 \leq t \leq t_{N-1}, \;\;\;  x \in \R^d, \\
        g(x), \quad &t_{N-1} < t \leq T, \;\;\; x \in \R^d,
    \end{aligned}\right.
\end{equation}
where $\phi_{\theta}$ is a neural network parameterized by $\theta$.
The adversarial network $\rho_{\eta}$ plays the role of test functions in classical Galerkin methods for solving PDEs.
By our experimental results, $\rho_{\eta}$ is not necessarily very deep; instead, it can be a shallow network with enough output dimensionality.
Following the multiscale neural network ideas in \cite{Liu2020Multi},
we consider a typical $\rho_{\eta}$ given by
\begin{equation}\label{eq_defrhonet}
    \rho_{\eta}(t, x) = \sin \br{\Lambda \br{W_1 t + W_2 x + b }} \in \R^r, 
\end{equation}
where $\eta := (W_1, W_2, b) \in \R^r \times \R^{r\times d} \times \R^r$ is the trainable parameter; $\Lambda(\cdot)$ is a scale layer defined by 
\begin{equation}\label{eq_defLamb}
    \Lambda(y_1, y_2, \cdots, y_r) = \br{c_1 y_1, c_2 y_2, \cdots, c_r y_r}^{\top} \in \R^r, \quad c_i := 1 + (i - 1) c
\end{equation}
for $y_i \in \R$ with $c > 0$ is a non-trainable hyperparameter; 
$\sin(\cdot)$ is the activation function applied to the outputs of $\Lambda$ in an element-wise manner.

Following the adversarial training approach proposed by \cite{Zang2020Weak}, the minimax problem \eqref{eq_minmax} can be solved by alternately minimizing $L^{\top}(v_{\theta}, \rho_{\eta}; \mathbb{A}_1) \cdot L(v_{\theta}, \rho_{\eta}; \mathbb{A}_2)$ over $\theta$ and maximizing $L^{\top}(v_{\theta}, \rho_{\eta}; \mathbb{A}_1) L(v_{\theta}, \rho_{\eta}; \mathbb{A}_2)$ over $\eta$ through stochastic gradient algorithms.
The detailed training procedure is outlined in Algorithm \ref{alg_pde}.

\begin{algorithm}[t]
    \caption{Weak-form DRDM for the quasilinear parabolic equation \eqref{eq_pde}}\label{alg_pde}
    \begin{algorithmic}[1]
        \Require 
        $v_{\theta}$/$\rho_{\eta}$: neural networks parameterized by $\theta$/$\eta$;
        $\delta_{1}$/$\delta_{2}$: learning rate for $v_{\theta}$/$\rho_{\eta}$;
        $I$: maximum number of iterations;
        $J$/$K$: number of $\theta$/$\eta$ updates per iteration;
        $I_0$ and $r\%$: interval for updating sample paths and the percentage of paths updated;
        $M$: total number of sample paths retained in memory;
        $P_{\xi}$: distribution of $\xi$ satisfying \eqref{eq_condi_xi}.
        
        \State Initialize $v_{\theta}$ and $\rho_{\eta}$
        \State Generate the samples $\{(X_n^{m}, \xi_n^{m}): 0 \leq n \leq N, 1 \leq m \leq M\}$ by \eqref{eq_pathX} with $\hat{v} = v_{\theta}$ and $\xi_n^m \sim P_{\xi}$ i.i.d. for all $n, m$
        
        \For{$i = 1, 2, \cdots, I-1$}

        \State Sample disjoint index subsets $\mathbb{A}_1, \mathbb{A}_2$ per \eqref{eq_defA1A2}
        
        \For{$j = 0, 1, \cdots ,J-1$}
        \State $\theta \leftarrow \theta - \delta_{1} \nabla_{\theta} \bbr{L^{\top}(v_{\theta}, \rho_{\eta}; \mathbb{A}_1) L(v_{\theta}, \rho_{\eta}; \mathbb{A}_2)}$ \Comment{$L$ is given by \eqref{eq_defL}}
        \EndFor
        \For{$k = 0, 1, \cdots, K-1$}
        \State $\eta \leftarrow \eta + \delta_{2} \nabla_{\eta} \bbr{L^{\top}(v_{\theta}, \rho_{\eta}; \mathbb{A}_1) L(v_{\theta}, \rho_{\eta}; \mathbb{A}_2)}$
        \EndFor
        
        \If{$i$ is divisible by $I_0$} 
        \State Update $r\%$ of the sample paths $\{X_n^m\}_{n=0}^N$ by \eqref{eq_pathX} with $\hat{v} = v_{\theta}$
        \EndIf
        
        \EndFor
        \Ensure $v_{\theta}$
    \end{algorithmic}
\end{algorithm}

\subsection{Extension to HJB equations}

We now consider the HJB equation, which is a typical high-dimensional PDE arising in dynamic programming for stochastic optimal control problems (SOCPs) \cite{Yong1999Stochastic}.
The HJB equation and its terminal condition take the form
\begin{equation}\label{eq_hjb}
    \inf_{\kappa \in U} \bbr{\mathcal{D}^{\kappa} v(t, x) + c(t, x,\kappa)} = 0, \quad (t, x) \in [0, T) \times \R^d
\end{equation}
with a terminal condition $v(T, x) = g(x)$, $x \in \R^d$, 
where $U \subset \R^m$ is a given control set, and $\mathcal{D}^{\kappa}$ denotes the second-order differential operator parameterized by $\kappa \in U$, defined as
\begin{equation}
    \mathcal{D}^{\kappa} := \partial_t + \mu^{\top}(t, x, \kappa) \partial_x + \frac{1}{2} \mathrm{Tr}\rbr{\sigma(t, x, \kappa) \sigma(t, x, \kappa)^{\top} \partial_{xx}}.
\end{equation}
Here, the functions $\mu$, $\sigma$, and $c$ are defined on $[0, T) \times \R^d \times U$ and take values in $\R^d$, $\R^{d\times q}$, and $\R$, respectively.

To facilitate the application of the DRDM, we adopt the idea of a policy improvement algorithm \cite{Al2022Extensions,cai2024socmartnet}, and decompose the HJB equation~\eqref{eq_hjb} into two subproblems:
\begin{gather}
    u(t, x) = \argmin_{\kappa \in U} \bbr{\mathcal{D}^{\kappa} v(t, x) + c(t, x, \kappa)}, \label{eq_utx}\\
    \mathcal{D}^u v(t, x) + c(t, x, u(t, x)) = 0, \label{eq_vtx} 
\end{gather}
for $(t, x) \in [0, T) \times \R^d$ with the shorthand $\mathcal{D}^u := \mathcal{D}^{u(t, x)}$.
In the context of SOCPs, the function $u$ denotes the optimal state-feedback control.

Applying the random difference approximation~\eqref{eq_rdo} to the operator $\mathcal{D}^u$ in \eqref{eq_utx} and \eqref{eq_vtx}, and then reformulating the resulting equations into weak forms as described in \cref{sec_galerkin}, we arrive at the following formulation for solving the HJB equation:
\begin{equation}\label{eq_wpmfhjb}
    \min_{\hat{v} \in \mathcal{V}} \; \sup_{\rho \in \mathcal{T}} \abs{L(\hat{u}, \hat{v}, \rho)}^2,  \quad \min_{\hat{u} \in \mathcal{U}_{\mathrm{ad}}} L(\hat{u}, \hat{v}, 1),
\end{equation}
where the sets $\mathcal{T}$ and $\mathcal{V}$ are given in \eqref{eq_defrho} and \eqref{eq1_defV}, respectively;
$\mathcal{U}_{\mathrm{ad}}$ denotes the set of admissible feedback controls $\hat{u}$ mapping $[0, T] \times \R^d$ to $U$;
the loss function $L$ is defined by
\begin{equation}\label{eq_defGhjb}
    L(\hat{u}, \hat{v}, \rho) := \int_0^{T-h} \E{\rho(t, X_t) R(t, X_t, \xi; \hat{u}, \hat{v})} \di t
\end{equation}
with $\xi$ being a random variable satisfying \eqref{eq_condi_xi}, and
$R$ the residual function defined by
\begin{equation}\label{eq_defDelMuv}
\begin{aligned}
    R\br{t, x, z; \hat{u}, \hat{v}} :=\;& \hat{v}\vbr{\big}{t+h, x + \mu(t, x, \hat{u}(t, x)) h + \sigma(t, x, \hat{u}(t, x)) \sqrt{h} z}\\
    & - \hat{v}(t, x) + h c(t, x, \hat{u}(t, x)). 
\end{aligned}
\end{equation}
The notation $1$ in $L(\hat{u}, \hat{v}, 1)$ denotes a constant test function with fixed output $1$.
In the second equation of \eqref{eq_wpmfhjb}, the optimal control function $\hat{u}$ is obtained by minimizing the averaged loss $L(\hat{u}, \hat{v}, 1)$ over the sampled time-space points $(t, x)$.
This approach, first proposed in \cite{cai2024socmartnet,Al2022Extensions,cai2024martnet}, avoids the computational cost of solving the minimization problem in \eqref{eq_utx} at every individual time-space point.

For numerical implementation, 
the functions $\hat{u}$, $\hat{v}$ and $\rho$ in \eqref{eq_wpmfhjb} are represented by neural networks $u_{\alpha}$, $v_{\theta}$, and $\rho_{\eta}$, parameterized by $\alpha$, $\theta$, and $\eta$, respectively.
Similar to \cref{sec_mini_batch}, the loss function $L$ in \eqref{eq_defGhjb} is approximated by its mini-batch version, defined as
\begin{equation}\label{eq_defGhjbA}
    L(\hat{u}, \hat{v}, \rho; \mathbb{A}_i) := \frac{T}{\abs{\mathbb{A}_i}} \sum_{(n, m) \in \mathbb{A}_i} \rho(t_n, X_n^m) R(t_n, X_n^m, \xi_n^m; \hat{u}, \hat{v}), \quad i = 1, 2,
\end{equation}
where $\mathbb{A}_i$ are still given by \eqref{eq_defA1A2}, and $\{X_{n}^m\}_{n=0}^N$ are samples generated by
\begin{equation}\label{eq_defhatXpil0_hjb}
    X_{n+1}^m = X_n^m + \mu\br{t_n, X_n^m, \hat{u}(t_n, X_n^m)}\,h + \sigma\br{t_n, X_n^m, \hat{u}(t_n, X_n^m)}\,\sqrt{h}\,\xi_{n+1}^{m}
\end{equation}
for $n = 0, \cdots, N-1$ and $m = 1, \cdots, M$, with each $\xi_{n}^{m}$ being i.i.d. samples of $\xi$ appearing in \eqref{eq_defGhjb}.
The sample paths $\{X_n^m\}_{n=0}^N$ can be generated offline and updated less frequently, as discussed in \Cref{rmk_offlinepath}.
The detailed training procedure for solving the HJB equation~\eqref{eq_hjb} is summarized in \Cref{alg_amnet}.

\begin{remark}
    An adversarial learning method for solving HJB equations was also proposed in \cite{wang2022is}, motivated by the theoretical result that $L^{\infty}$ loss is preferable to $L^2$ loss for HJB equations.
    Their approach differs from ours in \Cref{alg_amnet}, as our adversarial learning method is based on the minimax formulation~\eqref{eq_wpmfhjb} derived from a Galerkin framework, rather than the $L^{\infty}$ loss.
    Furthermore, similar to \cite[Theorem 4.3]{wang2022is}, it remains a meaningful direction for future research to provide a rigorous theoretical justification for the minimax formulation~\eqref{eq_wpmfhjb}, specifically, to prove that the error of $v$ for the HJB equation~\eqref{eq_hjb} can be bounded by its residual with respect to the minimax problem~\eqref{eq_wpmfhjb}.
\end{remark}

\begin{remark}
    In the context of SOCPs, the sampling method described in \eqref{eq_defhatXpil0_hjb} essentially uses a controlled SDE to guide the exploration of the state space $\R^d$. 
    This approach has been proposed by \cite{Li2024Neural}, which demonstrates that leveraging controlled SDEs for state-space exploration can significantly enhance the performance of SOCP solvers, particularly in problems with complex dynamics.
\end{remark}

\begin{algorithm}[t]
    \caption{Weak-form DRDM for solving the HJB equation~\eqref{eq_hjb}}\label{alg_amnet}
    \begin{algorithmic}[1]
        \Require 
        $v_{\theta}$/$u_{\alpha}$/$\rho_{\eta}$: neural networks parameterized by $\theta$/$\alpha$/$\eta$;
        $\delta_0$/$\delta_{1}$/$\delta_{2}$: learning rates for the networks $v_{\theta}$/$u_{\alpha}$/$\rho_{\eta}$;
        $I$: maximum number of iterations;
        $J$/$K$: number of $(\theta, \alpha)$/$\eta$ updates per iteration;
        $I_0$ and $r\%$: interval for updating spatial samples and the percentage of samples updated;
        $M$: total number of sample paths retained in memory;
        $P_{\xi}$: distribution of $\xi$ satisfying \eqref{eq_condi_xi}.
        
        \State Initialize the networks $u_{\alpha}$, $v_{\theta}$ and $\rho_{\eta}$
        \State Generate the samples $\{(X_n^{m}, \xi_n^{m}): 0 \leq n \leq N, 1 \leq m \leq M\}$ by \eqref{eq_defhatXpil0_hjb} with $\hat{u} = u_{\alpha}$ and $\xi_n^m \sim P_{\xi}$ i.i.d. for all $n, m$
        \For{$i = 0, 1, \cdots, I-1$}
        \State Sample disjoint index subsets $\mathbb{A}_1, \mathbb{A}_2$ per \eqref{eq_defA1A2}
        \For{$j = 0, 1, \cdots ,J-1$}
        \Comment{$L$ is given by \eqref{eq_defGhjbA}}
        \State $\theta \leftarrow \theta - \delta_{0} \nabla_{\theta} \bbr{L^{\top}(u_{\alpha}, v_{\theta}, \rho_{\eta}; \mathbb{A}_1) L(u_{\alpha}, v_{\theta}, \rho_{\eta}; \mathbb{A}_2)}$ 
        \State $\alpha \leftarrow \alpha - \delta_{1} \nabla_{\alpha} L(u_{\alpha}, v_{\theta}, 1; \mathbb{A}_1 \cup \mathbb{A}_2)$
        \EndFor
        \For{$k = 0, 1, \cdots, K-1$}
        \State $\eta \leftarrow \eta + \delta_{2} \nabla_{\eta} \bbr{L^{\top}(u_{\alpha}, v_{\theta}, \rho_{\eta}; \mathbb{A}_1) L(u_{\alpha}, v_{\theta}, \rho_{\eta}; \mathbb{A}_2)}$
        \EndFor
        
        \If{$i$ is divisible by $I_0$} 
        \State Update $r\%$ of the sample paths $\{X_n^m\}_{n=0}^N$ by \eqref{eq_defhatXpil0_hjb} with $\hat{u}=u_{\alpha}$
        \EndIf
        
        \EndFor
        \Ensure $u_{\alpha}$ and $v_{\theta}$
    \end{algorithmic}
\end{algorithm}

\section{Numerical experiments}\label{sec_numexp}

In this section, we present numerical experiments to demonstrate the performance of the weak-form DRDM in solving high-dimensional parabolic and HJB equations. 
The detailed experimental parameters are provided in \cref{sec_param_setting}.
In reporting the numerical results, the abbreviations ``RE'', ``SD'', ``RT'', ``vs'', and ``Iter.'' refer to ``relative error'', ``standard deviation'', ``running time'', ``versus'', and ``iteration'', respectively.

To visualize the numerical solutions, we plot the curves of both the numerical and exact solutions at $t = 0$ for $x$ lying on the curves $S_2$ and $S_3$ in $\R^d$, defined as
\begin{align}
    &S_2 := \bbr{s \B{1}_d: s \in [-1, 1]}, \quad \B{1}_d := (1, 1, \cdots, 1)^{\top} \in \R^{d}.\label{eq_defS2}\\
    &S_3 := \bbr{\B{l}(s): s \in [-1, 1]}, \quad \B{l}(s) := \vbr{\big}{l_1(s), l_2(s), \cdots, l_d(s)}^{\top} \in \R^d \label{eq_defS3}
\end{align}
where $l_i(s) := s \times \mathrm{sgn}\br{\sin(i)} + \cos(i + \pi s)$ and $\mathrm{sgn}(z) := -1, 0, 1$ for $z < 0$, $= 0$, and $> 0$, respectively. 
Here, $S_2$ is a straight line segment whose length grows with the dimension $d$, while $S_3$ is a curve that winds through the space $\R^d$. 

For the approximation $\hat{v}$ to the exact solution $v$, we define the relative $L^1$ and $L^{\infty}$ errors at time $t = t_n$ as follows:
\begin{align}
    \mathrm{RE}_1(t_n) &:= \vbr{\Big}{\sum_{x \in D_n} \abs{\hat{v}\br{t_n, x} - v(t_n, x)}} \big/ \vbr{\Big}{\sum_{x \in D_n} \abs{v(t_n, x)}}, \label{eq_defRelL1}\\
    \mathrm{RE}_{\infty}(t_n) &:= \vbr{\Big}{\max_{x \in D_n} \abs{\hat{v}\br{t_n, x} - v(t_n, x)}} \big/ \vbr{\Big}{\max_{x \in D_n} \abs{v(t_n, x)}},\label{eq_defRelLinf}
\end{align}
where $D_n = \{X_n^m\}_{m=1}^{\bar{M}}$ denotes the set of spatial points sampled according to \eqref{eq_pathX} for the quasilinear PDE~\eqref{eq_pde}, and \eqref{eq_defhatXpil0_hjb} for HJB equations.
Unless otherwise specified, $\mathrm{RE}_1$ and $\mathrm{RE}_{\infty}$ refer to the errors evaluated at $t = 0$.
The initial points $\{X_0^m\}_{m=1}^{\bar{M}}$ are chosen as uniformly spaced grid points along the curves $S_2$ and $S_3$. 
For evaluating the error trajectories $t_n \mapsto \mathrm{RE}_1(t_n)$ and $t_n \mapsto \mathrm{RE}_{\infty}(t_n)$ (see, e.g., the third column of \Cref{fig_HP9abcd}), we set $\bar{M} = 8$. 
For evaluating $\mathrm{RE}_1$ and $\mathrm{RE}_{\infty}$ as functions of the number of iterations (see, e.g., the fourth column of \Cref{fig_HP9abcd}), we use $\bar{M} = 1000$. 

\subsection{Convection-diffusion equation with a steep gradient solution}\label{sec_cde}

\begin{figure}[!t]
    \centering
    \subfloat[$s \mapsto v(0, s \B{1}_d)$, $c = 1$]{\includegraphics[width=0.34\textwidth]{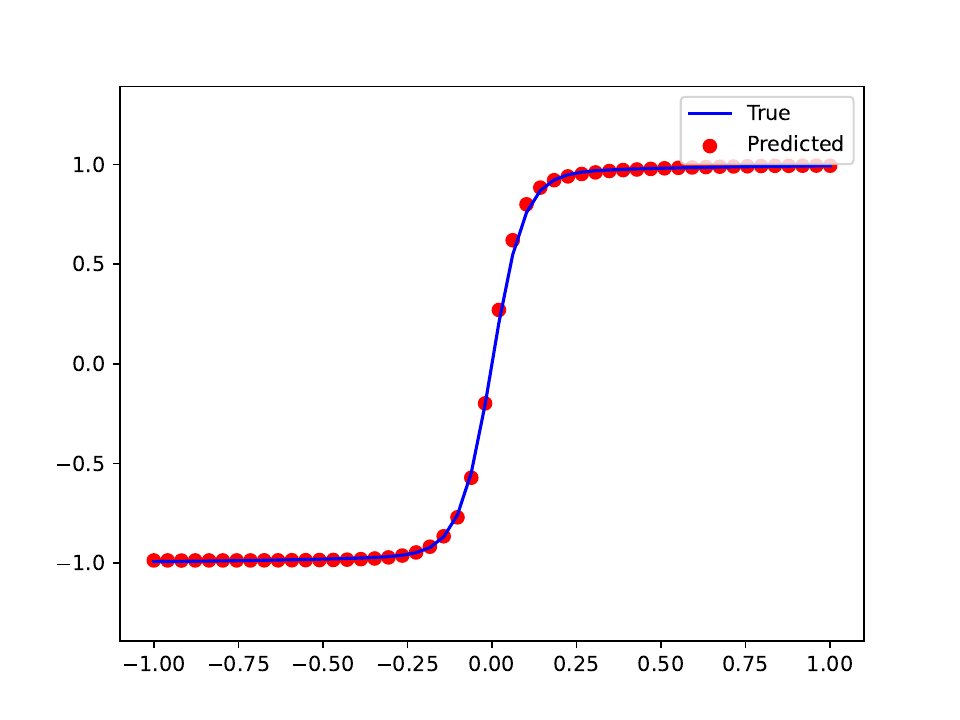}}
    \subfloat[Loss vs Iter., $c = 1$]{\includegraphics[width=0.31\textwidth]{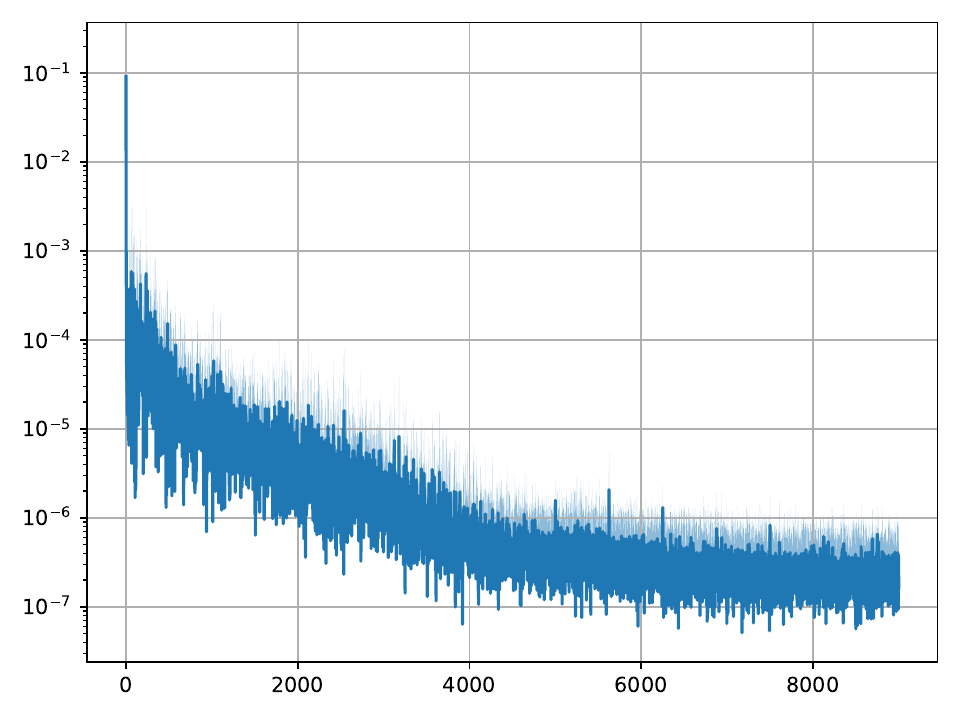}}
    \subfloat[RE$_1$ vs Iter., $c = 1$]{\includegraphics[width=0.31\textwidth]{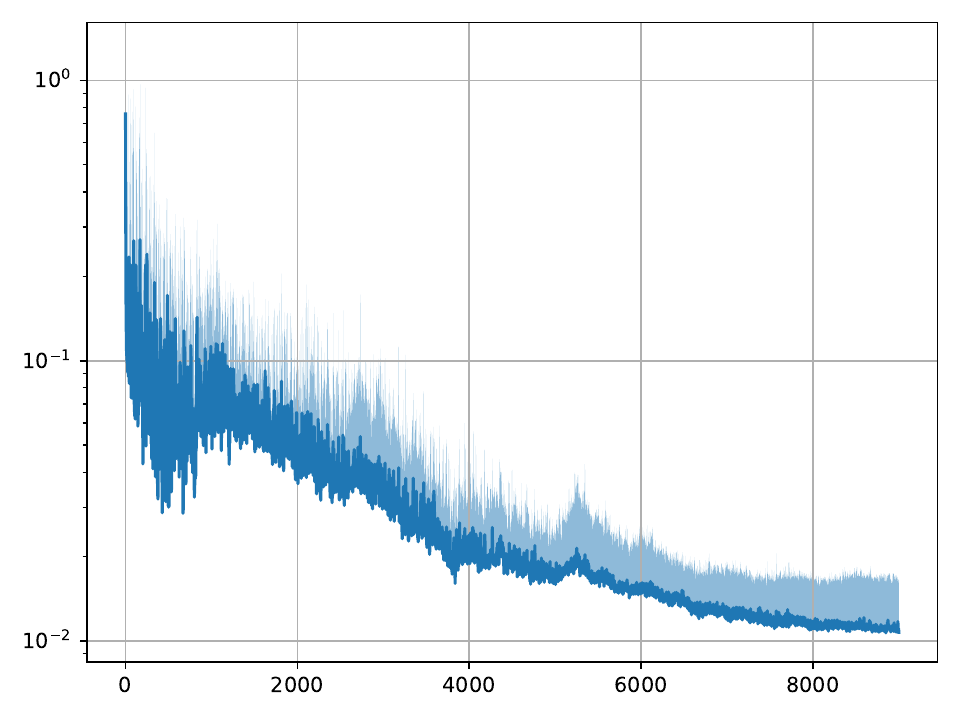}}
    \\
    \subfloat[$s \mapsto v(0, s \B{1}_d)$, $c = 5$]{\includegraphics[width=0.34\textwidth]{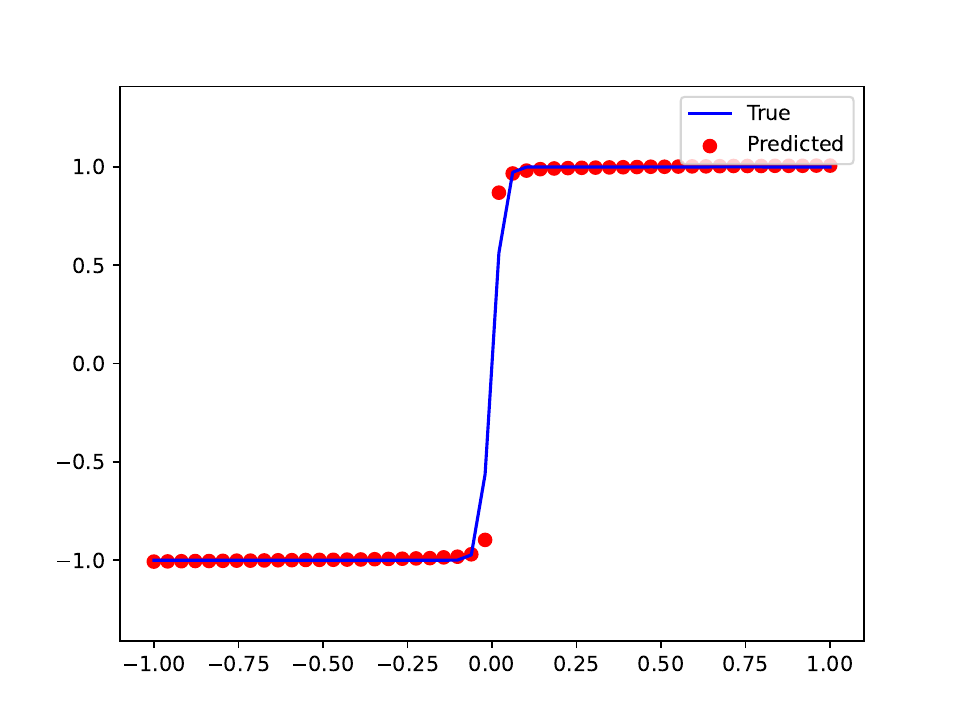}}
    \subfloat[Loss vs Iter., $c = 5$]{\includegraphics[width=0.31\textwidth]{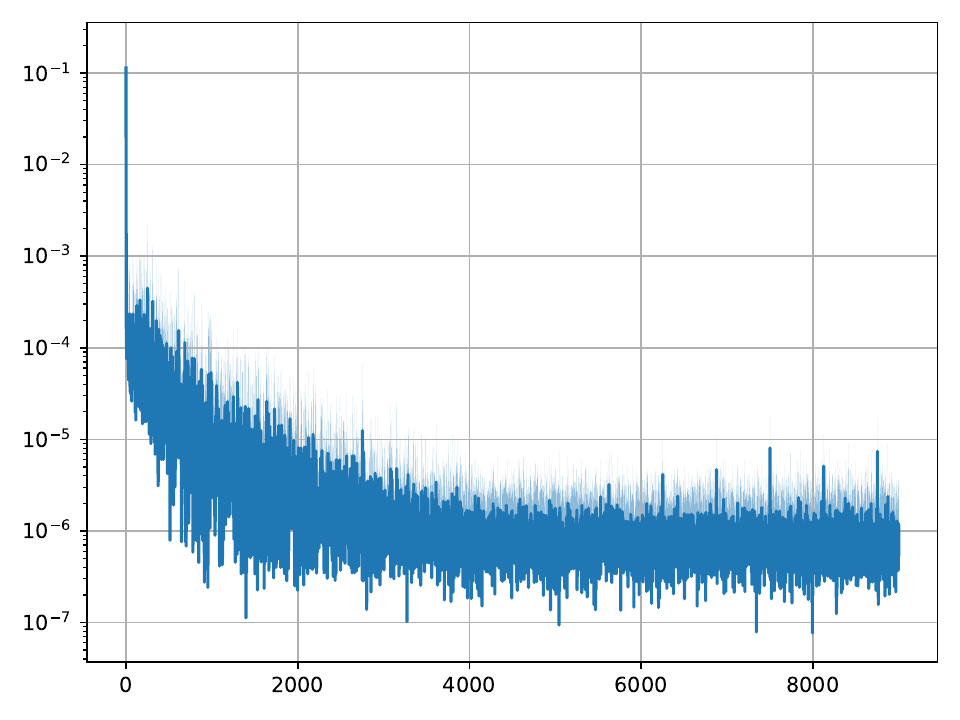}}
    \subfloat[RE$_1$ vs Iter., $c = 5$]{\includegraphics[width=0.31\textwidth]{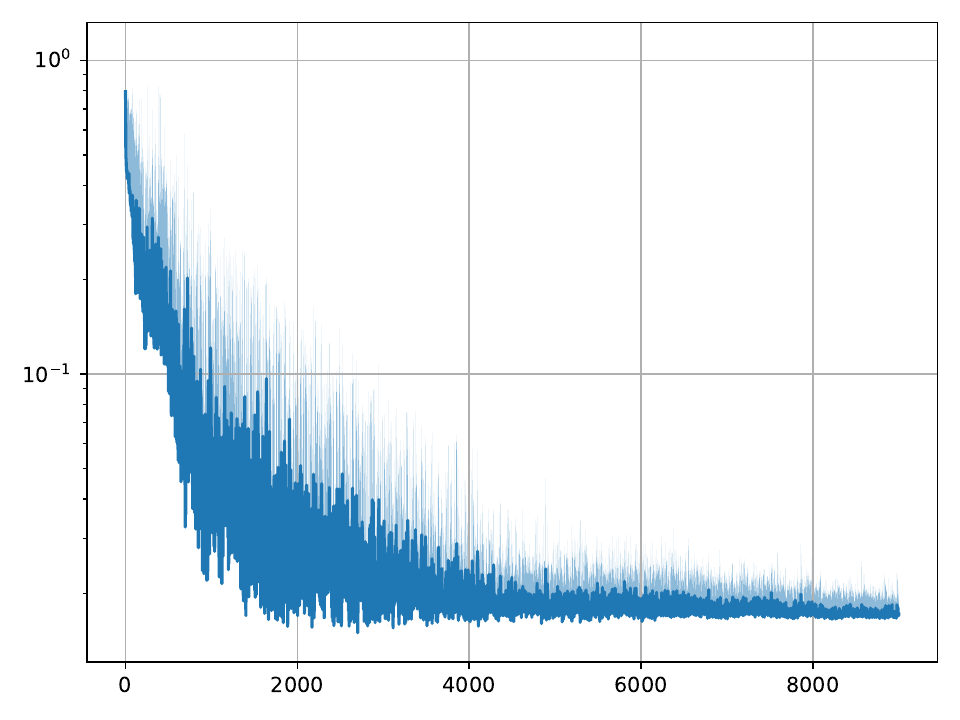}}
    \caption{
    Numerical results of \Cref{alg_pde} for the convection–diffusion equation \eqref{eq_linearPDE} with $d = 10^3$.
    The ``Loss'' in the second column refers to $\big|L^{\top}(v_{\theta}, \rho_{\eta}; \mathbb{A}_1)\,L(v_{\theta}, \rho_{\eta}; \mathbb{A}_2)\big|$ in \Cref{alg_pde}.
    The shaded regions in the second and third columns show the mean $\pm 2$ SD of the relative errors over 5 independent runs.
    }\label{fig_cde}  
\end{figure}

We consider the following linear parabolic PDE:
\begin{equation}\label{eq_linearPDE}
\left\{
\begin{aligned}
    &\vbr{\Big}{\partial_t + \mu^{\top} \partial_{x} + \frac{\bar{\sigma}^2}{2} \sum_{i=1}^d \partial_{x_i}^2} v(t, x) = 0, \quad (t, x) \in [0, T) \times \R^d,\\
    &v(T, x) = \frac{1}{d} \sum_{i=1}^d \bbr{\tanh(x_i) + \cos(10 x_i)}, \quad x \in \R^d,
\end{aligned}
\right.
\end{equation} 
where $T=2$, $\bar{\sigma}^2 = 0.1$, and the drift coefficient $\mu$ is given by
\begin{equation}\label{eq_defmu_cde}
    \mu(t, x) := c \times \vbr{\Big}{\tanh(10 x_1), \tanh(10 x_2), \cdots, \tanh(10 x_d)}^{\top}, 
\end{equation}
with $c > 0$ reported in the numerical results. 
This problem features a highly oscillatory terminal function $v(T, x)$.
Moreover, the drift coefficient $\mu(t, x)$ exhibits rapid transitions near $x = 0$. 
When $c$ is much larger than $\bar{\sigma}$, the solution exhibits a steep gradient near $(t,x)=(0,0)$; see \Cref{fig_cde,fig_linearPDE}.

According to the linear Feynman-Kac formula \cite[Theorem 8.2.1]{Oksendal2003Stochastic}, the exact solution is
\begin{equation}\label{eq_vtx_cde}
    v(t, x) = \E{v(T, X_T^{t, x})}, \quad (t, x) \in [0, T] \times \R^d,
\end{equation}
where $s \mapsto X_s^{t, x}$ is a $d$-dimensional stochastic process starting from $(t, x)$ and governed by the SDE
\begin{equation}\label{eq_Xstx}
    X_{s}^{t, x} = x + \int_t^s \mu(r, X_r^{t, x}) \di r + \int_t^s \bar{\sigma} \di B_r, \quad s \in [t, T]. 
\end{equation}
Here, $B$ denotes a $d$-dimensional standard Brownian motion.
The reference solution $v(t, x)$ is computed by approximating the expectation in \eqref{eq_vtx_cde} using the Monte Carlo method, where $10^6$ sample paths of $s \mapsto X_s^{t, x}$ are generated via the Euler-Maruyama scheme applied to \eqref{eq_Xstx} with a time step size of $T / 100$.\label{txt_stepsize}

We apply \Cref{alg_pde} to solve \eqref{eq_linearPDE} with $d = 10^3$ and $c = 1$ or $5$.
The numerical results are presented in \Cref{fig_cde}.
As shown in the first column of \Cref{fig_cde}, the numerical solutions closely match the exact solutions, and accurately capture the steep gradient at $x = 0$.
This performance is primarily due to the sampling method described in \cref{sec_pathsamp}, which is further analyzed in \cref{sec_ablation}.

\subsection{Quasilinear parabolic PDEs}\label{sec_qlpde}

We now evaluate the performance of \Cref{alg_pde} on several high-dimensional quasilinear parabolic (QLP) equations with varying degrees of nonlinearity. 
In particular, we consider the following QLP equations with $d=10^4$:
\begin{itemize}
    \item QLP-1: Given by \eqref{eq_pde} with a nonlinear diffusion term and an Allen-Cahn-type source term:
    \begin{equation}\label{eq_ac_sdgd}
        \mathcal{D} = \partial_t + \frac{v^2}{2} \sum_{i=1}^d \partial_{x_i}^2, \quad f(t, x, v) = v - v^3 + Q(t, x);
    \end{equation}

    \item QLP-2a: Given by \eqref{eq_pde} with nonlinearities in the drift, diffusion, and source terms:
    \begin{equation}\label{eq_ql_sdgd}
        \mathcal{D} 
        = \partial_t + \left(\frac{v}{2} - 1\right)\! \sum_{i=1}^d \partial_{x_i} 
        + \frac{v^2}{2} \sum_{i=1}^d \partial_{x_i}^2, \quad f(t, x, v) 
        = v^2 + Q(t, x).
    \end{equation}

    \item QLP-2b: A more complex variant of QLP-2a, where the diagonal diffusion coefficient is replaced by a dense one:
    \begin{equation}\label{eq_qld_sdgd}
    \begin{aligned}
        &\mathcal{D} 
        = \partial_t + \frac{1}{d} \left(\frac{v}{2} - 1\right) \sum_{i=1}^d \partial_{x_i}
        + \frac{1}{2 d^2} \sum_{i,j,k=1}^d \sigma_{ik}\sigma_{jk} \,\partial_{x_i}\partial_{x_j}, \\ 
        &f(t, x, v) = v^2 + Q(t, x), \quad \sigma_{ij} = \cos(x_i) + v\,\sin(x_j),  
    \end{aligned}
    \end{equation}
    for $i,j = 1,\cdots,d$.
    It is worth noting that high-dimensional QLP-2b poses a significant challenge for conventional PINNs, as these methods require explicit computation of the full Hessian of $v$ in the loss function.
\end{itemize}
In the above equations, we take $T = 1$.
The functions $Q(t, x)$ and the terminal condition $v(T, x) = g(x)$ are chosen such that the PDE admits an exact solution of the form
\begin{equation}\label{eq_vsdgd}
    v(t, x) = V\vbr{\big}{(t - 0.5) \B{1}_d + x}, \quad (t, x) \in [0, T] \times \R^d,
\end{equation}
where $\B{1}_d := (1, 1, \cdots, 1)^{\top} \in \R^d$, and $V$ is given by
\begin{equation}\label{eq_defVx}
    V(x) := \sum_{i=1}^{d-1} c_i K(x_i, x_{i+1}) + c_d K(x_{d}, x_1)
\end{equation}
with 
\begin{equation*}
    c_i := \br{1.5 - \cos(i \pi/d)}/d, \quad K(x_i, x_j) := \sin\vbr{\big}{x_i + \cos (x_j) + x_j \cos (x_i)}. 
\end{equation*}

\begin{table}[!t]
    \centering
    \caption{Numerical results of \Cref{alg_pde} for various quasilinear PDEs with $d = 10^4$ and $10^5$.
    The convergence histories are presented in \Cref{fig_pde_quasi}.
    }\label{tab_RESQLP}
    \resizebox{\textwidth}{!}{
    \begin{tabular}{c c c c c c c} 
    \toprule
    Equation & $d$ & Mean of $\mathrm{RE}_1$ & SD of $\mathrm{RE}_1$ & Mean of $\mathrm{RE}_{\infty}$ & SD of $\mathrm{RE}_{\infty}$ & RT (s) \\ [0.5ex] 
    \midrule
    QLP-1  & $10^4$  &  1.98E-2&  7.24E-3   & 3.91E-2 & 8.75E-3 & 1585  \\ 
    QLP-1  & $10^5$  &  1.29E-2&  1.40E-3   & 3.20E-2 & 3.66E-3& 3761   \\ 
    QLP-2a & $10^4$  &  2.82E-2&  1.28E-2   & 9.08E-2 & 1.94E-2 & 1589  \\
    QLP-2a & $10^5$  &  4.06E-2&  1.14E-3   & 1.37E-1 & 1.33E-2& 3773   \\
    QLP-2b & $10^4$  &  5.77E-2&  2.00E-3   & 1.12E-1 & 1.18E-2 & 1602  \\
    QLP-2b & $10^5$  &  5.19E-2&  7.69E-4 & 1.05E-1 &  1.09E-2 & 3822      \\
    \bottomrule
    \end{tabular}}
\end{table}
\begin{figure}[!t]
    \centering
    \resizebox{1.0\textwidth}{!}{
    \begin{tabular}{@{}c@{\hspace{1mm}}cccc@{}}
        & \textbf{$s \mapsto v(0, s\B{1}_d)$, $d=10^4$} & \textbf{$s \mapsto v(0, s\B{1}_d)$, $d=10^5$} & \textbf{RE vs Iter., $d=10^4$} & \textbf{RE vs Iter., $d=10^5$}  \\
        
        \adjustbox{valign=m}{\rotatebox[origin=c]{90}{\textbf{QLP-1}}} & 
        \adjustbox{valign=m}{\includegraphics[width=0.33\textwidth]{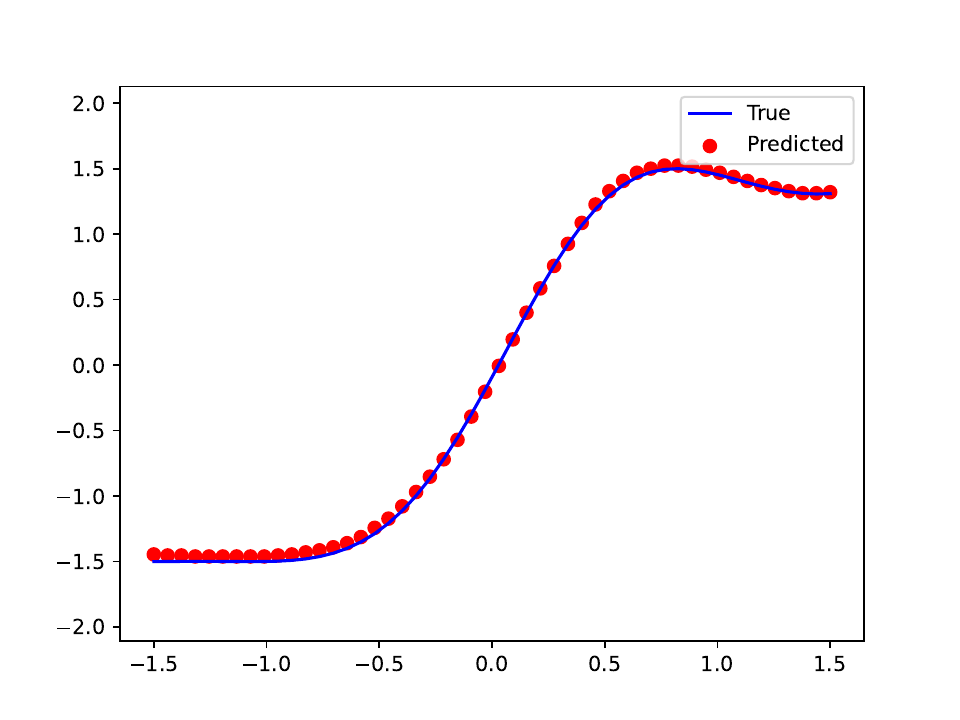}} &
        \adjustbox{valign=m}{\includegraphics[width=0.33\textwidth]{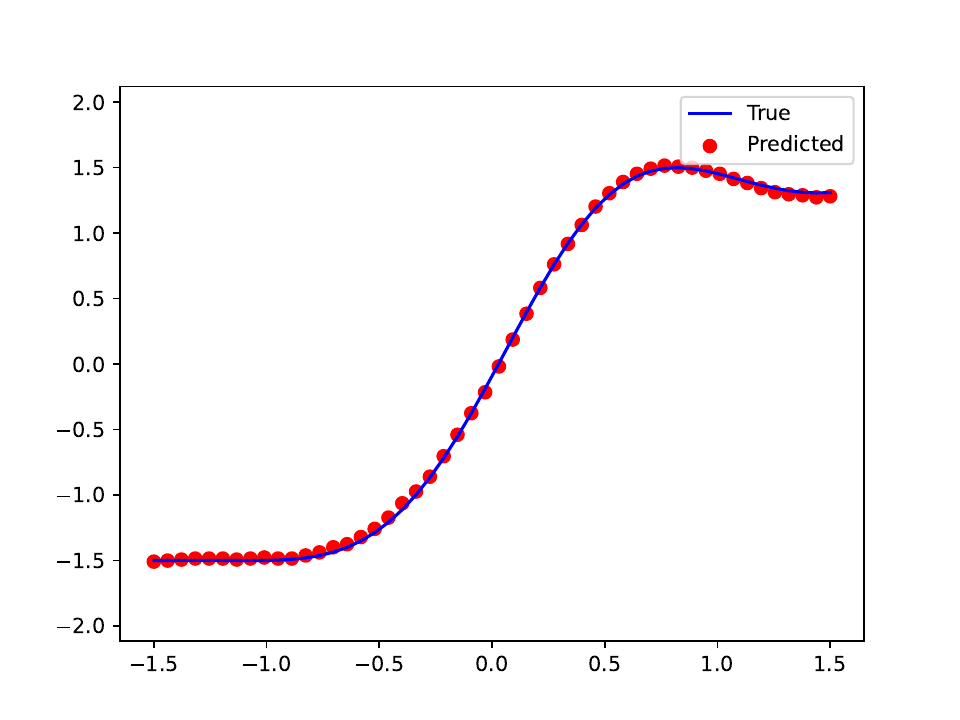}} &
        \adjustbox{valign=m}{\includegraphics[width=0.3\textwidth]{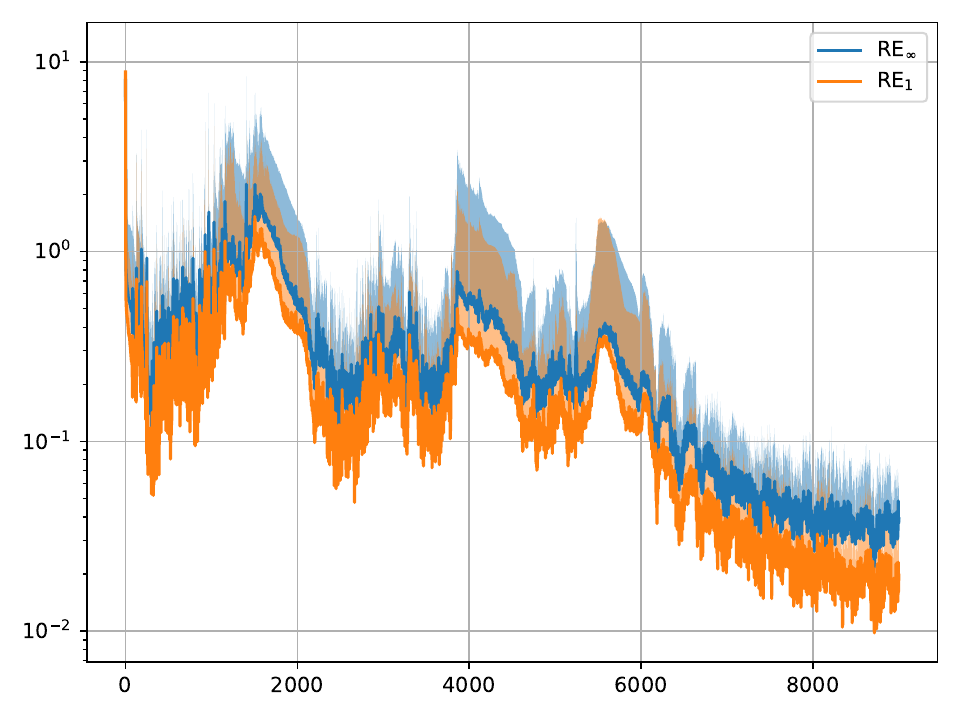}} &
        \adjustbox{valign=m}{\includegraphics[width=0.3\textwidth]{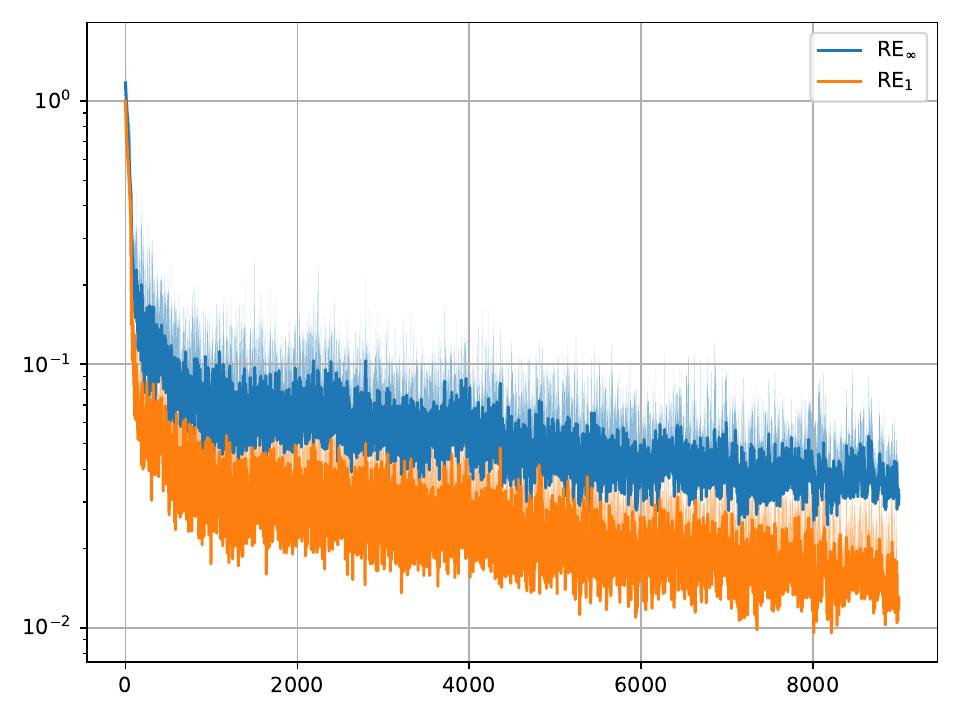}} \\
        \addlinespace[3pt]
        
        \adjustbox{valign=m}{\rotatebox[origin=c]{90}{\textbf{QLP-2a}}} & 
        \adjustbox{valign=m}{\includegraphics[width=0.33\textwidth]{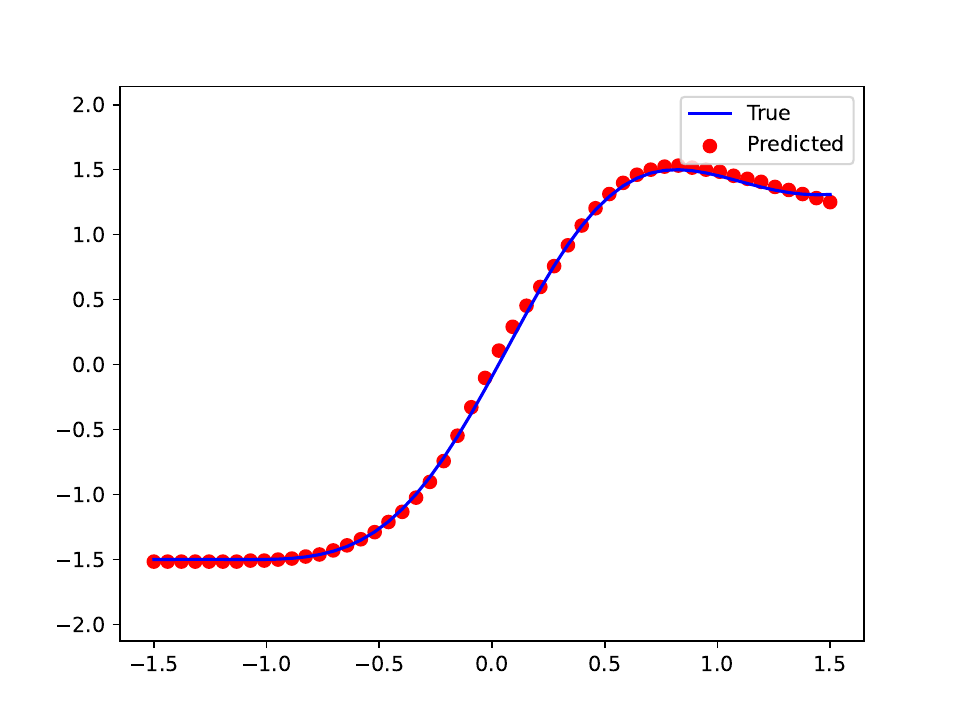}} &
        \adjustbox{valign=m}{\includegraphics[width=0.33\textwidth]{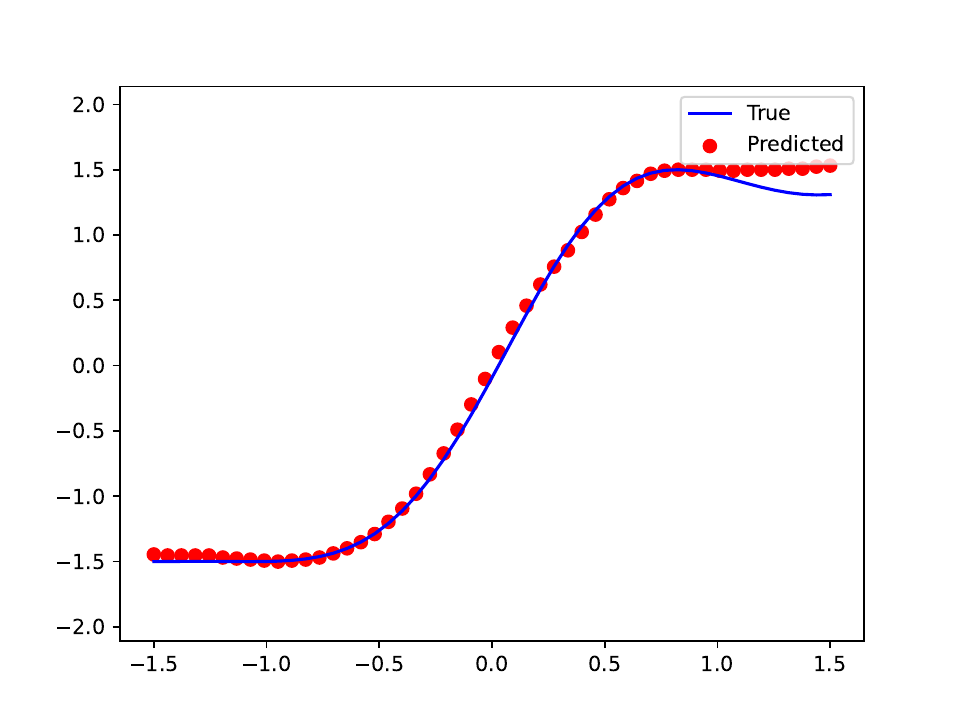}} &
        \adjustbox{valign=m}{\includegraphics[width=0.3\textwidth]{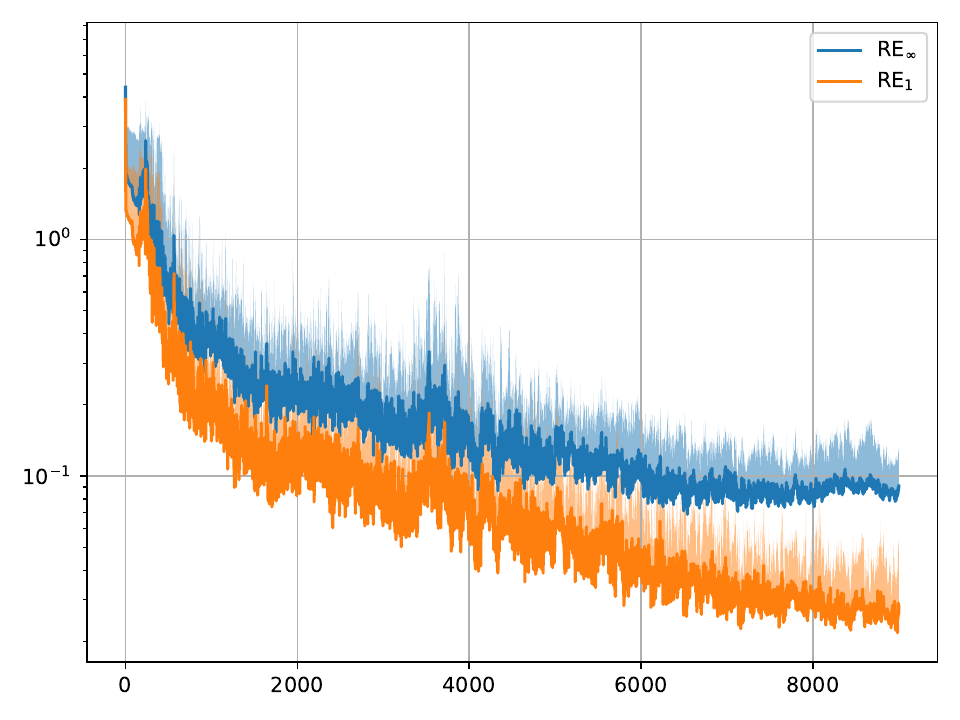}} &
        \adjustbox{valign=m}{\includegraphics[width=0.3\textwidth]{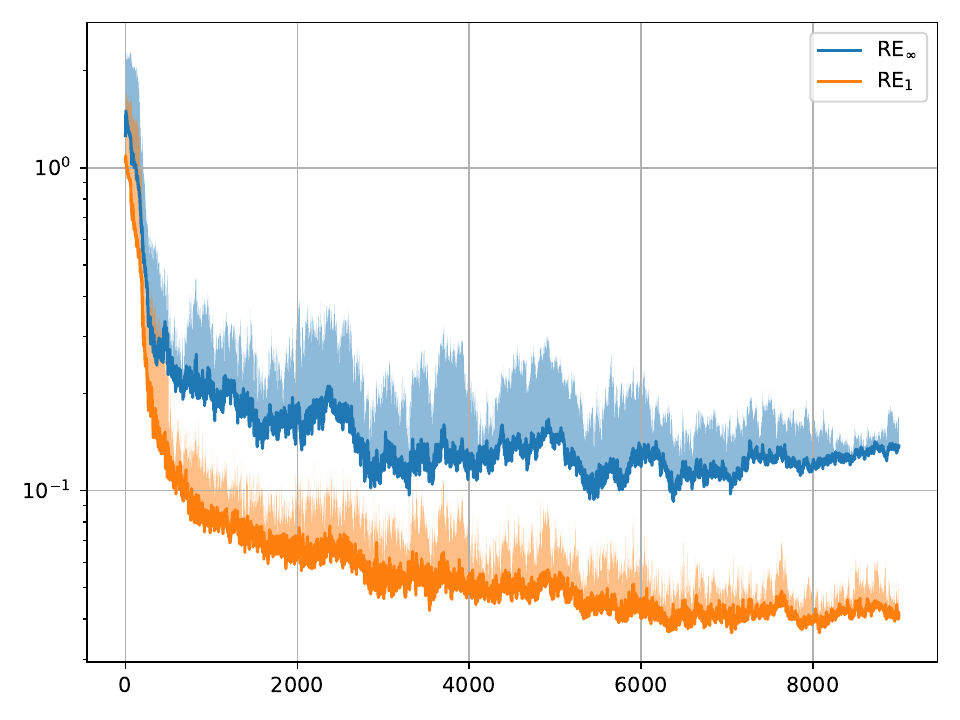}} \\
        \addlinespace[3pt]
        
        \adjustbox{valign=m}{\rotatebox[origin=c]{90}{\textbf{QLP-2b}}} & 
        \adjustbox{valign=m}{\includegraphics[width=0.33\textwidth]{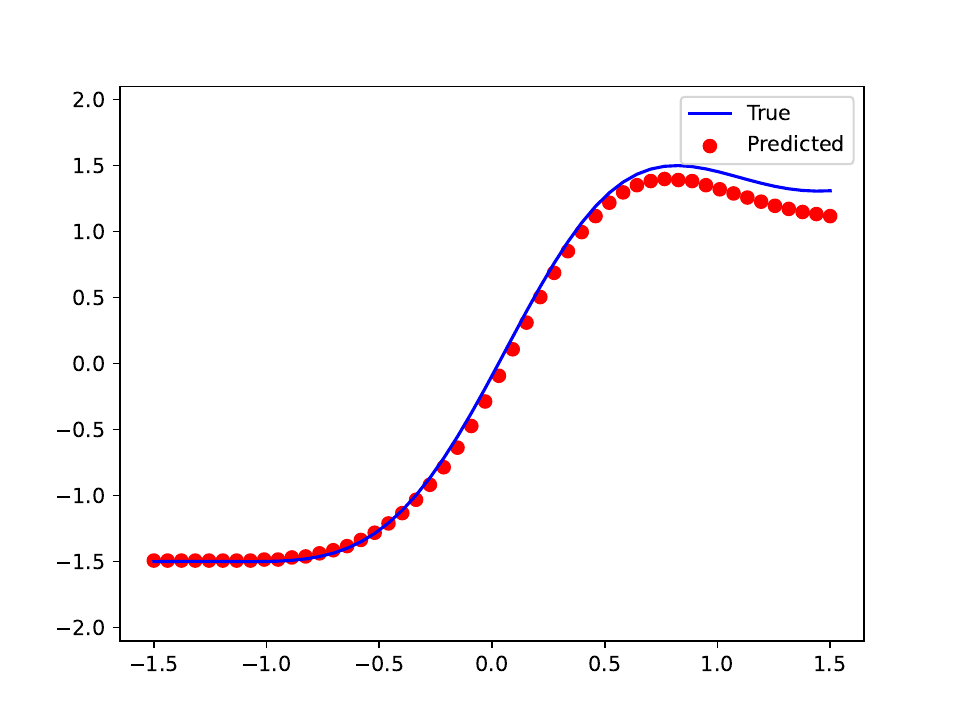}} &
        \adjustbox{valign=m}{\includegraphics[width=0.33\textwidth]{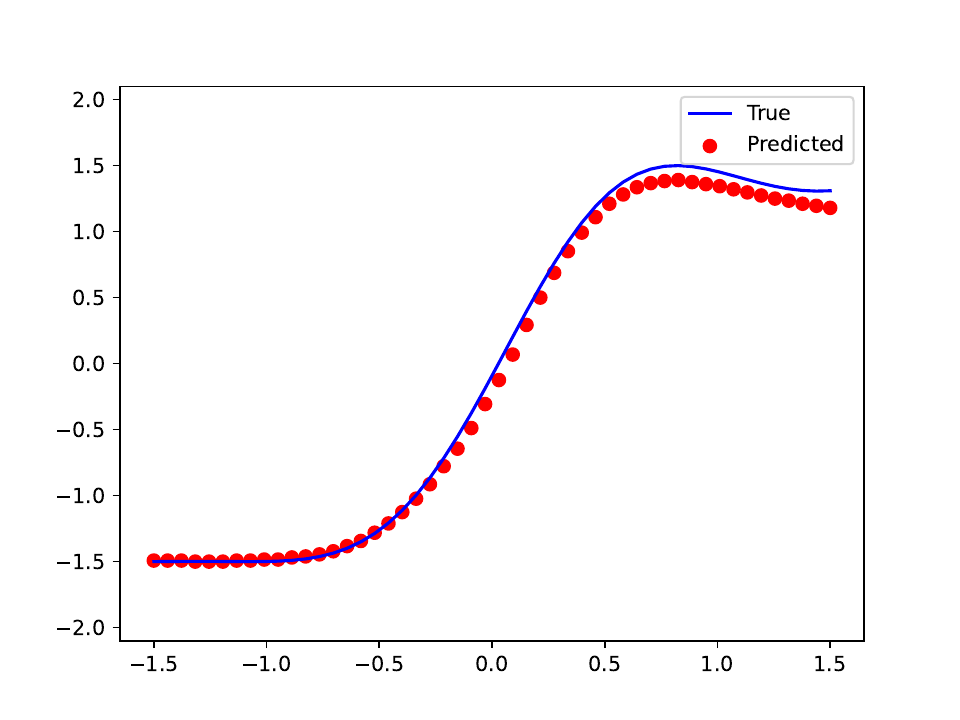}} &
        \adjustbox{valign=m}{\includegraphics[width=0.3\textwidth]{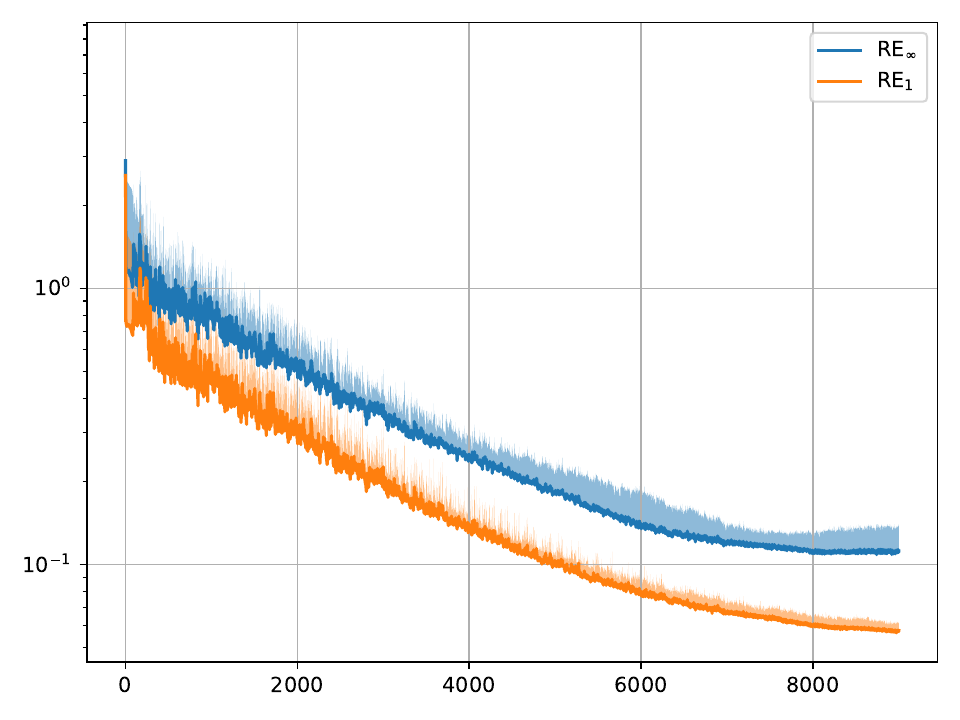}} &
        \adjustbox{valign=m}{\includegraphics[width=0.3\textwidth]{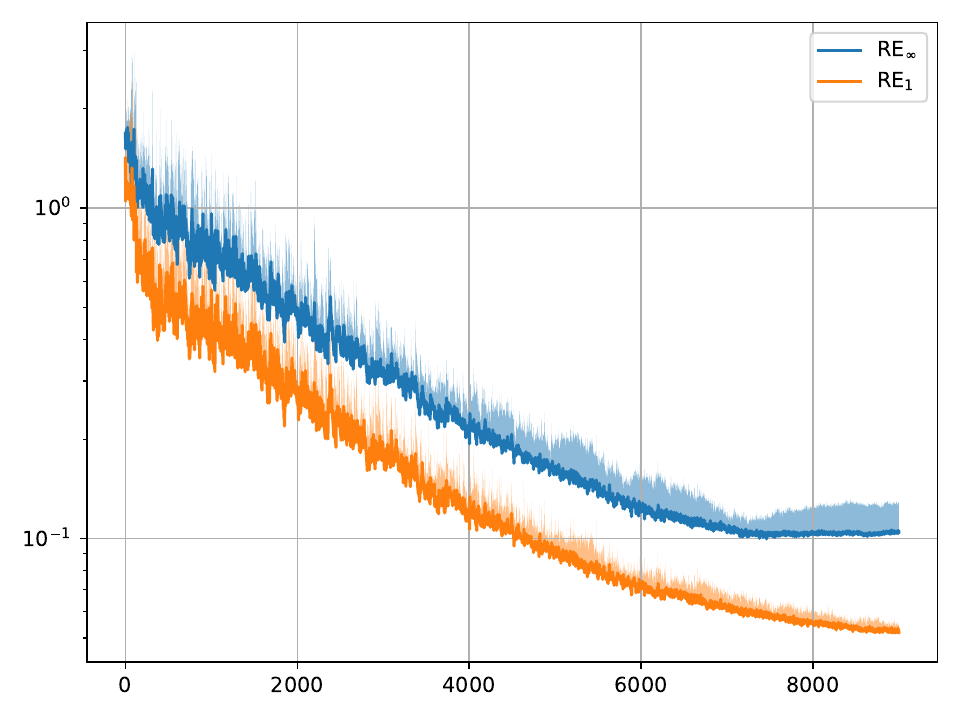}} \\
        \addlinespace[3pt]
    \end{tabular}
    }
    \caption{
    Numerical results of \Cref{alg_pde} for QLP-1, QLP-2a, and QLP-2b in \cref{sec_qlpde} for $d = 10^4$ and $10^5$.
    The shaded regions in the third and fourth columns represent the mean $\pm 2 \times \text{SD}$ of the REs over 5 independent runs, where all the REs are computed at $t = 0$.
    The mean and the SD of the REs, and the RT, are given in \Cref{tab_RESQLP}. 
    }\label{fig_pde_quasi}
\end{figure}

The function $V$ is adapted from the equation (28) in \cite{hu2024sdgd}, with modifications to the coefficients $c_i$ in \eqref{eq_defVx}. 
Unlike \cite{hu2024sdgd}, where $v(0, x)$ is computed for $x$ within the unit ball in $\R^d$, our setting evaluates $v(0, x)$ along the diagonal curve $S_2$ of the unit cube, as defined in \eqref{eq_defS2}. 
The choice of $c_i$ in \eqref{eq_defVx} ensures that $v(0, x)$ remains bounded for $x \in S_2$ as $d$ increases. 
Nevertheless, as discussed in \cite{hu2024sdgd}, the function $V$ still incorporates uneven coefficients $c_i$ and pairwise interactions between the variables $x_i$ and $x_{i+1}$. 
These features increase the complexity of the problem, making the exact solution $v(t, x)$ in \eqref{eq_vsdgd} highly nontrivial and challenging to approximate.
The corresponding numerical results for \Cref{alg_pde} are presented in \Cref{fig_pde_quasi} and \Cref{tab_RESQLP}, demonstrating that our method is both accurate and efficient for solving quasilinear equations in dimensions up to $d = 10^5$.

\subsection{HJB equation}\label{sec_hjb}

We consider the following HJB equation:
\begin{equation}\label{eq_hjb_test}
    \partial_t v + \inf_{\kappa \in \R^d} \bbr{(b + c \sigma \kappa)^{\top} \partial_x v + \frac{1}{2}\abs{\kappa}^2} + \frac{1}{2} \mathrm{Tr}\rbr{\sigma \sigma^{\top} \partial_{xx} v} = 0, 
\end{equation}
where $(t, x) \in [0, T) \times \R^d$, 
with the terminal condition 
\begin{equation*}
    v(T, x) = \ln\vbr{\Big}{1 + \frac{1}{d} \sum_{i=1}^d \br{x_i^2 + 0.5 \sin(10 x_i)}},\quad x \in \R^d, \quad T = 1.
\end{equation*}
Here, $b$ and $\sigma$ are functions of $(t, x)$, taking values in $\R^d$ and $\R^{d \times d}$, respectively, and $c = 2$ is a constant.
This HJB equation arises from the following stochastic optimal control problem (SOCP):
\begin{equation}\label{eq_socp}
    \min_{u \in \mathcal{U}_{\mathrm{ad}}} J(u), \quad J(u) := \int_0^T \frac{1}{2}\E{\abs{u(t, X_t^u)}^2} \di t + \E{g(X_T^u)},
\end{equation}
where $X^u: [0, T] \times \Omega \to \R^d$ is the state process defined by 
\begin{equation}\label{eq_xtu}
    X_t^u = X_0^u + \int_0^t \br{b + c \sigma u}(s, X_s^u) \di s + \int_0^t \sigma(s, X_s^u) \di B_s, \quad t \in [0, T]. 
\end{equation}
The admissible control set $\mathcal{U}_{\mathrm{ad}}$ consists of all Borel measurable functions from $[0, T] \times \R^d$ to $\R^d$ such that the SDE \eqref{eq_xtu} admits a unique strong solution and the expectations in \eqref{eq_socp} are finite.
By the Cole-Hopf transformation (see \cref{sec_constPDE}), the value function $v(t, x)$ is given by 
\begin{equation}\label{eq_vtx_hjb}
    v(t, x) = -\frac{1}{c^2}\ln\br{\E{\exp\br{-c^2 v(T, X_T^0)} \big\vert X_t^0 = x}}
\end{equation}
for $t \in [0, T)$ and $x \in \R^d$,
where $t \mapsto X_t^{0}$ is the state process $X_t^u$ given by \eqref{eq_xtu} with $u(t, x) \equiv 0$.
For reference solutions, we compute $v(t, x)$ by approximating the expectation using the Monte Carlo method, where $10^6$ samples of $X_T^{0}$ are generated via the Euler-Maruyama scheme applied to \eqref{eq_xtu} with a time step size of $T / 100$.
\label{txt_stepsize2}

We consider the following specific HJB equations with $d = 10^4$, which are given by \eqref{eq_hjb_test} with different parameters: 
\begin{itemize}
    \item HJB-1a: $b_i = \sin(t + i + x_{i+1} - 1)$ for $i=1,2,\cdots, d$ with $x_{d+1} := x_1$, and $\sigma = 0.5 I_d$.
    Here, the drift coefficient $b_i$ depends on $x_{i+1}$, introducing correlations among the components of $X_t^u$.

    \item HJB-1b: a variant of HJB-1a with $\sigma = 0.025 I_d$. The reduced diffusion coefficient in this case results in a solution that is less smooth compared to HJB-1a.

    \item HJB-2: $b_i = \sin(x_{i+1})$ and $\sigma_{ij} = 0.5 \delta_{ij} \tanh\fbr{(t - 0.5)^2 + d^{-1}\sum_{k=1}^d x_k^2}$ for $i, j = 1, 2, \cdots, d$. 
    In contrast to HJB-1a and HJB-1b, this equation features a variable diffusion coefficient, resulting in a solution whose smoothness varies with $t$ and $x$.
\end{itemize}

\begin{table}[t]
    \centering
    \caption{
    Numerical results of \Cref{alg_amnet} for HJB-1a and HJB-1b from \cref{sec_hjb} with $d = 10^4$.
    The convergence histories are shown in \Cref{fig_HP9abcd}.
    }\label{tab_ch_hjb}
    \resizebox{\textwidth}{!}{
    \begin{tabular}{l l l l l l} 
    \toprule
    Equation & Mean of $\mathrm{RE}_1$ & SD of $\mathrm{RE}_1$ & Mean of $\mathrm{RE}_{\infty}$ & SD of $\mathrm{RE}_{\infty}$ & RT (s) \\ [0.5ex] 
    \midrule
    HJB-1a & 8.76E{-3} & 8.09E{-4} & 3.76E{-2} & 5.74E{-3} & 2709 \\ 
    HJB-1b & 1.73E{-2} & 8.07E{-4} & 4.74E{-2} & 1.94E{-3} & 2707 \\
    \bottomrule
    \end{tabular}
    }
\end{table}

\begin{figure}[t]
    \centering
    \resizebox{1.0\textwidth}{!}{
    \begin{tabular}{@{}c@{\hspace{1mm}}cccc@{}}
        & \textbf{$s \mapsto v(0, s\B{1}_d)$} & \textbf{$s \mapsto v(0, \B{l}(s))$} & \textbf{$t \mapsto \mathrm{RE}(t)$} & \textbf{RE vs Iter.} \\
        
        \adjustbox{valign=m}{\rotatebox[origin=c]{90}{\textbf{HJB-1a}}} & 
        \adjustbox{valign=m}{\includegraphics[width=0.33\textwidth]{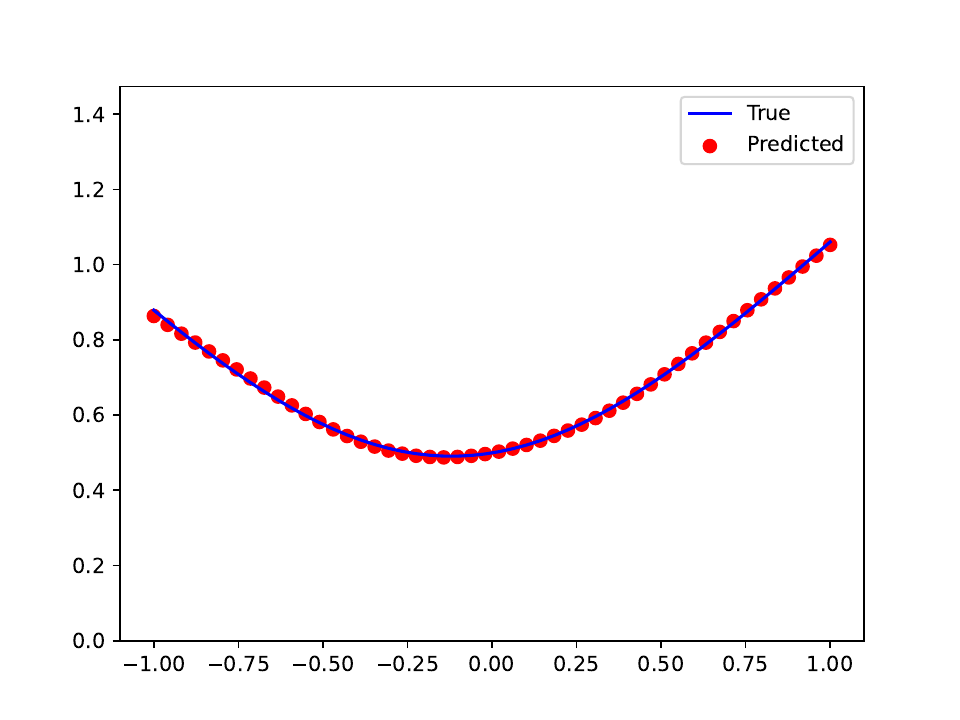}} &
        \adjustbox{valign=m}{\includegraphics[width=0.33\textwidth]{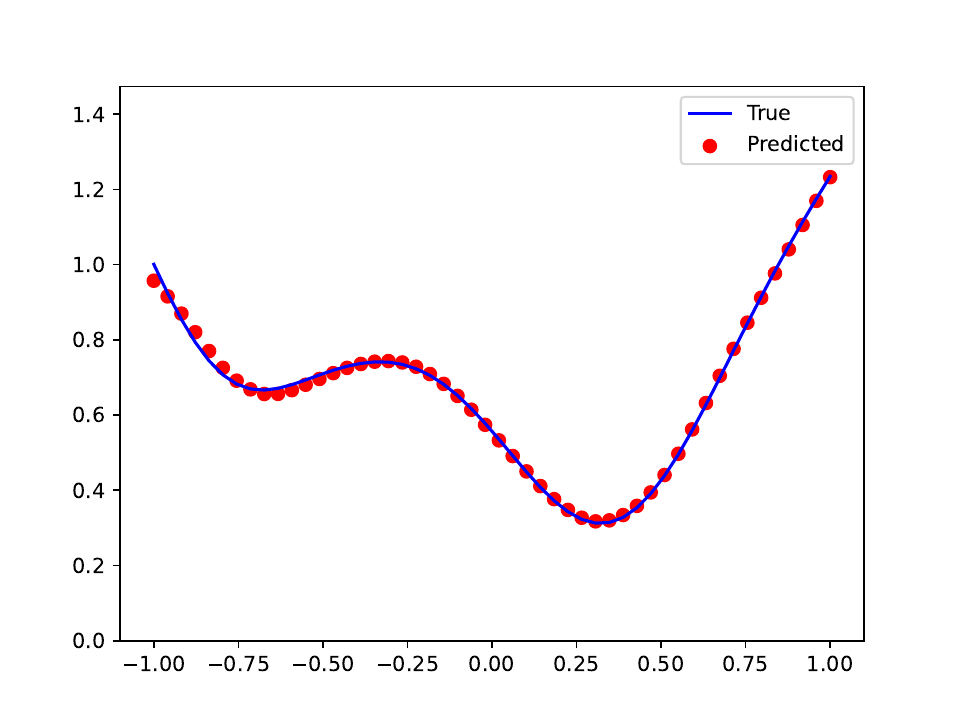}} &
        \adjustbox{valign=m}{\includegraphics[width=0.3\textwidth]{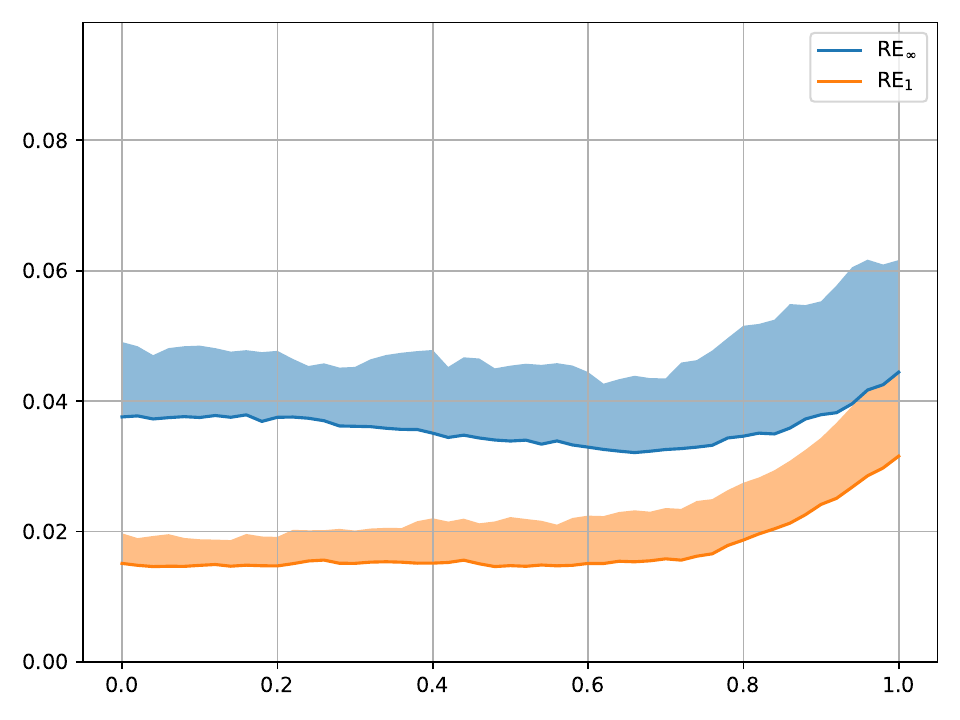}} &
        \adjustbox{valign=m}{\includegraphics[width=0.3\textwidth]{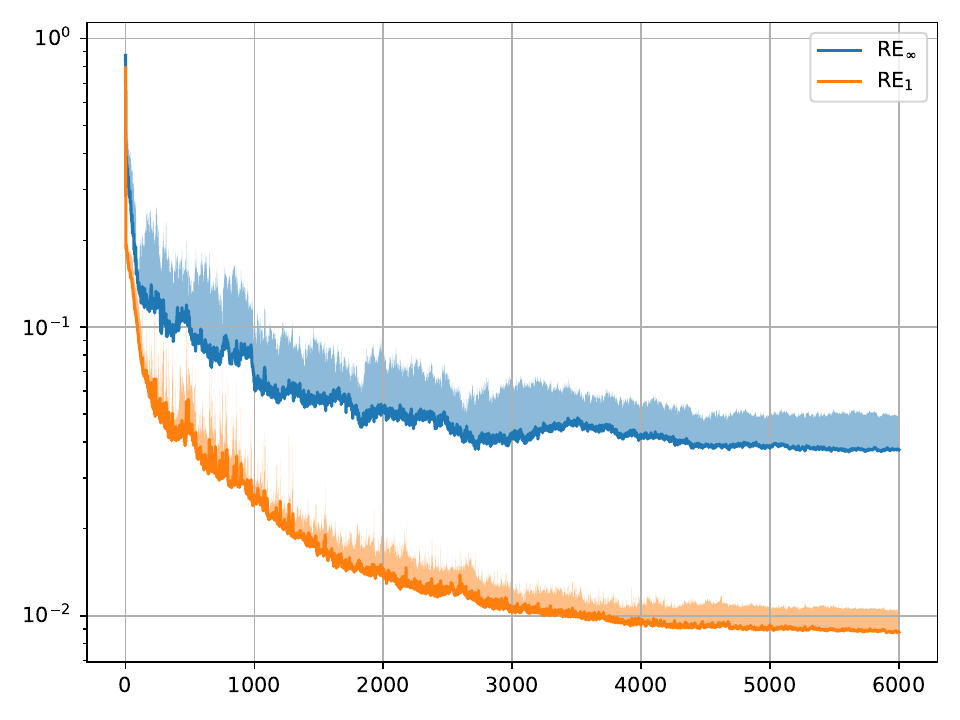}} \\
        \addlinespace[3pt]
        
        \adjustbox{valign=m}{\rotatebox[origin=c]{90}{\textbf{HJB-1b}}} & 
        \adjustbox{valign=m}{\includegraphics[width=0.33\textwidth]{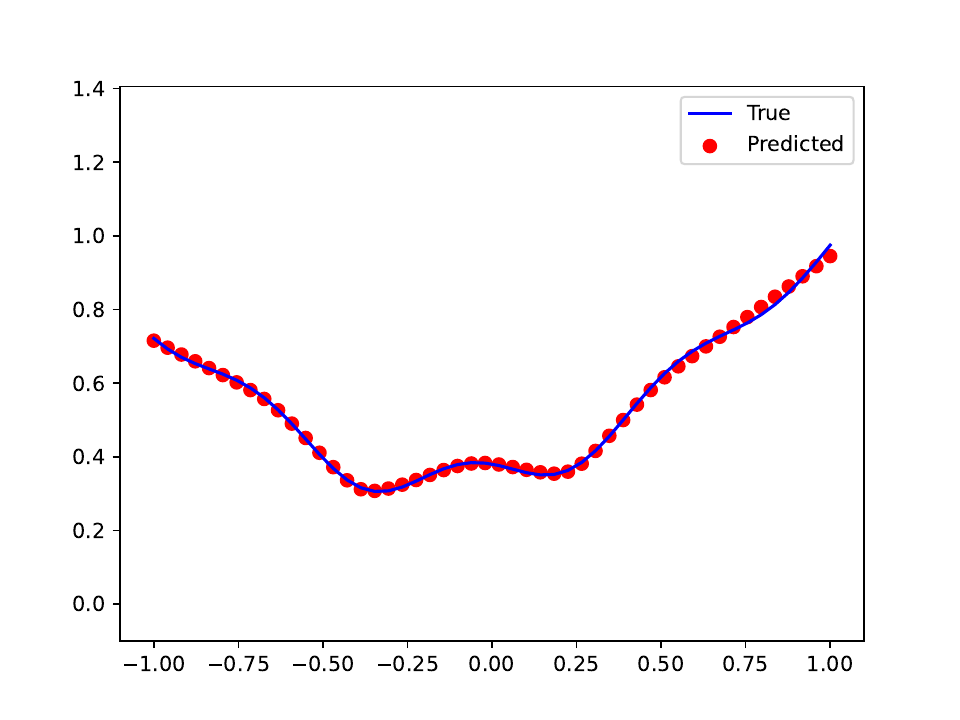}} &
        \adjustbox{valign=m}{\includegraphics[width=0.33\textwidth]{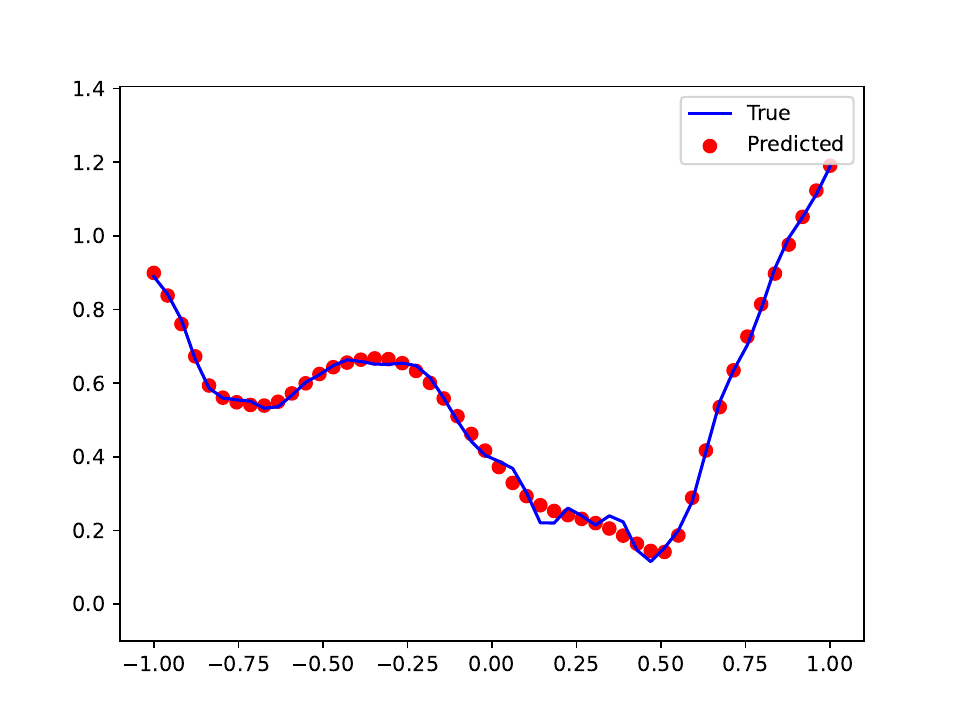}} &
        \adjustbox{valign=m}{\includegraphics[width=0.3\textwidth]{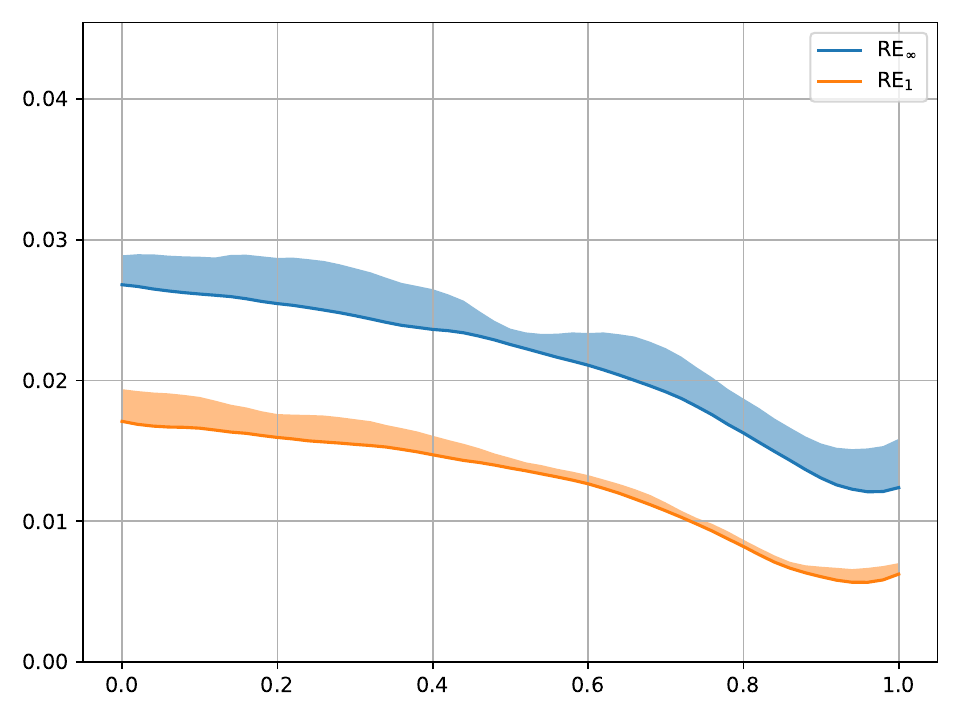}} &
        \adjustbox{valign=m}{\includegraphics[width=0.3\textwidth]{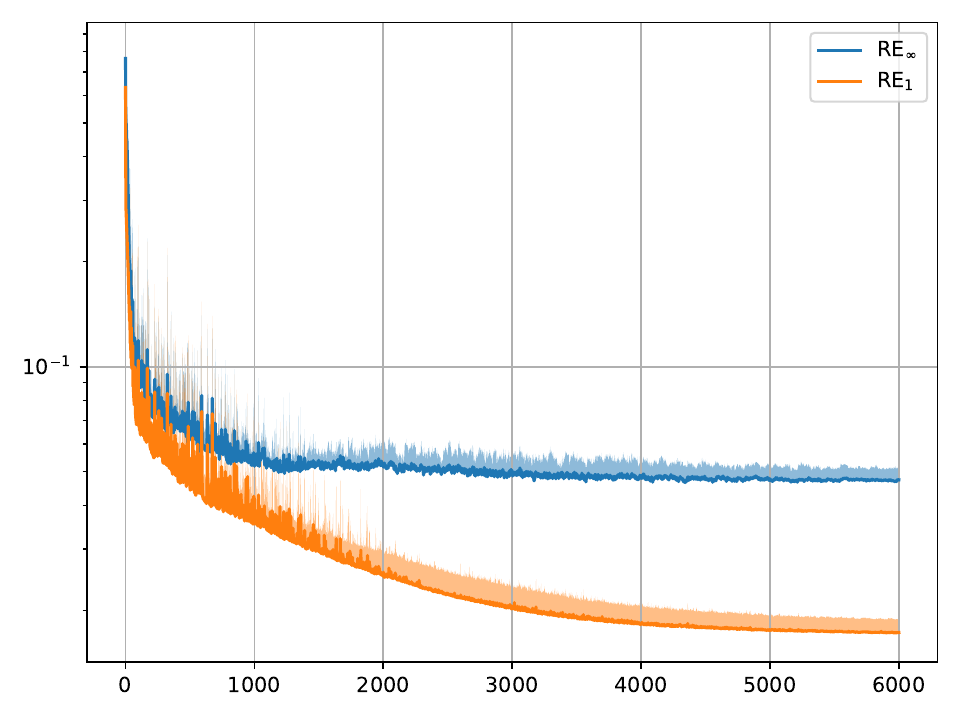}} \\
        \addlinespace[3pt]
    \end{tabular}
    }
    \caption{
    Numerical results of \Cref{alg_amnet} for HJB-1a and -1b from \cref{sec_hjb} with $d = 10^4$.
    The shaded regions in the third and fourth columns represent the mean $\pm 2 \times \text{SD}$ of the REs over 5 independent runs.
    The REs in the fourth column are computed at $t = 0$.
    The mean and the SD of the REs, and the RT, are given in \Cref{tab_ch_hjb}. 
    }\label{fig_HP9abcd}
\end{figure}

For HJB-1a and HJB-1b, we use standard fully connected DNNs with a width of $W = d + 10$.
The corresponding numerical results are shown in \Cref{fig_HP9abcd} and \Cref{tab_ch_hjb}.
As shown in the first two columns of \Cref{fig_HP9abcd}, the numerical solutions closely match the exact solutions. 
The third column of \Cref{fig_HP9abcd} demonstrates that our method maintains high accuracy along the regions explored by the sample paths generated via \eqref{eq_defhatXpil0_hjb}.  
The RTs in \Cref{tab_ch_hjb} further confirm that our method is highly efficient even for $d=10^4$. 

\begin{figure}[t]
    \centering
    \subfloat[$s \mapsto v(0, s \B{1}_d)$]{\includegraphics[width=0.3\textwidth]{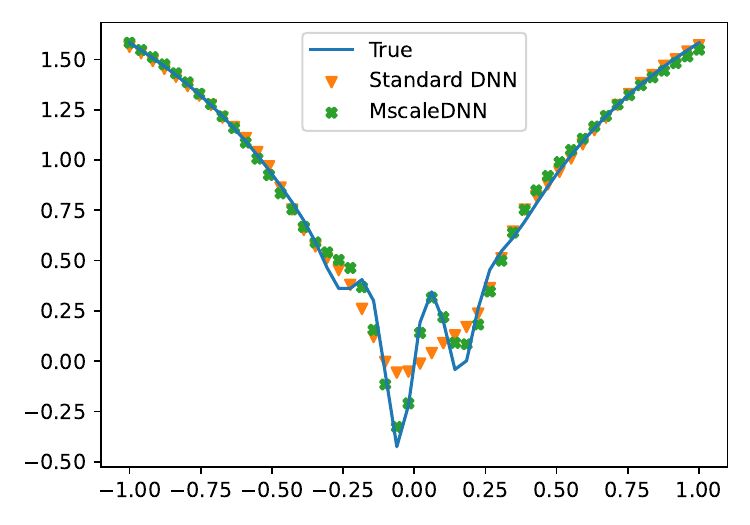}}
    \subfloat[$s \mapsto v(0, \B{l}(s))$]{\includegraphics[width=0.3\textwidth]{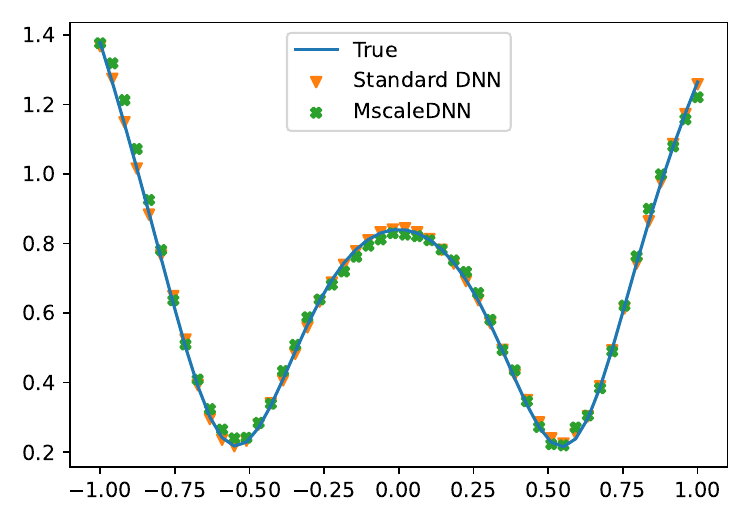}}
    \subfloat[RE$_{\infty}$ vs Iter.]{\includegraphics[width=0.3\textwidth]{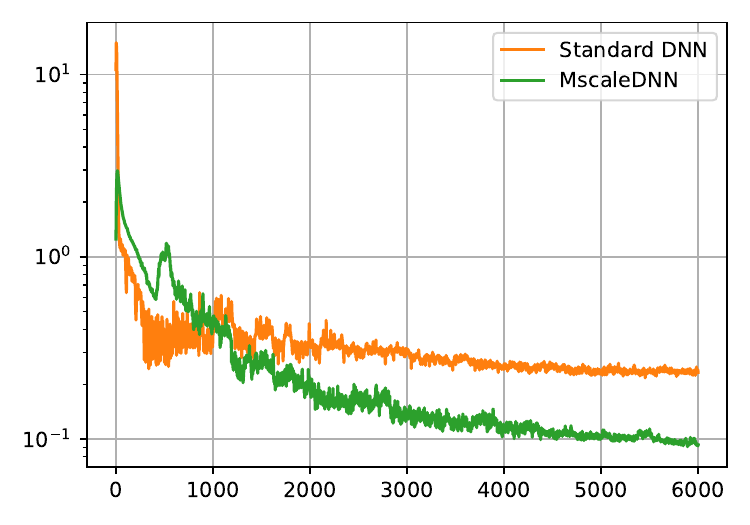}}
    \caption{
    Numerical results of \Cref{alg_amnet} for HJB-2 from \cref{sec_hjb} with $d = 10^4$.
    The RE$_{\infty}$ in (c) is evaluated at $t = 0$.
    The RT is 2953 seconds for Standard DNN and 2026 seconds for MscaleDNN, respectively.
    }\label{fig_hc2}
\end{figure}

For HJB-2, as shown in \Cref{fig_hc2} (a), the exact solution $s \mapsto v(0, s \B{1}_d)$ exhibits nonuniform smoothness: it is oscillatory near $s = 0$ and becomes smoother as $s$ moves away from zero. 
To improve function approximation for HJB-2, in addition to standard deep neural networks (DNNs) with width $W = d + 10$, we use the multi-scale DNNs (Mscale DNNs) \cite{Liu2020Multi} for the functions $u$ and $v$. 
Specifically, $u_{\alpha}$ is constructed as the MscaleDNN-2 architecture proposed in \cite[section 3.2]{Liu2020Multi}: 
\begin{equation}\label{eq_mscaleDNN}
    u_{\alpha}(t, x) = \frac{1}{6} \sum_{i=1}^6 u_{i, \alpha_i}(t, i x), \quad \alpha = \br{\alpha_1, \alpha_2, \cdots, \alpha_6},
\end{equation}
where each sub-network $u_{i, \alpha_i}$ is a fully connected DNN parameterized by $\alpha_i$ with width $\lceil W / 6\rceil$ and $W = d + 10$. 
The network $v_{\theta}$ is constructed in a similar manner.

The numerical results for HJB-2 with $d=10^4$ are presented in \Cref{fig_hc2}.
As shown in \Cref{fig_hc2} (a), the use of MscaleDNN enables our method to accurately capture the solution profile, even in regions with nonuniform smoothness.
Furthermore, using MscaleDNNs not only improves the relative error but also reduces the RT, 
as the adopted MscaleDNNs have fewer parameters compared to their standard versions.

\subsection{Comparison of sampling methods}\label{sec_ablation}

In this example, we examine how different spatial sampling strategies affect the performance of our method.
The benchmark problem is the convection-diffusion equation~\eqref{eq_linearPDE}, which features a steep gradient in its solution to highlight the importance of spatial sampling.
Since the effectiveness of our method in high-dimensional settings has been demonstrated in \cref{sec_cde}, this subsection focuses on the low-dimensional case with $d = 1$ to provide clearer visualization of spatial sampling. 
We consider two cases, $c = 1$ and $c = 5$, which lead to different levels of gradient steepness around $(t, x) = (0, 0)$.
Two strategies for generating spatial sampling points at time $t_n$ are compared:
\begin{itemize}
    \item \textbf{Dynamics}: The spatial samples at $t = t_n$ are generated by simulating the paths $X_n^m$ according to \eqref{eq_pathX}, where the initial points $X_0^m$ are drawn uniformly from $[-1, 1]$. This strategy exploits the dynamics of the PDE~\eqref{eq_linearPDE} and is theoretically supported by our analysis in \cref{sec_pathsamp}.

    \item \textbf{Plain}: The spatial samples at $t = t_n$ are given by $X_n^m = X_0^m + B_{t_n}^m$ for $m=1, 2, \cdots, M$, where $X_0^m$ and $B_{t_n}^m$ are i.i.d. samples from the uniform distribution on $[-1, 1]$ and the normal distribution $\mathrm{N}(0, t_n)$, respectively.
    This sampling method serves as a baseline for comparison.
\end{itemize}

\begin{figure}[t]
    \resizebox{\textwidth}{!}{
    \begin{tabular}{@{}c@{\hspace{1mm}}ccc@{}}
    & \textbf{Samples of $t \mapsto X_t$} & \textbf{True $t \mapsto v(t, X_t)$} & \textbf{True\&Predicted $x \mapsto v(0, x)$} \\
    
    \adjustbox{valign=m}{\rotatebox[origin=c]{90}{\textbf{Dynamics, $c=1$}}} &   \adjustbox{valign=m, raise=-0.03\height}{\includegraphics[trim=0mm 5mm 5mm 0mm,clip,width=0.42\textwidth]{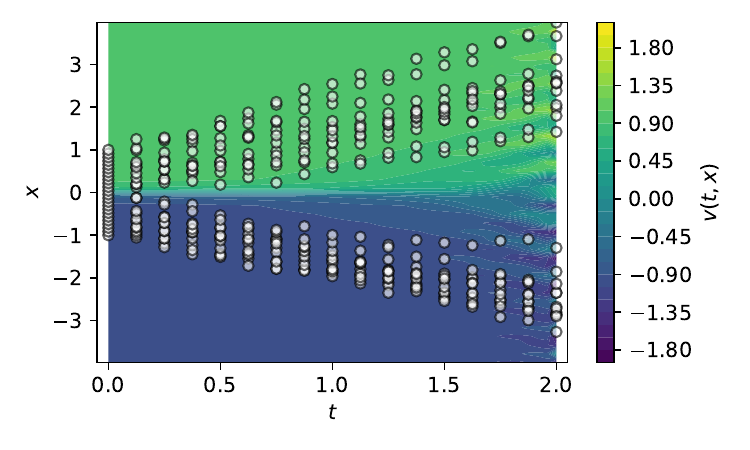}} &
    \adjustbox{valign=m}{\includegraphics[trim=10mm 5mm 15mm 10mm,clip,width=0.28\textwidth]{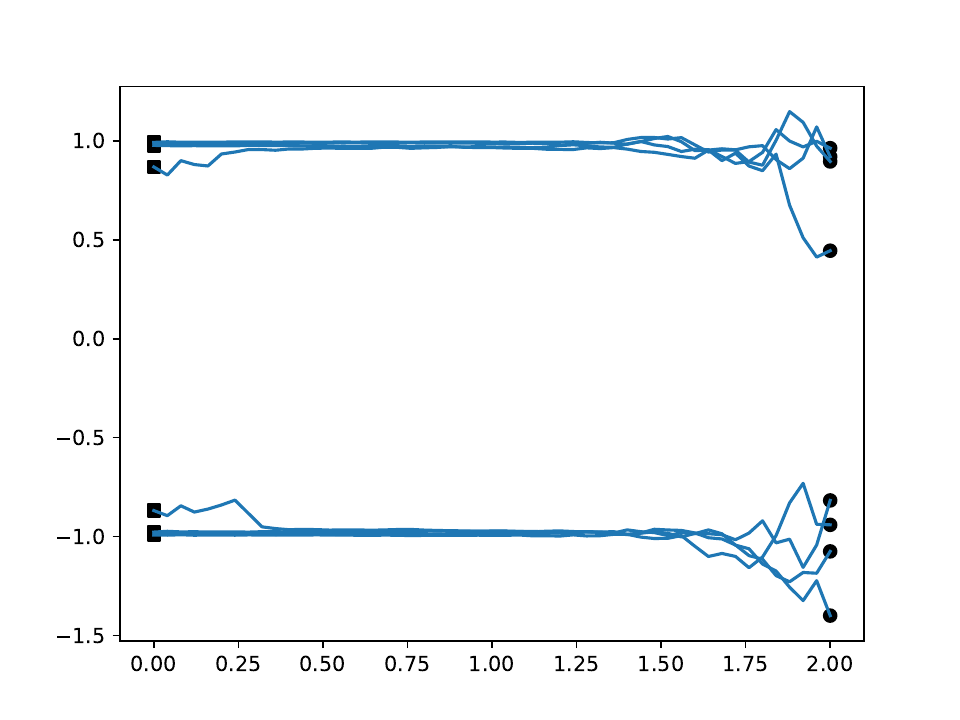}} &
    \adjustbox{valign=m}{\includegraphics[trim=10mm 5mm 15mm 10mm,clip,width=0.28\textwidth]{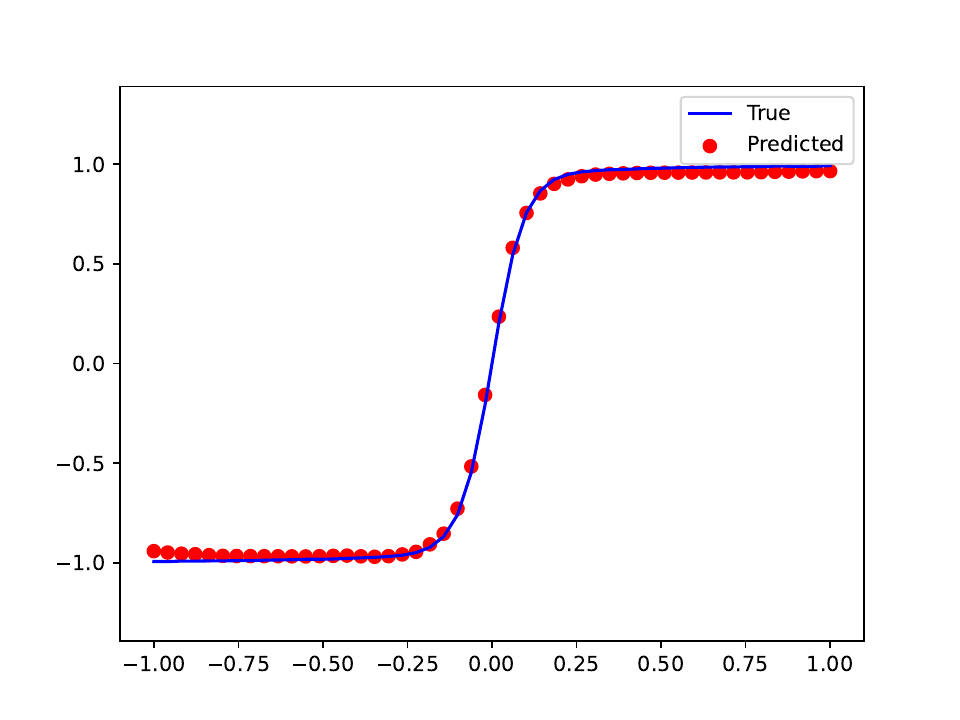}} \\
    \addlinespace[3pt]
    
    \adjustbox{valign=m}{\rotatebox[origin=c]{90}{\textbf{Dynamics, $c=5$}}} &   \adjustbox{valign=m, raise=-0.03\height}{\includegraphics[trim=0mm 5mm 5mm 0mm,clip,width=0.42\textwidth]{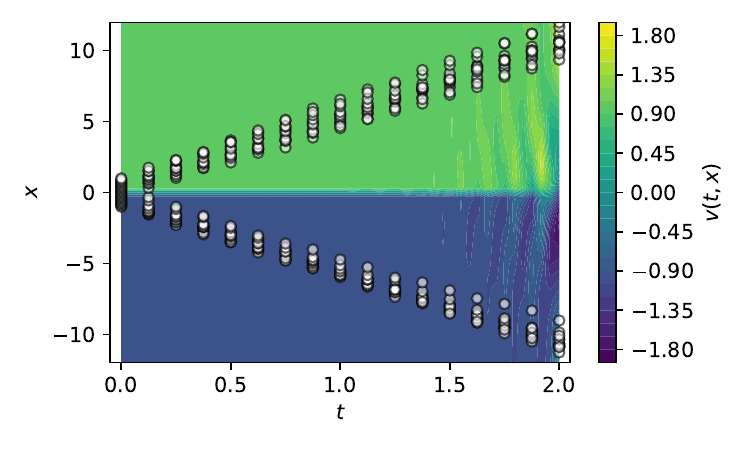}} &
    \adjustbox{valign=m}{\includegraphics[trim=10mm 5mm 15mm 10mm,clip,width=0.28\textwidth]{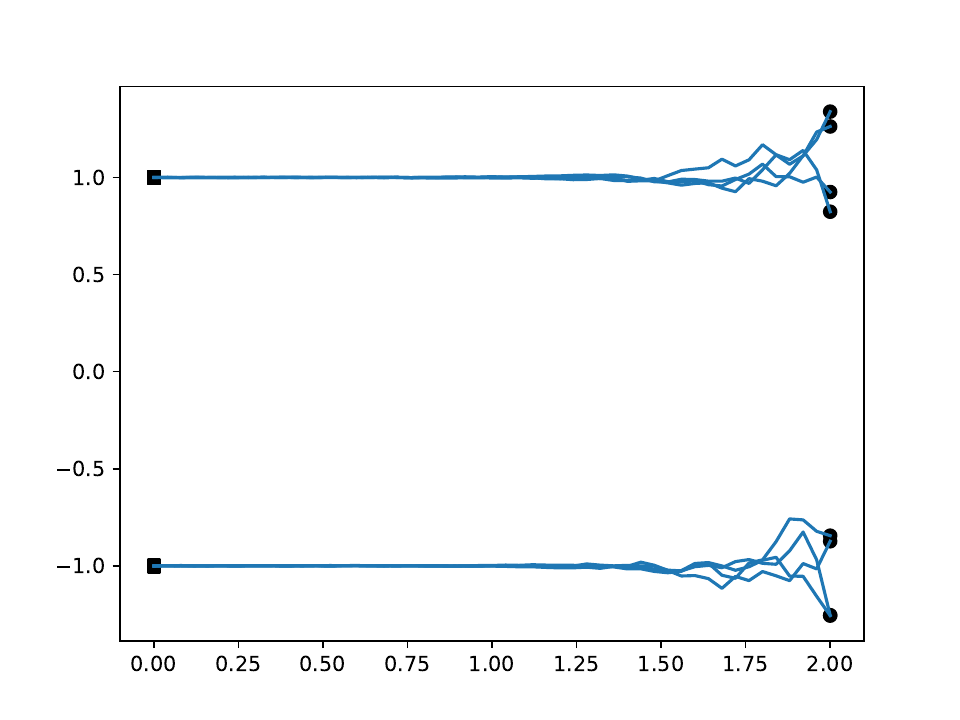}} &
    \adjustbox{valign=m}{\includegraphics[trim=10mm 5mm 15mm 10mm,clip,width=0.28\textwidth]{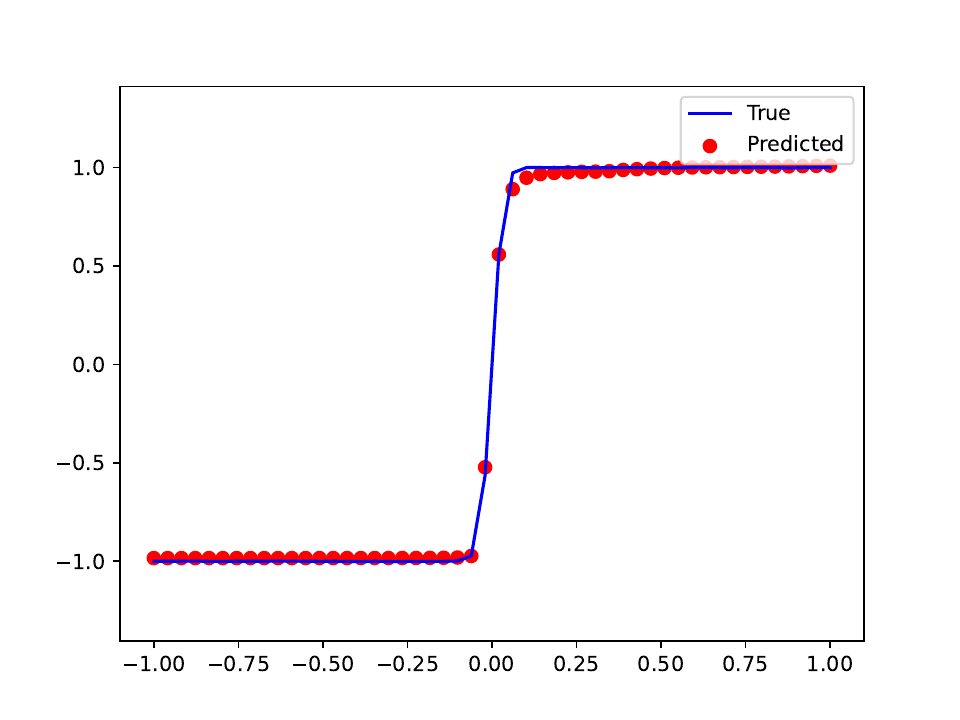}} \\
    \addlinespace[3pt]
    
    \adjustbox{valign=m}{\rotatebox[origin=c]{90}{\textbf{Plain, $c=1$}}} &   \adjustbox{valign=m, raise=-0.03\height}{\includegraphics[trim=0mm 5mm 5mm 0mm,clip,width=0.42\textwidth]{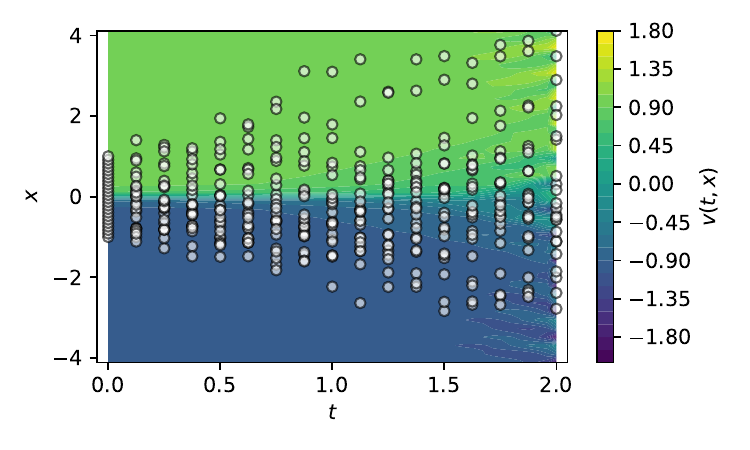}} &
    \adjustbox{valign=m}{\includegraphics[trim=10mm 5mm 15mm 10mm,clip,width=0.28\textwidth]{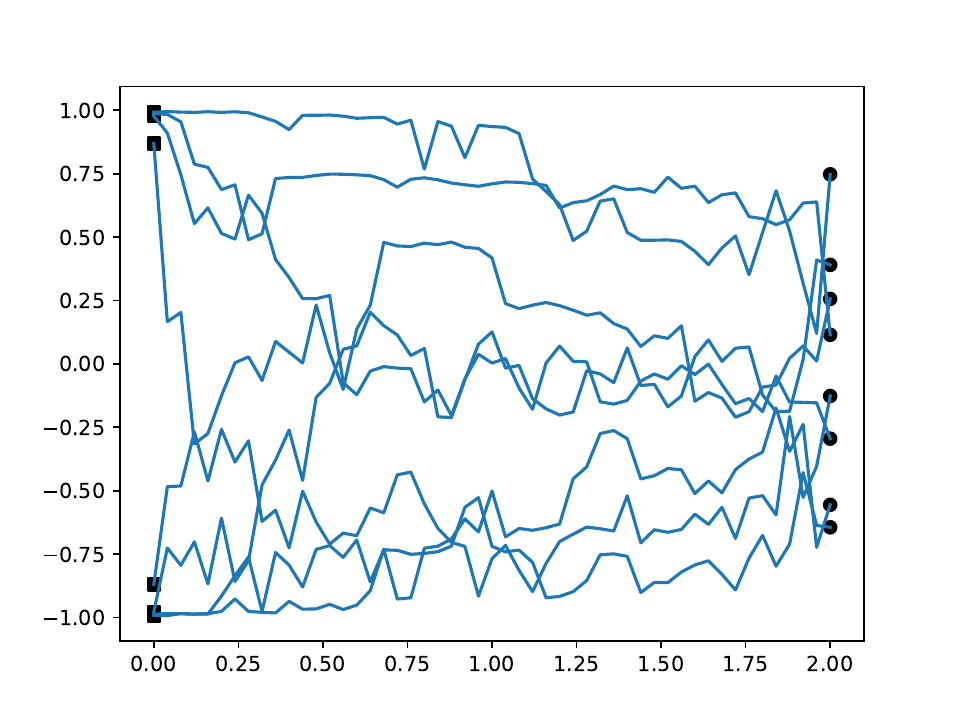}} &
    \adjustbox{valign=m}{\includegraphics[trim=10mm 5mm 15mm 10mm,clip,width=0.28\textwidth]{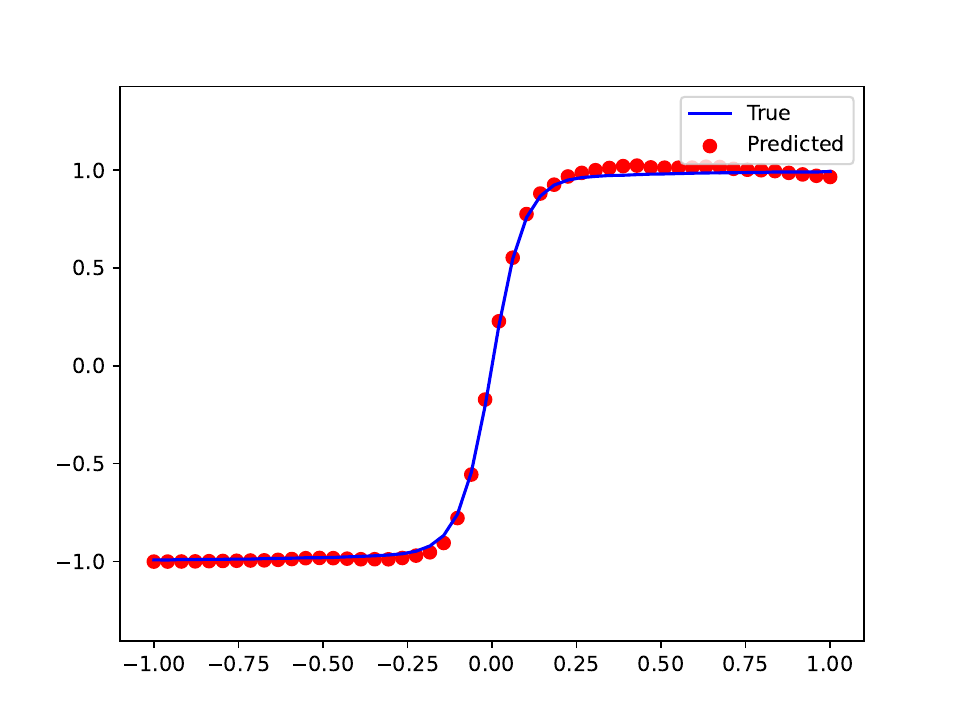}} \\
    \addlinespace[3pt]

    \adjustbox{valign=m}{\rotatebox[origin=c]{90}{\textbf{Plain, $c=5$}}} &   \adjustbox{valign=m, raise=-0.03\height}{\includegraphics[trim=0mm 5mm 5mm 0mm,clip,width=0.42\textwidth]{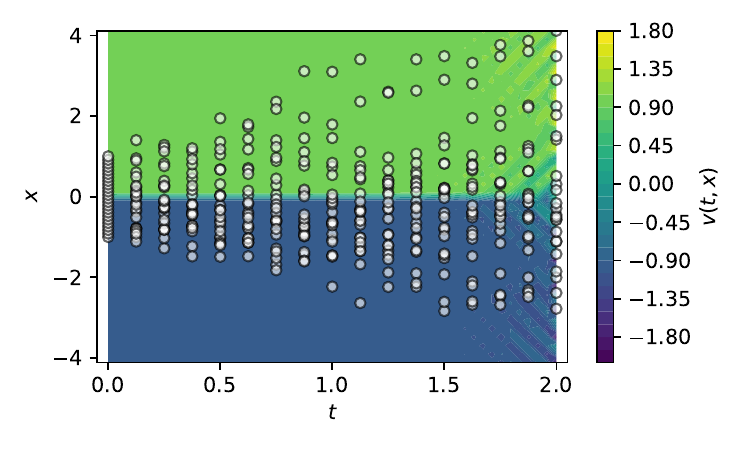}} &
    \adjustbox{valign=m}{\includegraphics[trim=10mm 5mm 15mm 10mm,clip,width=0.28\textwidth]{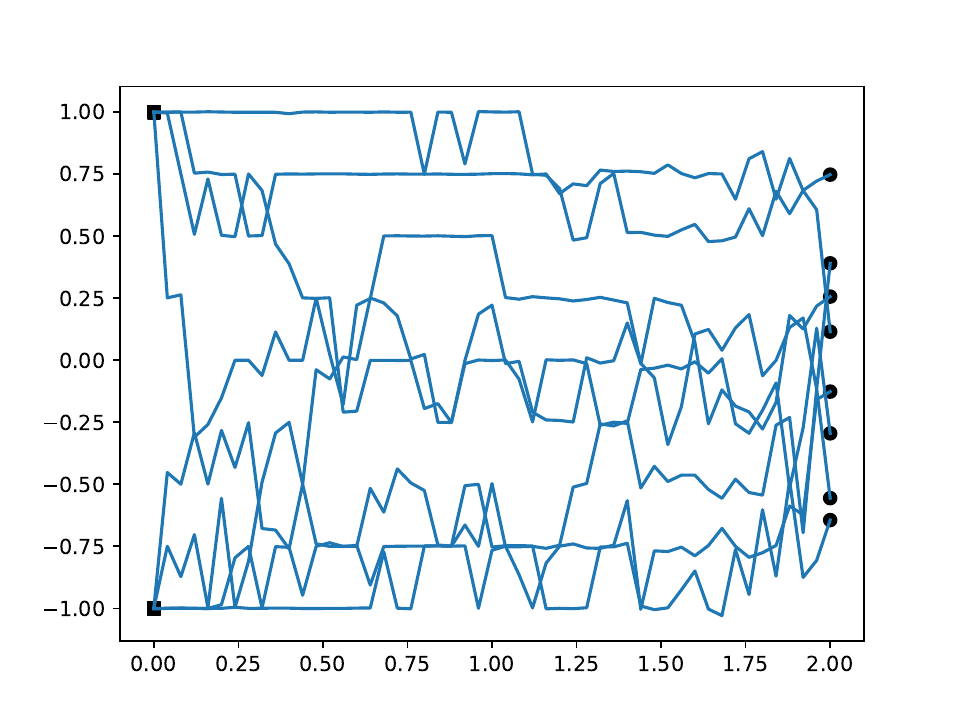}} &
    \adjustbox{valign=m}{\includegraphics[trim=10mm 5mm 15mm 10mm,clip,width=0.28\textwidth]{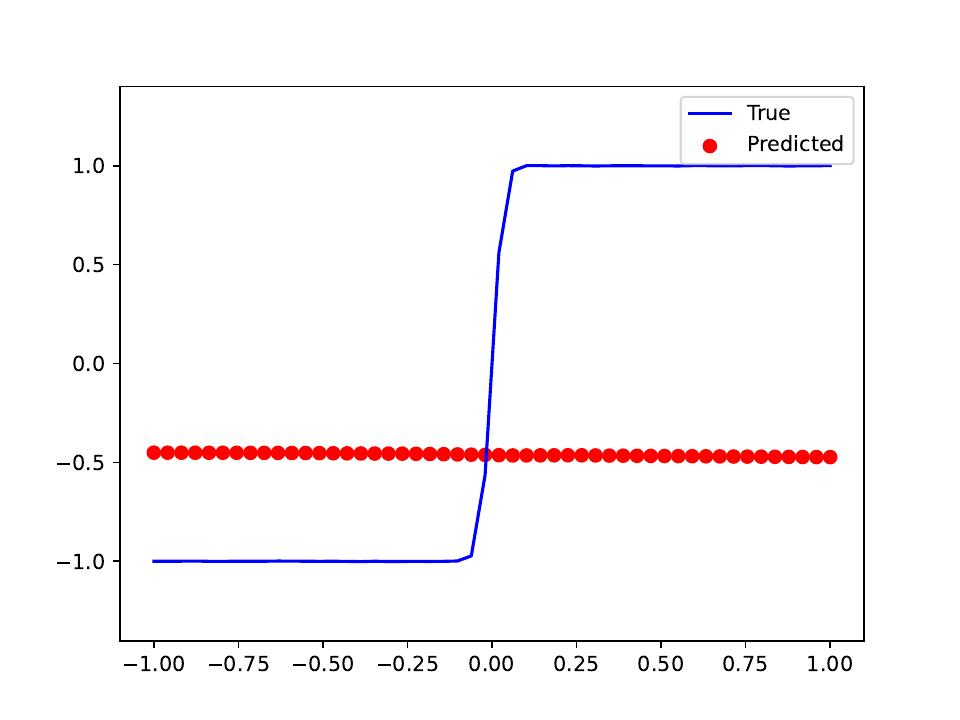}} \\
    \addlinespace[3pt]
    \end{tabular}
}
    \caption{ 
    Numerical results of \Cref{alg_pde} for the linear parabolic PDE~\eqref{eq_linearPDE} with $d = 1$.
    In the first column, the scatter plots show the spatial samples, where the color indicates the exact solution $v(t, x)$.
    The second column shows the landscape of the exact solution $v(t, x)$ along the sampling path $t \mapsto X_t$.
    The third column presents the numerical results for $x \mapsto v(t, x)$ at $t = 0$.
    }\label{fig_linearPDE}
\end{figure}

\begin{figure}
    \centering
    \includegraphics[width=0.65\textwidth]{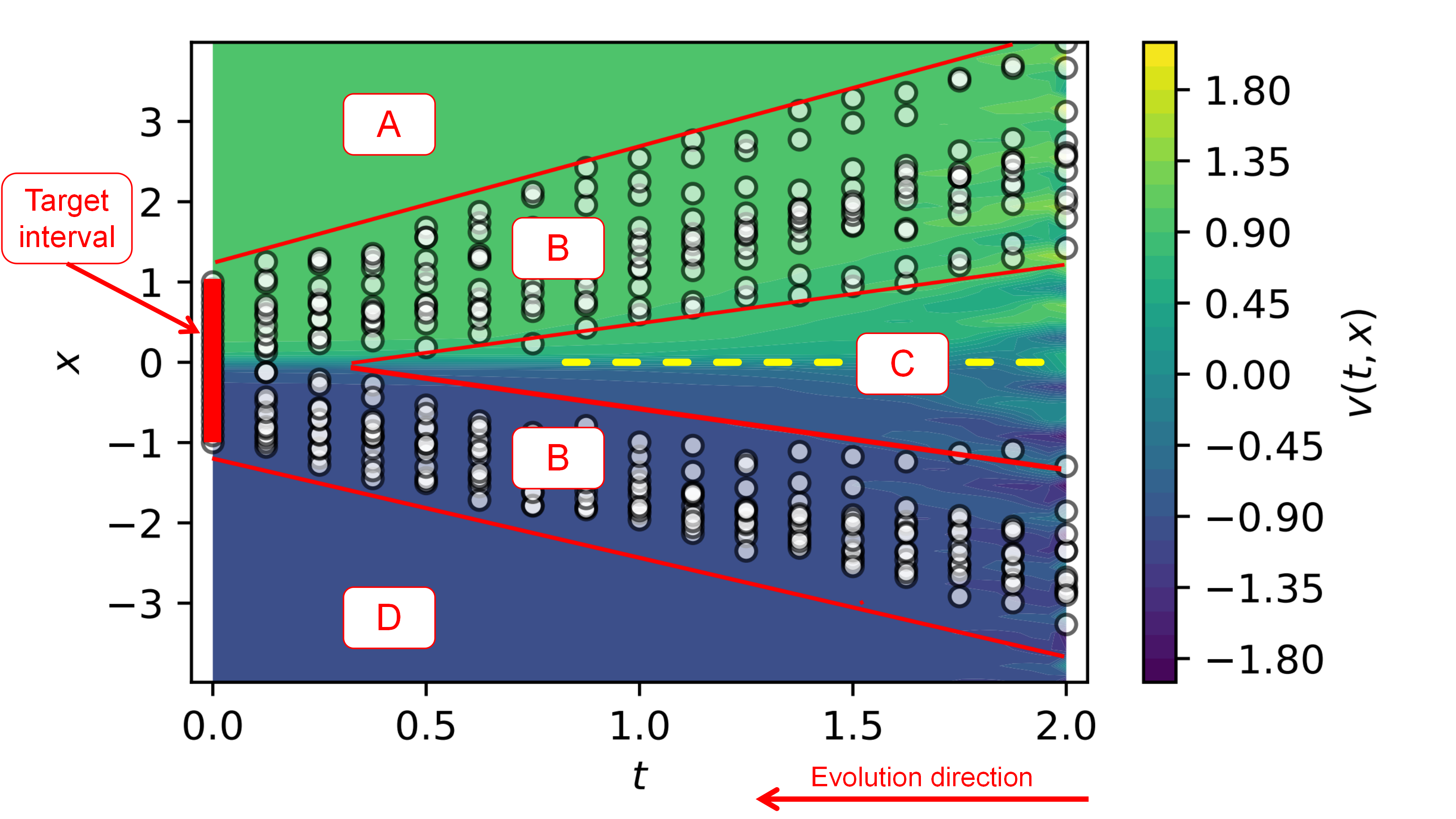}
    \caption{
    Region partition for the subfigure in the first row, first column of \Cref{fig_linearPDE}.
    The red solid line partitions the domain into four regions, labeled A, B, C, and D.
    The yellow dashed lines in region~C indicate the line $x = 0$, where the drift function $\mu$ in \eqref{eq_defmu_cde} transitions from negative to positive. 
    The scatter plots denote the spatial samples generated by the dynamics sampling for the linear parabolic PDE~\eqref{eq_linearPDE} with $d = 1$ and $c = 1$, where the color represents the exact solution $v(t, x)$.
    The red solid line labeled ``target interval'' indicates the spatial interval over which the third column of \Cref{fig_linearPDE} plots the solution profile. 
    }\label{fig_region_part}
\end{figure}

The corresponding numerical results are shown in \Cref{fig_linearPDE}.
In the first column of \Cref{fig_linearPDE}, the spatial samples generated by the dynamics sampling (first two rows) tend to avoid the steep gradient (i.e., the area near $x = 0$), especially for larger $c$ values. 
In contrast, the plain sampling (last two rows) produces more samples in the steep gradient, since the Gaussian distribution of $B_{t_n}^m$ is centered at $x = 0$.
Interestingly, as shown in the third column, the dynamics sampling yields better results than the plain sampling, even though the latter provides more samples in the steep gradient. 

To provide some intuition for this phenomenon,\label{txt_samp_insights}
we focus on the subfigure in the first row, first column of \Cref{fig_linearPDE} and partition the domain into four regions A, B, C, and D, as illustrated in \Cref{fig_region_part}.
Recalling the sign of the diffusion term in \eqref{eq_linearPDE}, it describes a convection-diffusion process evolving backward in time, from $t = 2$ to $t = 0$, as indicated by the red arrow in the bottom-right corner of \Cref{fig_region_part}.
As $t$ decreases, two mechanisms govern the evolution of the solution profile $v(t, x)$:
i) the drift term $\mu$ given in \eqref{eq_defmu_cde} transports the solution and tends to compress it toward $x=0$, 
thereby sharpening the gradient near $x=0$;
ii) the diffusion term smooths the solution profile, 
and the steep gradient region around $x=0$ is particularly sensitive to diffusion effects.

These drift-diffusion mechanisms also govern error propagation through the error equation~\eqref{eq_error_eq}. 
In particular, the residual $\hat{r}$ acts as a source term that generates the error $\epsilon$, which is then transported backward in time by the drift and smoothed by the diffusion.
This leads to distinct error-propagation behaviors across the four regions shown in \Cref{fig_region_part}.
In regions A and D, points are too far from the target interval $x \in [-1,1]$ given the remaining time $t$, so errors can hardly be transported into the target interval by $t = 0$. 
In region C, points lie too close to the target interval: the drift compresses the errors into the narrow region around the yellow dashed lines, 
producing high-frequency components in the error profile.
These high-frequency components are highly sensitive to diffusion and are rapidly smoothed out before they can reach $t = 0$.
Consequently, errors generated in regions A, C, and D have little influence on $v(0, x)$ over the target interval.
In contrast, region B lies at an appropriate distance from the target interval, allowing errors to be transported into it without being subject to excessive compression and diffusion.
Therefore, residuals in region B are the most influential.

We now illustrate why the dynamics sampling is effective.
Overall, as shown in \Cref{fig_linearPDE,fig_region_part}, the dynamics sampling places fewer samples in regions A, C, and D, and instead concentrates on the more influential region B. 
This behavior aligns with the intuition above and provides further justification for the error-residual relationship in \eqref{eq_v0x0}.
In contrast, the plain sampling tends to oversample region C: although this region contains the steep gradient, residuals 
located there have limited impact on the error at $t = 0$.
In terms of information propagation, the dynamics of how the terminal data at $t=T$ affects the final profile of the solution at $t=0$ is consistent with the interpretation of the Feynman-Kac formula for the parabolic equation as propagation of the terminal data through time-backward running infinitely many ``generalized characteristic lines'' (borrowing the characteristics concept of the hyperbolic PDE corresponding to a pure drift system).

Finally, we turn to the center column of \Cref{fig_linearPDE}, which shows the profile of the exact solution $v(t, x)$ at $t = 0$ along the sampling path $t \mapsto X_t$.
Notably, the solution is significantly smoother along the path generated by the dynamics sampling compared to the plain one.
This observation aligns with the method of characteristics.
Dynamics-based sampling may thus reduce the difficulty of function approximation in high-dimensional spaces, which is an interesting direction for future research.

\subsection{Comparison of PINNs, RS-PINNs, and strong/weak-form DRDM}\label{sec_compare}

We compare the proposed DRDM~\eqref{eq_minmax} against two baselines: PINNs \cite{Raissi2019Physics} and the hybrid version of RS-PINNs \cite{Hu2025Bias} (abbreviated as PINNs and RS-PINNs, respectively).
Both baselines have demonstrated strong performance on high-dimensional PDEs.
To highlight the efficiency gains achieved by the weak-form reformulation in \eqref{eq_minmax}, 
we also evaluate the strong-form DRDM, whose loss is defined in \eqref{eq_Lrdm}, with the algorithm described in \Cref{alg_pde_strf} in the appendix. 
Although the strong-form DRDM is derived from the random difference method, \Cref{alg_pde_strf} is essentially a variant of the shotgun method \cite{zhang2025shotgun}, differing mainly in its offline path sampling and in not using antithetic variates for variance reduction. 
Thus, we also refer to it as the shotgun method in the numerical results.
 
The benchmark problem is from \cite[subsection 5.4]{Hu2025Bias} and includes:
\begin{itemize}
    \item Allen-Cahn equation:
    \begin{equation}\label{eq_ac}
        \Delta v(x) + v(x) - v(x)^{3} = f(x), \quad x \in \mathbb{B}^d,
    \end{equation}
    \item Sine-Gordon equation:
    \begin{equation}\label{eq_sg}
        \Delta v(x) + \sin(v(x)) = f(x), \quad x \in \mathbb{B}^d. 
    \end{equation}
\end{itemize}
In both problems, the domain is the unit ball $\mathbb{B}^d := \{x \in \R^d : \abs{x} < 1\}$; 
the boundary condition is $v(x) = 0$ for all $x \in \partial \mathbb{B}^d$;
the source term $f$ is chosen so that the exact solution is
\begin{equation}
    v(x) = \br{1 - \abs{x}^2} \vbr{\Big}{\sum_{i=1}^{d-1} c_i \sin\vbr{\big}{x_i + \cos(x_{i+1}) + x_{i+1} \cos(x_i)}},
\end{equation}
where $c_1, c_2, \cdots, c_{d-1}$ are i.i.d. samples from $\mathrm{N}(0, 1)$.

To ensure a fair comparison and avoid discrepancies due to hyperparameter tuning,
we directly use the results reported in \cite[subsection 5.4.2]{Hu2025Bias} for PINNs and RS-PINNs,
and implement the strong- and weak-form DRDM to closely match their experimental setup.
Detailed parameter settings are given in \cref{sec_param_setting}.
The results are summarized in \Cref{tab_compare}.

\begin{table}[t]
\centering
\caption{
Relative $L_2$ errors, running times (RTs), and memory usage for PINNs \cite{Raissi2019Physics}, hybrid version of RS-PINNs \cite{Hu2025Bias}, strong- and weak-form DRDM across Allen-Cahn (AC) and Sine-Gordon (SG) equations under various high dimensions.
Results for PINNs and RS-PINNs are taken directly from \cite[Table 6]{Hu2025Bias}, with the units of RTs converted to minutes to facilitate reading.
}\label{tab_compare}
\resizebox{\textwidth}{!}{
\begin{tabular}{cclcccccc}
\toprule
Method & Equation & Metric & $d=10^2$ & $10^3$ & $5 \times 10^3$ & $10^4$ & $10^5$ \\
\midrule
\multirow{4}{*}{\begin{tabular}[c]{@{}c@{}}PINNs\\ \cite{Raissi2019Physics}\end{tabular}}
 & AC & Rel. $L_2$ error & 7.187E-3 & 5.617E-4 & 1.773E-3 & N.A. & N.A. \\
 & SG & Rel. $L_2$ error & 7.192E-3 & 5.642E-4 & 1.782E-3 & N.A. & N.A. \\
 & & RT (minute) & 3 & 285 & 1832.4 & N.A. & N.A. \\
 & & Memory (MB) & 1328 & 4425 & 56563 & $>$81252 & $>$81252 \\
\midrule
\multirow{4}{*}{\begin{tabular}[c]{@{}c@{}}RS-PINNs\\ \cite{Hu2025Bias}\end{tabular}}
 & AC & Rel. $L_2$ error & 7.923E-3 & 5.504E-4 & 1.802E-3 & 1.860E-3 & 2.192E-3 \\
 & SG & Rel. $L_2$ error & 7.835E-3 & 6.744E-4 & 1.795E-3 & 1.854E-3 & 2.176E-3 \\
 & & RT (minute) & 1.8 & 7.2 & 31.8 & 66 & 720 \\
 & & Memory (MB) & 1413 & 1815 & 3593 & 5789 & 45599 \\
 \midrule
 \multirow{4}{*}{\begin{tabular}[c]{@{}c@{}}Strong-form\\DRDM\\ (Shotgun \cite{zhang2025shotgun}) \end{tabular}}
 & AC & Rel. $L_2$ error & 1.311E-2 & 6.171E-3 & 4.714E-3 & 3.112E-3 & 7.423E-4 \\
 & SG & Rel. $L_2$ error & 2.009E-2 & 6.191E-3 & 4.715E-3 & 3.112E-3 & 2.232E-2 \\
 & & RT (minute)        &  1.04      & 1.11      & 1.81      & 2.70      & 48.64  \\
 & & Memory (MB)        &  75        & 254       & 1204      & 2389      & 23739  \\
 \midrule
 \multirow{4}{*}{\begin{tabular}[c]{@{}c@{}}Weak-form\\ DRDM\\ (ours) \end{tabular}}
 & AC & Rel. $L_2$ error & 3.963E-2 & 6.211E-3 & 4.699E-3 & 3.118E-3 & 7.381E-4 \\
 & SG & Rel. $L_2$ error & 4.172E-2 & 6.306E-3 & 4.701E-3 & 3.118E-3 & 4.769E-3 \\
 & & RT (minute)        & 1.22      & 1.26     & 1.43     & 1.59     & 13.26  \\
 & & Memory (MB)        & 57        & 185      & 859      & 1698     & 16829  \\
\bottomrule
\end{tabular}
}
\end{table}

We make the following observations from \Cref{tab_compare}:
\begin{itemize}
\item 
As indicated by Rel. $L_2$ errors, 
both forms of DRDM generally yield lower accuracy than PINNs and RS-PINNs.
The reduced accuracy of DRDM is likely caused by the truncation error introduced by the random difference approximation in \eqref{eq_rdo}.
Both PINNs and RS-PINNs do not suffer from this error. 

\item 
In terms of RT and memory usage, both forms of DRDM provide significant advantages over PINNs, and their advantages become more pronounced as the dimension $d$ increases.
As reported in \cite{Hu2025Bias}, PINNs run out of memory for $d \geq 10^4$, due to the costly automatic differentiation for second-order derivatives.
In contrast, DRDM avoids this bottleneck and successfully handles $d = 10^5$ with favorable RT and memory usage.

\item
Across all tested dimensions, the strong-form DRDM requires less RT and memory than RS-PINNs, with the advantage in RT becoming more pronounced at higher dimensions.
This is somewhat surprising, as both methods share a similar spirit of employing Monte Carlo sampling to approximate the differential operator. 
We suspect the difference stems from the practical numerical implementations of the two methods. 

\item 
In most cases, the two forms of DRDM achieve similar accuracy. 
The primary differences are in runtime and memory usage. 
For $d \leq 10^3$, the strong form is slightly more efficient, since the weak form requires additional costs from gradient ascent steps.
As $d$ grows, the weak form gradually shows its advantage, which stems from requiring fewer evaluations of $v_{\theta}$ per gradient step.
In particular, in the strong-form loss~\eqref{eq_Lrdm}, each spatial sample $x$ needs $\ell$ independent draws of $\xi$ to approximate the expectation; in our experiments, $\ell = 128$.
In the weak-form loss~\eqref{eq_minmax}, the expectation is taken jointly over $(X_t, \xi)$, so each $X_t$ sample requires only one draw of $\xi$.
Thus, for a single gradient step with mini-batch size $n_b$, the strong form evaluates $v_{\theta}$ $n_b + n_b\ell$ times, whereas the weak form evaluates it only $2 n_b$ times; cf. \eqref{eq_defR}.

\item 
Compared with the results in \Cref{tab_RESQLP,tab_ch_hjb},
for the same dimensions $d = 10^4$ and $10^5$,
the weak-form DRDM achieves substantially shorter RTs in \Cref{tab_compare}.
This discrepancy is mainly due to the different network sizes used. 
As described in \cref{sec_param_setting}, the experiments reported in \Cref{tab_RESQLP,tab_ch_hjb} use 6-layer DNNs with widths of at least $10^4$ in order to achieve satisfactory accuracy for time-dependent quasilinear parabolic PDEs and HJB equations.
In contrast, the networks used in \Cref{tab_compare} are 4-layer DNNs with width 100, which have been shown in \cite{Hu2025Bias} to be sufficient for attaining satisfactory accuracy on the tested elliptic PDEs \eqref{eq_ac} and \eqref{eq_sg}.
\end{itemize}

\section{Conclusion}\label{sec_conclusion}

In this work, we proposed a deep random difference method (DRDM) for solving high-dimensional quasilinear parabolic PDEs and HJB equations. 
The DRDM approximates the convection-diffusion operator using only first-order differences, thereby avoiding explicit computation of Hessians and reducing reliance on automatic differentiation. 
By integrating DRDM within a Galerkin framework, we significantly reduce the variance of mini-batch estimation, enabling efficient and highly parallelizable implementations. 
The method is further extended to HJB equations.
We provide rigorous error estimates and demonstrate, through numerical experiments, that DRDM can efficiently and accurately solve quasilinear parabolic PDEs and HJB equations in dimensions up to $10^5$ and $10^4$, respectively.

The DRDM recovers existing martingale deep learning methods~\cite{cai2024socmartnet,cai2023deepmartnet,cai2024martnet} from a new perspective that is closer to the PINN framework and requires minimal reliance on stochastic calculus.
Our method offers several distinctive advantages:
\begin{itemize}
    \item The loss function does not require AD for any PDE derivatives, resulting in efficient computation and reduced memory usage.
    
    \item The loss function can be computed in parallel across both time and space, further enhancing efficiency.
    
    \item The sample distribution can be tailored to the dynamics of the convection-diffusion operator, enhancing the method's ability to handle problems with complex solution structures.
    
    \item The truncation error of the DRDM achieves first-order accuracy in time, which is higher than the half-order accuracy of methods based on the pathwise properties of the Euler-Maruyama scheme.
\end{itemize}
Future work will explore the application of DRDM to more challenging problems, including convection-diffusion equations with complex boundaries, committor functions, and more complex and non-convex SOCPs.

\appendix

\section{Error estimates}\label{sec_theory}

In this section, a theoretical analysis is presented for the RDM formulation \eqref{eq_ER0} applied to the quasilinear parabolic PDE~\eqref{eq_pde}.
We first derive a bound on the local truncation error, then prove zero-stability, and finally combine these results to obtain a global error estimate.
A discussion of the global error estimate is provided at the end of this section.

Some frequently used notations are introduced as follows.
\label{txt_notation}
Let $\abs{\cdot}$ denote the Euclidean norm on $\R^d$ and the Frobenius norm on $\R^{d\times q}$. Specifically, for any $x \in \R^d$, $\abs{x} := (\sum_{i=1}^d x_i^2)^{1/2}$, and for any $A \in \R^{d\times q}$, $\abs{A} := (\sum_{i=1}^d \sum_{j=1}^q A_{ij}^2)^{1/2}$, so that $\abs{Ax} \leq \abs{A}\,\abs{x}$.
Let $C^{p,q}$ denote the set of functions that are $p$-times continuously differentiable in $t \in [0, T)$ and $q$-times continuously differentiable in $x \in \R^d$. 
Let $\mathcal{M}_{2, 4}$ denote the set of multi-indices $\alpha = (\alpha_0, \alpha_1, \cdots, \alpha_d) \in \{0, 1, 2, \cdots\}^{d+1}$ satisfying $\alpha_0 \in \{1, 2\}$ and $\sum_{i=1}^d \alpha_i \in \{1, 2, 3, 4\}$.
For any $\alpha \in \mathcal{M}_{2, 4}$, define the differential operator
\begin{equation}\label{eq_defDalp}
    D^{\alpha} := \partial_t^{\alpha_0} \partial_{x_1}^{\alpha_1} \partial_{x_2}^{\alpha_2} \cdots \partial_{x_d}^{\alpha_d}. 
\end{equation}
For a function $(t,x) \mapsto \phi(t,x)$, the notation $\partial \phi/\partial t$ refers to the partial derivative of $\phi$ with respect to $t$, viewed as a function of $(t,x)$. 
By contrast, $\partial_t$ denotes the partial-derivative operator with respect to $t$. 
An example illustrating the distinction is $\partial_t \phi(t^2, x) = (\partial \phi/\partial t)(t^2, x) \times 2t$. 
The notation $\partial \phi/ \partial_{x_i}$ and $\partial_{x_i}$ are understood similarly.

To quantify errors on the time partition \eqref{eq_timepart}, for any Borel-measurable function $\phi: [0, T] \times \R^d \to \R$, we define the $p$-th moment functional $\mathbb{M}_n^p[\phi]$ by
\begin{equation}\label{eq_defMtp}
    \mathbb{M}_n^p[\phi] := \int_{\R^d} \abs{\phi(t_n, x)}^{p}\, P_{X_n}(x)\, \di x, \quad p \geq 1,
\end{equation}
where $x \mapsto P_{X_n}(x)$ is the PDF of $X_n^m$, the discrete-time stochastic process defined in \eqref{eq_pathX}.
The superscript $m$ in $X_n^m$ indicates the $m$-th independent sample path (not an exponent).
We suppress the notational dependence on $m$ in $\mathbb{M}_n^p[\cdot]$ and $P_{X_n}(\cdot)$, since $\{X_n^m\}_{m=1}^M$ are identically distributed for all $m$, and these quantities do not depend on $m$.

\subsection{Local truncation error}

To analyze the local truncation error of the random difference approximation~\eqref{eq_rdo}, we introduce the following assumption.

\begin{assumption}\label{assum_polygrow}
The solution $v$ of the PDE~\eqref{eq_pde} satisfies $v \in C^{2,4}$. 
The approximator $\hat{v}: [0, T] \times \R^d \to \R$ of $v$ is Borel-measurable.
The functions $\mu$, $\sigma$, $v$ and $\hat{v}$ satisfy
\begin{align}
    \abs{\mu(t, x, v(t, x))} + \abs{\mu(t, x, \hat{v}(t, x))} &\leq C_{\mathrm{g}} \big(1 + \abs{x}^{p_{\mu}}\big), \label{eq_lgrow_mu} \\ 
    \abs{\sigma(t, x, v(t, x))} + \abs{\sigma(t, x, \hat{v}(t, x))} &\leq C_{\mathrm{g}} \big(1 + \abs{x}^{p_{\sigma}}\big), \label{eq_lgrow_sgm}\\
    \sum_{\alpha \in \mathcal{M}_{2, 4}} \abs{D^{\alpha} v(t, x)} &\leq C_{\mathrm{g}} \big(1 + \abs{x}^{p_v}\big), \label{eq_polygrow} 
\end{align}
for all $t \in [0, T]$ and $x \in \R^d$, where $C_{\mathrm{g}} \geq 1$ and $p_{\mu}, p_{\sigma}, p_v \geq 0$ are constants independent of $t$, $x$ and $\hat{v}$.
The random vector $\xi$ in \eqref{eq_rdo} satisfies \eqref{eq_condi_xi} and has finite componentwise moments of order $\bar{m}$:
\begin{equation}\label{eq_finite_moment}
    \E{\abs{\xi_i}^{\bar{m}}} < \infty, \quad \bar{m} := \max\{2 p_v, 8\}, \quad i = 1, 2, \cdots, q.
\end{equation}
\end{assumption}

In \Cref{assum_polygrow}, the requirement $v \in C^{2,4}$ is mild for parabolic PDEs; see \Cref{rmk_reg} for details.
The Borel measurability of $\hat{v}$ is imposed to ensure that the expectation in \eqref{eq_ER0} is well defined. 
The growth conditions \eqref{eq_lgrow_mu}--\eqref{eq_lgrow_sgm} are satisfied whenever $\mu(t, x, y)$, $\sigma(t, x, y)$, $v(t, x)$, and $\hat{v}(t, x)$ grow at most polynomially in $x \in \R^d$ and $y \in \R$.
The polynomial-growth bounds for $v$ and its derivatives in \eqref{eq_polygrow} can be verified
using standard PDE theory.
For instance, combining \cite[p.~412, Proposition 8.3]{Taylor2023PDEIII} with \cite[p.~6, Proposition 2.4]{Taylor2023PDEIII} implies that $v$ and all derivatives indexed by $\mathcal{M}_{2,4}$ are bounded on $[0, T] \times \R^d$, provided that the PDE~\eqref{eq_pde} satisfies the strong parabolicity condition \cite[p.~409, (8.2)]{Taylor2023PDEIII} and that the terminal function $x \mapsto v(T, x)$ belongs to $\mathbb{H}^{k,2}(\R^d)$ for some $k > d/2 + 6$.
Here $\mathbb{H}^{k, 2}(\R^d)$ denotes the Sobolev space of functions on $\R^d$ with square-integrable weak derivatives up to order $k$.
The moment condition \eqref{eq_finite_moment} is mild and holds for common choices of $\xi$, such as a standard Gaussian random vector.

The following \Cref{thm_pre_consist} is a key result and will be used in the analyses of local truncation error and stability of the RDM formulation~\eqref{eq_ER0}.

\begin{theorem}\label{thm_pre_consist}
Let \cref{assum_polygrow} hold. For any $(t, x) \in [0, T) \times \R^d$ and any step size $h$ with $0 < h < \min\{1,\, T - t\}$, it holds that
\begin{equation}
    \big|\hat{\mathcal{D}}_h v(t, x) - \hat{\mathcal{D}} v(t, x)\big| \leq C_{\mathrm{loc}} h \big(1 + |x|^{\bar{p}}\big),
\end{equation}
where
\begin{align}
    \hat{\mathcal{D}} v(t, x) &:= \partial_t v(t, x) + \hat{\mu} \partial_x v(t, x) + \frac{1}{2} \mathrm{Tr}[\hat{\sigma} \hat{\sigma}^{\top} \partial_{xx} v(t, x)], \label{eq_hatD}\\
    \hat{\mathcal{D}}_h v(t, x) &:= \frac{1}{h} \E{v(t + h, x + \hat{\mu} h + \hat{\sigma} \sqrt{h} \xi) - v(t, x)}, \label{eq_hatDh}\\
    \bar{p} &:= p_v + 3 p_{\mu} \max\{1,\, p_v\} + 3 p_{\sigma} \max\{2,\, p_v\},\label{eq_defbarp}
\end{align}
with $\hat{\mu} := \mu(t, x, \hat{v}(t, x))$ and $\hat{\sigma} := \sigma(t, x, \hat{v}(t, x))$.
The constant $C_{\mathrm{loc}} > 0$ is independent of $\hat{v}$, $t$, $h$, $T$, and $x$; see \eqref{eq_defCD} for its explicit form.
\end{theorem} 

\begin{proof}
        
    For $s \in [0, T-t]$, define $\hat{V}(s, z)$ as the analogue of $V(s, z)$ in \eqref{eq_defVsz}, obtained by replacing $(\mu, \sigma)$ with $(\hat{\mu}, \hat{\sigma})$ given in the statement of \Cref{thm_pre_consist}; namely,
    \begin{equation}\label{eq_hatVsz}
        \hat{V}(s, z) := v(t + s, x + \hat{\mu}\, s + \hat{\sigma}\, z).
    \end{equation}
    Then, by direct analogy with the derivations in \Cref{eq_expandV,eq1_expandV,eq2_expandV}, we obtain the following expansions:
    \begin{align}
        &\E{\hat{V}(h, \sqrt{h}\xi)} = \hat{V}(h, 0) + \frac{h}{2} \sum_{i=1}^q \frac{\partial^2 \hat{V}}{\partial z_i^2} (h, 0) + R_1, \label{eq_exp1}\\
        &\hat{V}(h, 0) = \hat{V}(0, 0) + \br{\partial_t + \hat{\mu}^{\top} \partial_x} v(t, x) h + R_2, \label{eq_exp2} \\
        &\sum_{i=1}^q \frac{\partial^2 \hat{V}}{\partial z_i^2} (h, 0) = \mathrm{Tr}\rbr{\hat{\sigma} \hat{\sigma}^{\top} \partial_{xx} v(t, x)} + R_3, \label{eq_exp3}
    \end{align}
    where the remainder terms $R_1$, $R_2$, and $R_3$ are defined as
    \begin{equation}\label{eq_defR1R2R3}
    \begin{aligned}
        &R_1 := \frac{h^2}{4!} \E{\vbr{\Big}{\xi^{\top} \frac{\partial}{\partial z}}^4 \hat{V}(h, c_1 \sqrt{h}\xi)}, \\ 
        &R_2 := \frac{h^2}{2} \frac{\partial^2 \hat{V}}{\partial s^2}(c_2 h, 0), \quad R_3 := h \sum_{i=1}^q \frac{\partial^3 \hat{V}}{\partial z_i \partial z_i \partial s}(c_3 h, 0).
    \end{aligned}
    \end{equation}
    Here, $c_1$ is a random variable depending on $\sqrt{h}\xi$, while $c_2$ and $c_3$ are deterministic constants depending on $h$; all take values in $[0, 1]$. 
    Further, by analogy with \Cref{eq_rdo}, the expansions \eqref{eq_exp1}--\eqref{eq_exp3} lead to
    \begin{equation}\label{eq_expandV3}
        \hat{\mathcal{D}}_h v(t, x) = \hat{\mathcal{D}}v(t, x) + \frac{R_1 + R_2}{h} + \frac{1}{2} R_3
    \end{equation}
    with $\hat{\mathcal{D}}$ and $\hat{\mathcal{D}}_h$ defined in \eqref{eq_hatD} and \eqref{eq_hatDh}, respectively.
    
    We now derive bounds for $R_1$, $R_2$, and $R_3$, and combine them to obtain an explicit bound on the remainder in \eqref{eq_expandV3}. The proof proceeds in four steps.    
    
    \subsection*{(i) Notation and derivative bounds.}
    
    To facilitate the estimation of the remainder terms $R_1$, $R_2$, and $R_3$, we first introduce some commonly used notation:
    \begin{align}
        S_{\mu} &:= 1 + \max_{i=1,\cdots, d} \abs{\hat{\mu}_i} + \sum_{i=1}^d \abs{\hat{\mu}_i}^{p_v}, \quad M_{p_v} := \max\{1, 3^{p_v-1}\}, \label{eq_defSmu} \\
        S_{\sigma} &:= 1 + \sum_{j = 1}^q \max_{i=1,\cdots, d} \abs{\hat{\sigma}_{ij}}^2 + \sum_{i=1}^d \max_{j=1,\cdots,q} \abs{\hat{\sigma}_{ij}}^{p_v}, \label{eq_defSsgm}\\
        C_{\xi} &:= \Big\{1 + \sqrt{\mathbb{E}\Big[\abs{\xi}^{2 p_v}\Big]}\Big\} \times \Big\{1 + \vbr{\Big}{\sum_{i,j,k,l=1}^q \E{\abs{\xi_i \xi_j \xi_k \xi_l}^2}}^{1/2} \Big\}.  \label{eq_defCxi} 
    \end{align}
    Here, $S_{\mu}$ and $S_{\sigma}$ are functions of $(t, x)$, where $\hat{\mu}$ and $\hat{\sigma}$ are evaluated at $(t, x, \hat{v}(t, x))$; $C_{\xi}$ is a constant independent of $t$, $h$, and $x$, and is finite by the moment condition in \eqref{eq_finite_moment};
    the definition of $M_{p_v}$ ensures that $(a + b + c)^{p_v} \leq M_{p_v} (a^{p_v} + b^{p_v} + c^{p_v})$ for any $a, b, c \geq 0$, which will be frequently used in the following estimates.
    
    Next, we derive bounds for the derivatives of $v$, which in turn control those of $\hat{V}$ in \eqref{eq_hatVsz}.
    For any $(t, x) \in [0, T) \times \R^d$, $s \in [0, T - t) \cap [0, 1]$, and $z \in \R^q$, we have
    \begin{align}
        &\sum_{\alpha \in \mathcal{M}_{2, 4}} \abs{(D^{\alpha} v)(t + s, x + \hat{\mu} s + \hat{\sigma} z)} \leq C_{\mathrm{g}} \vbr{\big}{1 + \big|x + \hat{\mu} s + \hat{\sigma} z \big|^{p_v}} \notag \\
        \leq\;& C_{\mathrm{g}} M_{p_v} \br{1 + \abs{x}^{p_v} + \abs{\hat{\mu}}^{p_v} + \abs{\hat{\sigma} z}^{p_v}} \notag \\
        \leq\;& C_{\mathrm{g}} M_{p_v} \br{1 + \abs{x}^{p_v} + S_{\mu} + S_{\sigma} \abs{z}^{p_v}},\label{eq_estDv}
    \end{align}
    where the first and second inequalities follow from \eqref{eq_polygrow} and the Cauchy-Schwarz inequality, respectively. 
    
    \subsection*{(ii) Estimate of $R_1$.}
    
    To estimate $R_1$ in \eqref{eq_defR1R2R3}, we deduce
    \begin{equation}\label{eq_estR1}
        \abs{R_1} \leq \frac{h^2}{4!} \sum_{i,j,k,l=1}^q \mathbb{E}\Big[\Big|\xi_i \xi_j \xi_k \xi_l\Big| \Big|\frac{\partial^4 \hat{V}}{\partial {z_i} \partial {z_j} \partial {z_k} \partial {z_l}}(h, c_1 \sqrt{h}\xi)\Big|\Big] \leq \frac{h^2}{4!} C_{\xi} \mathbb{I}_1,
    \end{equation}
    where the second inequality follows from the Cauchy-Schwarz inequality with $C_{\xi}$ given by \eqref{eq_defCxi}, and
    $\mathbb{I}_1$ given by 
    \begin{align}
        \mathbb{I}_1 :=&\vbr{\bigg}{\sum_{i,j,k,l=1}^q \mathbb{E}\bigg[\Big|\frac{\partial^4 \hat{V}}{\partial {z_i} \partial {z_j} \partial {z_k} \partial {z_l}}(h, c_1 \sqrt{h}\xi)\Big|^2\bigg]}^{1/2} \notag\\
        \leq\;& \Big(d \sum_{i,j,k,l=1}^d \,\sum_{i^{\prime},j^{\prime},k^{\prime},l^{\prime}=1}^q \Big|\hat{\sigma}_{ii^{\prime}} \hat{\sigma}_{jj^{\prime}} \hat{\sigma}_{kk^{\prime}} \hat{\sigma}_{ll^{\prime}}\Big|^2 \times\notag\\
        &\qquad\qquad\qquad \mathbb{E}\Big[\big|\frac{\partial^4 v}{\partial_{x_{i}} \partial_{x_{j}} \partial_{x_{k}} \partial_{x_{l}}}(t + h, x + \hat{\mu} h + \hat{\sigma} c_1 \sqrt{h}\xi) \big|^2\Big]\Big)^{1/2} \notag\\
        \leq\;& \sqrt{d} S_{\sigma}^2 \times \mathbb{I}_2 \label{eq_estI1}
    \end{align}
    where the second line follows from the definition \eqref{eq_hatVsz} of $\hat{V}(s, z)$ and the Cauchy-Schwarz inequality; the last inequality follows from the Hölder inequality with $S_{\sigma}$ given in \eqref{eq_defSsgm}, and $\mathbb{I}_2$ given by
    \begin{align}
        \mathbb{I}_2 :=\;&\vbr{\bigg}{\sum_{i,j,k,l=1}^d \E{\abs{\frac{\partial^4 v}{\partial_{x_{i}} \partial_{x_{j}} \partial_{x_{k}} \partial_{x_{l}}}(t + h, x + \hat{\mu} h + \hat{\sigma} c_1 \sqrt{h}\xi)}^2 }}^{1/2} \notag\\
        \leq\;& \vbr{\Big}{\mathbb{E}\Big[\vbr{\Big}{\sum_{\alpha \in \mathcal{M}_{2, 4}} |D^{\alpha} v|(t + h, x + \hat{\mu} h + \hat{\sigma} c_1 \sqrt{h} \xi)}^2\Big]}^{1/2} \notag\\
        \leq\;& C_{\mathrm{g}} M_{p_v} \vbr{\big}{\E{\vbr{\big}{1 + \abs{x}^{p_v} + S_{\mu} + S_{\sigma} \abs{\xi}^{p_v}}^2}}^{1/2} \notag\\
        \leq\;& 4 C_{\mathrm{g}} M_{p_v} \vbr{\big}{1 + \abs{x}^{p_v} + S_{\mu} + S_{\sigma} C_{\xi}}.\label{eq_estI2}
    \end{align}
    Here the third line follows from \eqref{eq_estDv} and $|c_1| \sqrt{h} < 1$; the last line follows from the Minkowski inequality and the definition of $C_{\xi}$ in \eqref{eq_defCxi}. 
    Inserting \eqref{eq_estI2} into \eqref{eq_estI1}, and further inserting the resulting equation into \eqref{eq_estR1}, we obtain
    \begin{equation}\label{eq2_estR1}
        \abs{R_1} \leq \frac{h^2}{4!} \sqrt{d} C_{\xi} C_{\mathrm{g}} M_{p_v} S_{\sigma}^2 \vbr{\big}{1 + \abs{x}^{p_v} + S_{\mu} + S_{\sigma} C_{\xi}}.
    \end{equation}
    
    \subsection*{(iii) Estimates of $R_2$ and $R_3$.}
    We introduce the convention
    \begin{equation*}
        \hat{\mu}_0 := 1, \quad \frac{\partial v}{\partial x_0} := \frac{\partial v}{\partial t},
    \end{equation*}
    so that the derivative of $\hat{V}(s,0)$ with respect to $s$ in \eqref{eq_hatVsz} can be written as
    \begin{equation}\label{eq_PVPs}
        \frac{\partial \hat{V}}{\partial s}(s, 0)
        = \sum_{i=0}^d \hat{\mu}_i \frac{\partial v}{\partial x_i}\big(t + s,\, x + \hat{\mu} s\big).
    \end{equation}
    Substituting the above expression into the definition of $R_2$ in \eqref{eq_defR1R2R3}, we obtain
    \begin{align}
        \abs{R_2} \leq\;& \frac{h^2}{2}\abs{\frac{\partial^2 \hat{V}}{\partial s^2}(c_2 h, 0)} = \frac{h^2}{2} \abs{\sum_{i,j=0}^d \hat{\mu}_i \hat{\mu}_j \frac{\partial^2 v}{\partial x_i \partial x_j} (t + c_2 h, x + \hat{\mu} c_2 h)} \notag\\
        \leq\;& \frac{h^2}{2} S_{\mu}^2 C_{\mathrm{g}} M_{p_v} \vbr{\big}{1 + \abs{x}^{p_v} + S_{\mu}}, \label{eq_estR2}
    \end{align}
    where $\hat{\mu}$ is evaluated at $(t, x, \hat{v}(t, x))$; 
    the inequality follows from the definition of $S_{\mu}$ in \eqref{eq_defSmu}, and the bound \eqref{eq_estDv} with $s = c_2 h$ and $z = 0$.
    
    To estimate $R_3$ in \eqref{eq_defR1R2R3}, we first use the definition \eqref{eq_hatVsz} of $\hat{V}(s, z)$ to deduce
    \begin{align}
        \abs{R_3} \leq\;& h \Bigg| \sum_{i,j=1}^d \sum_{k=1}^q \sum_{l=0}^d \hat{\sigma}_{ik} \hat{\sigma}_{jk} \hat{\mu}_l \frac{\partial^3 v}{\partial_{x_i} \partial_{x_j} \partial_{x_l}} (t + c_3 h, x + \hat{\mu} c_3 h) \Bigg| \notag\\
        \leq\;& h \max_{i,j,l=1,\cdots,d} \sum_{k=1}^q \abs{\hat{\sigma}_{ik} \hat{\sigma}_{jk}} \abs{\hat{\mu}_l} \times \sum_{i,j,l=0}^d \Big|\frac{\partial^3 v}{\partial_{x_i} \partial_{x_j} \partial_{x_l}} (t + c_3 h, x + \hat{\mu} c_3 h) \Big| \notag \\
        \leq\;& h S_{\sigma} S_{\mu} \times C_{\mathrm{g}} M_{p_v} \vbr{\big}{1 + \abs{x}^{p_v} + S_{\mu}}, \label{eq_estR3}
    \end{align}
    where the first line uses \eqref{eq_PVPs} to express $\partial \hat{V}/\partial s$;
    the second line uses the Hölder inequality; 
    the third line follows from the Cauchy-Schwarz inequality applied to $\sum_{k=1}^q \abs{\hat{\sigma}_{ik}\hat{\sigma}_{jk}}$, together with the bound \eqref{eq_estDv} (with $s = c_3 h$ and $z = 0$), to bound the sum of derivatives of $v$.

    \subsection*{(iv) Estimate of the remainder term.}
    
    Combining \eqref{eq2_estR1}, \eqref{eq_estR2}, and \eqref{eq_estR3},
    we obtain 
    \begin{equation}\label{eq_R1R2R3}
    \begin{aligned}
        \abs{\frac{R_1 + R_2}{h} + \frac{1}{2} R_3} \leq\;& 2 h \sqrt{d} C_{\xi} C_{\mathrm{g}} M_{p_v} S_{\mu}^2 S_{\sigma}^2 \vbr{\big}{1 + \abs{x}^{p_v} + S_{\mu} + S_{\sigma} C_{\xi}}\\
        \leq\;& 6 h \sqrt{d} C_{\mathrm{g}} M_{p_v} C_{\xi}^2 S_{\mu}^3 S_{\sigma}^3 \vbr{\big}{1 + \abs{x}^{p_v}},
    \end{aligned}
    \end{equation}
    where we have used the inequalities $M_{p_v}, C_{\xi}, S_{\mu}, S_{\sigma} \geq 1$ to simplify the expression.
    It remains to establish bounds on $S_{\mu}$ and $S_{\sigma}$.
    Inserting the polynomial growth conditions \eqref{eq_lgrow_mu} and \eqref{eq_lgrow_sgm} into the right-hand sides of \eqref{eq_defSmu} and \eqref{eq_defSsgm}, we obtain
    \begin{align}
        S_{\mu} &\leq 1 + \abs{\hat{\mu}} + d \abs{\hat{\mu}}^{p_v} \leq C_{\hat{\mu}\hat{\sigma}} \br{1 + \abs{x}^{p_{\mu} \max\{1,\, p_v\}}}, \label{eq_bdSmu} \\
        S_{\sigma} &\leq 1 + q \abs{\hat{\sigma}}^2 + d \abs{\hat{\sigma}}^{p_v} \leq C_{\hat{\mu}\hat{\sigma}} \br{1 + \abs{x}^{p_{\sigma} \max\{2,\, p_v\}}}, \label{eq_bdSsgm}
    \end{align}
    where $C_{\hat{\mu}\hat{\sigma}}$ is a constant defined by
    \begin{equation*}
        C_{\hat{\mu}\hat{\sigma}} := \br{1 + (q + d) C_{\mathrm{g}}^{\max\{2,\, p_v\}}} M_{p_v}.
    \end{equation*}
    Inserting \eqref{eq_bdSmu} and \eqref{eq_bdSsgm} into \eqref{eq_R1R2R3}, we obtain 
    \begin{equation}\label{eq_finalR1R2R3}
        \abs{\frac{R_1 + R_2}{h} + \frac{1}{2} R_3} \leq C_{\mathrm{loc}} h \vbr{\big}{1 + \abs{x}^{\bar{p}}},
    \end{equation}
    where $\bar{p}$ is defined in \eqref{eq_defbarp} and 
    \begin{equation}\label{eq_defCD}
        C_{\mathrm{loc}} := 6 \times 2^7 \times \sqrt{d} C_{\mathrm{g}} C_{\xi}^2 M_{p_v}^7 \vbr{\big}{1 + (q + d) C_{\mathrm{g}}^{\max\{2,\, p_v\}}}^6
    \end{equation}
    with $M_{p_v}$ and $C_{\xi}$ given by the second equation in \eqref{eq_defSmu} and \eqref{eq_defCxi}, respectively.
    Inserting \eqref{eq_finalR1R2R3} into \eqref{eq_expandV3}, we obtain the desired result. 
\end{proof}

The following corollary, obtained by setting $\hat{v}=v$ in \Cref{thm_pre_consist}, establishes that the local truncation error of the random difference approximation \eqref{eq_rdo} is first order in the step size $h$.
\begin{corollary}\label{coro_consist}
    Let \cref{assum_polygrow} hold. 
    For any $(t, x) \in [0, T) \times \R^d$ and any $0 < h < \min\{1,\, T - t\}$, it holds that
    \begin{equation}
        \big|\mathcal{D}_h v(t, x) - \mathcal{D} v(t, x)\big| \leq C_{\mathrm{loc}}\, h \big(1 + |x|^{\bar{p}}\big),
    \end{equation}
    where $\mathcal{D}$ and $\mathcal{D}_h$ are defined in \eqref{eq_defD} and \eqref{eq_rdo}, respectively; the constants $C_{\mathrm{loc}}$ and $\bar{p}$ are the same as those in \Cref{thm_pre_consist}.
\end{corollary}

\subsection{Stability analysis}

Our stability analysis relies on the following two assumptions.
\begin{assumption}\label{assum_Lip}
There exist constants $C_{\mu}$, $C_{\sigma}$, $C_f$, and $C_{\mathrm{nl}} \geq 0$ such that the functions $\mu$, $\sigma$, $f$, and $v$ appearing in \eqref{eq_pde} satisfy 
\begin{align}
    &\abs{\mu(t, x, y_1) - \mu(t, x, y_2)} \leq C_{\mu} \abs{y_1 - y_2}, \label{eq_lipmu} \\
    &\abs{\sigma(t, x, y_1) - \sigma(t, x, y_2)} \leq C_{\sigma} \abs{y_1 - y_2}, \label{eq_lipsgm}\\
    &\abs{f(t, x, y_1) - f(t, x, y_2)} \leq C_f \abs{y_1 - y_2}, \label{eq_lipf}\\
    &C_{\mu} \abs{\partial_x v(t, x)} + C_{\sigma} \abs{\sigma(t, x, y)} \abs{\partial_{xx} v(t, x)} \leq C_{\mathrm{nl}}, \label{eq_defCnl}
\end{align}
for all $t \in [0, T]$, $x \in \R^d$, and $y_1, y_2, y \in \R$. 
\end{assumption}

\begin{assumption}\label{assum_em}
The sample paths $\{X_n^m\}_{n=0}^N$ generated by the Euler-Maruyama scheme~\eqref{eq_pathX} satisfy
\begin{equation*}
    \max_{0 \leq n \leq N} \E{\abs{X_n^m}^{2 \bar{p}}} \leq C_{\mathrm{EM}}
\end{equation*}
for $m = 1, 2, \cdots, M$, where $C_{\mathrm{EM}} \geq 1$ is a constant independent of $m$, $M$, $N$, $h$ and $T$, and $\bar{p}$ is defined in~\eqref{eq_defbarp}.
\end{assumption}

The Lipschitz conditions \eqref{eq_lipmu}--\eqref{eq_lipf} in \Cref{assum_Lip} are standard in the numerical analysis of PDEs.
Condition \eqref{eq_defCnl} reflects the nonlinearity of the operator $\mathcal{D}$ in \eqref{eq_defD}.
If $\mathcal{D}$ is linear, the constants $C_{\mu}$ and $C_{\sigma}$ in \eqref{eq_lipmu} and \eqref{eq_lipsgm} can be set to zero, and then \Cref{assum_Lip} holds trivially with $C_{\mathrm{nl}}=0$.
Alternatively, \eqref{eq_defCnl} is satisfied if $\abs{\partial_x v}$ and $\abs{\sigma}\,\abs{\partial_{xx} v}$ are bounded, a condition commonly met in practice.

\Cref{assum_em} is introduced to avoid technical complications arising in the stability analysis of the Euler-Maruyama scheme~\eqref{eq_pathX}, which lies beyond the main scope of this work.
A simple sufficient condition for \Cref{assum_em} is that the functions $\sum_{\alpha \in \mathcal{M}_{2, 4}} \abs{D^{\alpha} v}$, $\mu$, and $\sigma$ are bounded, which implies $p_{\mu} = p_{\sigma} = p_v = 0$, and $\bar{p} = 0$ recalling \eqref{eq_defbarp}.
In this case, \Cref{assum_em} holds trivially.
More general sufficient conditions are well studied in the field of numerical SDEs. 
For example, \cite[Lemma 2.1]{Yuan2008note} provides a result for arbitrary $\bar{p}$ when $\mu(t, x, \hat{v}(t, x))$ and $\sigma(t, x, \hat{v}(t, x))$ grow at most linearly in $x$.

The following lemma establishes the stability of the formulation \eqref{eq_ER0}.

\begin{lemma}\label{lem_priest}
    Under \Cref{assum_polygrow,assum_Lip,assum_em}, it holds that 
    \begin{equation}\label{eq_estMn2Delv}
        \max_{0 \leq n\leq N-1} \mathbb{M}_n^2\rbr{\Delta v} \leq C_{\mathrm{st1}} \exp\br{C_{\mathrm{st2}} T} \Big\{\mathbb{M}_{N}^2\rbr{\Delta v} + \frac{T}{N}\sum_{n=0}^{N-1} \mathbb{M}_n^2\rbr{\Delta R} + T h^2 \Big\}
    \end{equation}
    for $0 < h \leq \min\{1, (24 C_{\mathrm{nl}}^2 + 12 C_f^2)^{-1}\}$, where  
    \begin{align}
        \Delta v(t, x) &:= v(t, x) - \hat{v}(t, x), \label{eq_defdeltav} \\
        \Delta R(t, x) &:= \E{R(t, x, \xi; v)} - \E{R(t, x, \xi; \hat{v})} \label{eq_defRi}
    \end{align}
    with $R$ given in \eqref{eq_defR}, and the functionals $\mathbb{M}_n^2\rbr{\cdot}$ are defined in \eqref{eq_defMtp}. 
    In \eqref{eq_estMn2Delv}, $C_{\mathrm{st1}}$ and $C_{\mathrm{st2}}$ are stability constants independent of $h$, $N$ and $\hat{v}$, given in \eqref{eq_Cst1} and \eqref{eq_estCh}, respectively.
\end{lemma}

\begin{proof}
    
    
    Let $(t, x) \in [0, T - h) \times \R^d$ be arbitrary.
    Substituting the definition of $R$ in \eqref{eq_defR} into \eqref{eq_defRi}, we obtain
    \begin{equation}\label{eq_DeltaR}
        \Delta R(t, x) = \frac{\Delta V(t, x) - \Delta v(t, x)}{h} - \Delta f(t, x),
    \end{equation}
    where $\Delta f(t, x) := f(t, x, v(t, x)) - f(t, x, \hat{v}(t, x))$; $\Delta v$ is defined in \eqref{eq_defdeltav}, and
    \begin{equation}\label{eq_defVh}
        \Delta V(t, x) := \E{v(t + h, x + \xi_h)} - \mathbb{E}[\hat{v}(t + h, x + \hat{\xi}_h)]. 
    \end{equation}
    Here $\xi_h$ is given by \eqref{eq_def_xih}; $\hat{\xi}_h$ is defined analogously, i.e., 
    \begin{equation}\label{eq_defhatxi}
        \hat{\xi}_h := \hat{\mu} h + \hat{\sigma} \sqrt{h} \,\xi
    \end{equation}
    with $\hat{\mu} = \mu(t, x, \hat{v}(t, x))$ and $\hat{\sigma} = \sigma(t, x, \hat{v}(t, x))$.
    The equation \eqref{eq_defVh} can be further decomposed as 
    \begin{equation}\label{eq_decompDeltaV}
        \Delta V(t, x) = \Delta V_1(t, x) + \Delta V_2(t, x),
    \end{equation}
    where 
    \begin{align}
        \Delta V_1(t, x) :=\;& \E{v(t + h, x + \xi_h) - v(t + h, x + \hat{\xi}_h)}, \label{eq_defDeltaV1} \\ 
        \Delta V_2(t, x) :=\;& \E{v(t + h, x + \hat{\xi}_h) - \hat{v}(t + h, x + \hat{\xi}_h)} \notag\\
        =\;& \E{\Delta v(t + h, x + \hat{\xi}_h)}. \label{eq_defDeltaV2}
    \end{align}
    
    Inserting \eqref{eq_decompDeltaV} into \eqref{eq_DeltaR}, and rearranging the terms, we obtain
    \begin{equation}\label{eq_EDeltaV}
        \Delta v(t, x) = \Delta V_1(t, x) + \Delta V_2(t, x) - h \Delta f(t, x) - h \Delta R(t, x).
    \end{equation}
    Squaring both sides of \eqref{eq_EDeltaV}, and using the inequality $(a + b)^2 \leq (1+h) a^2 + \br{1+h^{-1}} b^2$ for $a, b > 0$, we obtain
    \begin{equation}\label{eq_bdDelv}
    \begin{aligned}
        \abs{\Delta v(t, x)}^2 \leq\;& (1 + h) \abs{\Delta V_2(t, x)}^2 + 3 \br{1 + h^{-1}} \abs{\Delta V_1(t, x)}^2 \\
        &+ 3\br{1 + h^{-1}} h^2 \bbr{\abs{\Delta R(t, x)}^2 + \abs{\Delta f(t, x)}^2}.
    \end{aligned}
    \end{equation}
    Starting from \eqref{eq_bdDelv}, the remainder of the proof is devoted to establishing a bound for $|\Delta v|$ at $t = t_n$ in terms of $|\Delta v|$ at $t = t_{n+1}$, and then applying a discrete Gronwall inequality to obtain the final result.
    We proceed in four steps.
    
    \subsubsection*{(i) Estimate of $\Delta V_1(t, x)$.}
    
    The main results of \Cref{thm_pre_consist} and \Cref{coro_consist} imply that
    \begin{align*}
        \big| h^{-1} \big\{\mathbb{E}[v(t + h, x + \hat{\xi}_h)] - v(t, x)\big\} - \hat{\mathcal{D}} v(t, x) \big| \leq C_{\mathrm{loc}} h \big(1 + |x|^{\bar{p}}\big),\\
        \big| h^{-1} \bbr{\E{v(t + h, x + \xi_h)} - v(t, x)} - \mathcal{D} v(t, x) \big| \leq C_{\mathrm{loc}} h \big(1 + |x|^{\bar{p}}\big),
    \end{align*}
    where $\mathcal{D}$ and $\hat{\mathcal{D}}$ are defined in \eqref{eq_defD} and \eqref{eq_hatD}, respectively; the constant $C_{\mathrm{loc}}$ and exponent $\bar{p}$ are given in \Cref{thm_pre_consist}.
    Combining the above two inequalities, the term $\Delta V_1(t, x)$ in \eqref{eq_defDeltaV1} can be bounded by
    \begin{equation}\label{eq_bdDeltaV1}
        \abs{\Delta V_1(t, x)} \leq h \big|(\mathcal{D} - \hat{\mathcal{D}}) v(t, x) \big| + 2 C_{\mathrm{loc}} h^2 \big(1 + |x|^{\bar{p}}\big). 
    \end{equation}
    
    Using the definitions in \eqref{eq_defD} and \eqref{eq_hatD}, the term $|(\mathcal{D} - \hat{\mathcal{D}}) v(t, x)|$ can be bounded by
    \begin{align}
        \big|(\mathcal{D} - \hat{\mathcal{D}}) v(t, x) \big| 
        \leq& \abs{\mu - \hat{\mu}} \abs{\partial_x v} + \frac{1}{2} \br{\abs{\sigma} + \abs{\hat{\sigma}}} \abs{\sigma - \hat{\sigma}} \abs{\partial_{xx} v} \notag\\
        \leq\;& C_{\mu} \abs{\Delta v(t, x)} \abs{\partial_x v} + \frac{1}{2} (\abs{\sigma} + \abs{\hat{\sigma}}) C_{\sigma} \abs{\Delta v(t, x)} \abs{\partial_{xx} v} \notag\\
        \leq\;& C_{\mathrm{nl}} \abs{\Delta v(t, x)}, \label{eq_DelDh}
    \end{align}
    where the second line follows from the Lipschitz conditions \eqref{eq_lipmu} and \eqref{eq_lipsgm}, and the last line uses the definition of $C_{\mathrm{nl}}$ in \eqref{eq_defCnl}.
    Inserting \eqref{eq_DelDh} into \eqref{eq_bdDeltaV1} we obtain 
    \begin{equation}\label{eq2_bdDeltaV1}
        \abs{\Delta V_1(t, x)} \leq C_{\mathrm{nl}} h \abs{\Delta v(t, x)} + 2 C_{\mathrm{loc}} h^2 \big(1 + |x|^{\bar{p}}\big). 
    \end{equation}
    
    \subsubsection*{(ii) Estimate of $\Delta v(t, x)$. }
    
    Using \eqref{eq2_bdDeltaV1}, the second term on the right side of \eqref{eq_bdDelv} can be bounded by 
    \begin{equation}\label{eq3_bdDeltaV1}
    \begin{aligned}
        &3 \br{1 + h^{-1}} \abs{\Delta V_1(t, x)}^2 \\
        \leq\;& 3 \br{1 + h^{-1}} \bbr{2 C_{\mathrm{nl}}^2 h^2 \abs{\Delta v(t, x)}^2 + 16 C_{\mathrm{loc}}^2 h^4 \big(1 + |x|^{2 \bar{p}}\big)} \\
        \leq\;& 12 C_{\mathrm{nl}}^2 h \abs{\Delta v(t, x)}^2 + 96 C_{\mathrm{loc}}^2 h^3 \big(1 + |x|^{2 \bar{p}}\big),
    \end{aligned}
    \end{equation}
    where the second line follows from the inequality $(a + b)^2 \leq 2 a^2 + 2 b^2$ for $a, b > 0$, and the last line uses $0 < h \leq 1$.
    Using the Lipschitz condition~\eqref{eq_lipf}, the last term on the right side of \eqref{eq_bdDelv} can be bounded by
    \begin{equation}\label{eq_estDelf}
        3 \br{1 + h^{-1}} h^2 \abs{\Delta f(t, x)}^2 \leq 6 h C_f^2 \abs{\Delta v(t, x)}^2,
    \end{equation}
    where we have used $0 < h \leq 1$ to simplify the expression.
    Inserting \eqref{eq3_bdDeltaV1} and \eqref{eq_estDelf} into \eqref{eq_bdDelv}, we obtain 
    \begin{equation}\label{eq2_estDelv2}
    \begin{aligned}
        \abs{\Delta v(t, x)}^2 \leq\;& (1 + h) \br{\abs{\Delta V_2(t, x)}^2 + 3 h \abs{\Delta R(t, x)}^2} \\
        &+ C_0 h \abs{\Delta v(t, x)}^2 + 96 C_{\mathrm{loc}}^2 h^3 \big(1 + |x|^{2 \bar{p}}\big) \\
        \leq\;& C_0 h \abs{\Delta v(t, x)}^2 + (1 + h) \vbr{\big}{\abs{\Delta V_2(t, x)}^2 + \mathbb{I}_3},
    \end{aligned}
    \end{equation}
    where $C_0$ and $\mathbb{I}_3$ are defined by
    \begin{equation}\label{eq_defI3}
        C_0 := 12 C_{\mathrm{nl}}^2 + 6 C_f^2, \quad \mathbb{I}_3 := 3 h \abs{\Delta R(t, x)}^2 + 96 C_{\mathrm{loc}}^2 h^3 \big(1 + |x|^{2 \bar{p}}\big). 
    \end{equation}
    
    Rearranging \eqref{eq2_estDelv2} gives
    \begin{equation}\label{eq_ineq_1h}
        \br{1 - C_0 h}\abs{\Delta v(t, x)}^2 \leq (1 + h) \vbr{\big}{\abs{\Delta V_2(t, x)}^2 + \mathbb{I}_3}.   
    \end{equation}
    Combining the condition on $h$ in the statement of \Cref{lem_priest} with the definition of $C_0$ in \eqref{eq_defI3}, we have
    \begin{equation}\label{eq_estCh}
        1 - C_0 h \geq \frac{1}{2} > 0, \quad \frac{1 + h}{1 - C_0 h} \leq 1 + C_{\mathrm{st2}} h, \quad C_{\mathrm{st2}} := 2 + 24 C_{\mathrm{nl}}^2 + 12 C_f^2. 
    \end{equation}
    Then \eqref{eq_ineq_1h} leads to 
    \begin{equation}\label{eq_estDelv2}
    \begin{aligned}
        \abs{\Delta v(t, x)}^2 \leq\;& \frac{1 + h}{1 - 2\br{1 + h} h C_f^2} \vbr{\big}{\abs{\Delta V_2(t, x)}^2 + \mathbb{I}_3} \\
        \leq\;& \vbr{\big}{1 + C_{\mathrm{st2}} h} \vbr{\big}{\abs{\Delta V_2(t, x)}^2 + \mathbb{I}_3},
    \end{aligned}
    \end{equation}
    where the second inequality follows from the second inequality in \eqref{eq_estCh}.

    \subsubsection*{(iii) Bounding $\mathbb{M}_n^2\rbr{\Delta v}$ by $\mathbb{M}_{n+1}^2\rbr{\Delta v}$. }
    
    Using \eqref{eq_estDelv2}, set $t = t_n$, multiply both sides by the PDF $P_{X_n}(x)$ introduced in \eqref{eq_defMtp}, and integrate over $x \in \R^d$ to obtain
    \begin{equation}\label{eq_estMn2}
        \mathbb{M}_n^2\rbr{\Delta v} \leq (1 + C_{\mathrm{st2}} h) \br{\mathbb{M}_n^2\rbr{\Delta V_2} + \mathbb{I}_{4, n}},
    \end{equation}
    where $\mathbb{M}_n^2\rbr{\,\cdot\,}$ is defined in \eqref{eq_defMtp}, and $\mathbb{I}_{4, n}$ denotes the integral of $\mathbb{I}_3$ in \eqref{eq_defI3} evaluated at $t = t_n$ with respect to $P_{X_n}(x)$, i.e.,
    \begin{equation}\label{eq_defI4n}
    \begin{aligned}
        \mathbb{I}_{4, n} :=\;& 3 h \mathbb{M}_n^2\rbr{\Delta R} + 96 C_{\mathrm{loc}}^2 h^3 \int_{\R^d} \br{1 + |x|^{2 \bar{p}}} P_{X_n}(x) \di x \\
        \leq\;& 3 h \mathbb{M}_n^2\rbr{\Delta R} + 96 C_{\mathrm{loc}}^2 (1 + C_{\mathrm{EM}}) h^3
    \end{aligned}
    \end{equation}
    with $C_{\mathrm{EM}}$ given in \Cref{assum_em}.
    

    Now we estimate $\mathbb{M}_n^2\rbr{\Delta V_2}$ in \eqref{eq_estMn2}.
    Using the definition of $\Delta V_2$ in \eqref{eq_defDeltaV2} and the Jensen inequality, we have
    \begin{equation}\label{eq_estMnV2}
    \begin{aligned}
        \mathbb{M}_n^2\rbr{\Delta V_2} &\leq \int_{\R^d} \E{\abs{\Delta v(t_{n+1}, x + \hat{\xi}_h)}^2} P_{X_n}\br{x} \, \di x \\
        &= \int_{\R^d} \int_{\R^d} \abs{\Delta v(t_{n+1}, z)}^2 P_{x + \hat{\xi}_h}\br{z} P_{X_n}\br{x} \, \di z \, \di x,
    \end{aligned}
    \end{equation}
    where $P_{x + \hat{\xi}_h}(\cdot)$ denotes the PDF of $x + \hat{\xi}_h$ with $\hat{\xi}_h$ defined in \eqref{eq_defhatxi}. 
    
    By comparing \eqref{eq_defhatxi} with \eqref{eq_pathX}, it is clear that the PDF $P_{x+\hat{\xi}_h}(\cdot)$ equals the conditional PDF $z \mapsto P_{X_{n+1}^m}\br{z \vert X_n^m = x}$ of $X_{n+1}^m$ given $X_n^m = x$. 
    Consequently, 
    \begin{equation}\label{eq_PDFjoint}
        P_{x+\hat{\xi}_h}\br{z}\, P_{X_n}\br{x} = P_{(X_n, X_{n+1})}\br{x, z}, \quad x, z \in \R^d,
    \end{equation}
    where $P_{(X_n, X_{n+1})}$ denotes the joint PDF of $(X_n^m, X_{n+1}^m)$. 
    As in \eqref{eq_defMtp}, we omit the superscript $m$ since $P_{(X_n, X_{n+1})}$ does not depend on $m$.
    Inserting \eqref{eq_PDFjoint} into \eqref{eq_estMnV2}, we obtain
    \begin{equation}\label{eq_estMnp1}
    \begin{aligned}
        \mathbb{M}_n^2\rbr{\Delta V_2} \leq\;& \int_{\R^d} \int_{\R^d} \abs{\Delta v(t_{n+1}, z)}^2 P_{(X_n, X_{n+1})}\br{x, z} \di z \di x\\
        =\;& \int_{\R^d} \abs{\Delta v(t_{n+1}, z)}^2  P_{X_{n+1}}\br{z}\di z\\
        =\;& \mathbb{M}_{n+1}^2\rbr{\Delta v}. 
    \end{aligned}
    \end{equation}
    Inserting \eqref{eq_estMnp1} into \eqref{eq_estMn2}, we obtain 
    \begin{equation}\label{eq_estMn2DelV}
        \mathbb{M}_n^2\rbr{\Delta v} \leq (1 + C_{\mathrm{st2}} h) \br{\mathbb{M}_{n+1}^2\rbr{\Delta v} + \mathbb{I}_{4, n}}.
    \end{equation}
    
    \subsubsection*{(iv) Estimate of $\mathbb{M}_n^2\rbr{\Delta v}$. }
    
    Let $k \in \{0, 1, \cdots, N-1\}$ be arbitrary.
    Sum \eqref{eq_estMn2DelV} over $n = k, k+1, \cdots, N-1$, leading to
    \begin{equation*}
        \sum_{n=k}^{N-1} \mathbb{M}_n^2\rbr{\Delta v} \leq (1 + C_{\mathrm{st2}} h) \vbr{\Big}{\sum_{n=k+1}^{N} \mathbb{M}_{n}^2\rbr{\Delta v} + \sum_{n=k}^{N-1} \mathbb{I}_{4, n}}.    
    \end{equation*}
    Canceling $\sum_{n=k+1}^{N-1} \mathbb{M}_n^2\rbr{\Delta v}$ on both sides, we obtain
    \begin{equation}\label{eq_Mk2Delv}
        \mathbb{M}_k^2\rbr{\Delta v} \leq C_{\mathrm{st2}} h \sum_{n=k+1}^{N} \mathbb{M}_{n}^2\rbr{\Delta v} + \mathbb{M}_{N}^2\rbr{\Delta v} + \mathbb{I}_5,
    \end{equation}
    where $\mathbb{I}_5 := (1 + C_{\mathrm{st2}} h) \sum_{n=0}^{N-1} \mathbb{I}_{4, n}$ satisfies
    \begin{equation}\label{eq_estI5}
    \begin{aligned}
        \mathbb{I}_5 \leq\;& (1 + C_{\mathrm{st2}} h) \sum_{n=0}^{N-1} \vbr{\Big}{3 h \mathbb{M}_n^2\rbr{\Delta R} + 96 C_{\mathrm{loc}}^2 (1 + C_{\mathrm{EM}}) h^3} \\
        \leq\;& 3 (1 + C_{\mathrm{st2}}) h \sum_{n=0}^{N-1} \mathbb{M}_n^2\rbr{\Delta R} + 96 C_{\mathrm{loc}}^2 (1 + C_{\mathrm{EM}}) (1 + C_{\mathrm{st2}}) T h^2.
    \end{aligned}
    \end{equation}
    Here the first line follows from inserting the bound in \eqref{eq_defI4n}, and the second line simplifies the expression using $h \leq 1$.
    
    Applying the backward discrete Gronwall
    inequality \cite[Lemma 3]{zhao2010stable} to \eqref{eq_Mk2Delv}, we obtain
    \begin{equation}\label{eq_Gronwall}
    \begin{aligned}
        \mathbb{M}_k^2\rbr{\Delta v} \leq\;& \exp\br{C_{\mathrm{st2}} T} \bbr{(1 + C_{\mathrm{st2}} h)\mathbb{M}_{N}^2\rbr{\Delta v} + \mathbb{I}_5} \\
        \leq\;& C_{\mathrm{st1}} \exp\br{C_{\mathrm{st2}} T} \Big\{\mathbb{M}_{N}^2\rbr{\Delta v} + h \sum_{n=0}^{N-1} \mathbb{M}_n^2\rbr{\Delta R} + T h^2 \Big\},
    \end{aligned}
    \end{equation}
    where the second line follows from \eqref{eq_estI5} with 
    \begin{equation}\label{eq_Cst1}
        C_{\mathrm{st1}} := \max\bbr{1 + C_{\mathrm{st2}},\, 3 (1 + C_{\mathrm{st2}}),\, 96 C_{\mathrm{loc}}^2 (1 + C_{\mathrm{EM}}) (1 + C_{\mathrm{st2}})}. 
    \end{equation}
    Here $C_{\mathrm{loc}}$, $C_{\mathrm{EM}}$, and $C_{\mathrm{st2}}$ are given in \eqref{eq_defCD}, \Cref{assum_em}, and \eqref{eq_estCh}, respectively.
    As the right-hand side of \eqref{eq_Gronwall} is independent of $k$, \eqref{eq_estMn2Delv} follows.
    
\end{proof}

\subsection{Global truncation error analysis}

Combining \Cref{coro_consist} and \Cref{lem_priest}, we can obtain the following convergence result of the formulation \eqref{eq_ER0}.

\begin{theorem}\label{thm_convergence}
    Under \Cref{assum_polygrow,assum_Lip,assum_em}, it holds that 
    \begin{equation}\label{eq_converge2}
        \max_{0 \leq n\leq N-1} \mathbb{M}_n^2\rbr{v - \hat{v}} \\
        \leq C_1 \exp\br{C_2 T} \Big(\mathbb{M}_{N}^2\rbr{v - \hat{v}} + L_{\mathrm{rdm}, \pi}(\hat{v}) + T h^2\Big),
    \end{equation}
    for $0 < h \leq \min\{1, (24 C_{\mathrm{nl}}^2 + 12 C_f^2)^{-1}\}$, 
    where $\mathbb{M}_n^2\rbr{\,\cdot\,}$ is defined in \eqref{eq_defMtp}; 
    $L_{\mathrm{rdm}, \pi}(\hat{v})$ is the RDM loss given in \eqref{eq_Lrdm} discretized on the time grid given in \eqref{eq_timepart}, i.e.,
    \begin{equation}\label{eq_defL2v}
        L_{\mathrm{rdm}, \pi}(\hat{v}) := \frac{T}{N}\sum_{n=0}^{N-1} \int_{\R^d} \big|\E{R(t_n, x, \xi; \hat{v})} \big|^2 P_{X_{n}}\br{x} \di x 
    \end{equation}
    with $x \mapsto P_{X_n}(x)$ the PDF introduced in \eqref{eq_defMtp}, and $R$ defined in \eqref{eq_defR};
    $C_1$ and $C_2$ are constants independent of $h$, $N$, and $\hat{v}$; see \eqref{eq_C1C2} for their expressions. 
\end{theorem}

\begin{proof}

    Inserting \eqref{eq_defdeltav} and \eqref{eq_defRi} into \eqref{eq_estMn2Delv} in \Cref{lem_priest}, we obtain
    \begin{equation}\label{eq_estglo}
        \max_{0 \leq n\leq N-1} \mathbb{M}_n^2\rbr{v - \hat{v}} \leq C_{\mathrm{st1}} \exp\br{C_{\mathrm{st2}} T} \Big\{\mathbb{M}_{N}^2\rbr{v - \hat{v}} + \mathbb{I}_6 + T h^2 \Big\},
    \end{equation}
    where 
    \begin{equation}\label{eq_defI6}
    \begin{aligned}
        \mathbb{I}_6 :=\;& \frac{T}{N}\sum_{n=0}^{N-1} \int_{\R^d} \big|\E{R(t_n, x, \xi; v) - R(t_n, x, \xi; \hat{v})}\big|^2 P_{X_{n}}\br{x} \di x \\
        \leq\;& 2 L_{\mathrm{rdm}, \pi}(v) + 2 L_{\mathrm{rdm}, \pi}(\hat{v}). 
    \end{aligned}
    \end{equation}
    Here the second line follows from $\abs{a - b}^2 \leq 2 \abs{a}^2 + 2 \abs{b}^2$ for any $a, b \in \R$, together with the definition in \eqref{eq_defL2v}.
    
    To estimate $L_{\mathrm{rdm}, \pi}(v)$, we substitute \eqref{eq_pde} into \eqref{eq_defR} to eliminate the term $f(t, x, v(t, x))$, yielding
    \begin{align}
        \abs{\E{R(t, x, \xi; v)}} =\;& \Big|\mathbb{E}\Big[ \frac{v(t + h, x + \xi_h) - v(t, x)}{h} \Big] - \mathcal{D} v(t, x)\Big| \notag \\
        =\;& \abs{\mathcal{D}_h v(t, x) - \mathcal{D} v(t, x)} \leq C_{\mathrm{loc}}\, h \big(1 + |x|^{\bar{p}}\big),\label{eq_ERtoDh}
    \end{align}
    where $\mathcal{D}_h$ is given by \eqref{eq_rdo}; 
    the last inequality follows from \Cref{coro_consist}. 
    Inserting \eqref{eq_ERtoDh} into \eqref{eq_defL2v} with $v$ in place of $\hat{v}$, we obtain
    \begin{equation}\label{eq_Lrdmv}
    \begin{aligned}
        L_{\mathrm{rdm}, \pi}(v) \leq\;& C_{\mathrm{loc}}^2 \frac{h^2 T}{N}\sum_{n=0}^{N-1} \int_{\R^d} \big(1 + |x|^{\bar{p}}\big)^2 P_{X_{n}}\br{x} \di x \\
        \leq\;& C_{\mathrm{loc}}^2 \frac{h^2 T}{N} \sum_{n=0}^{N-1} \br{2 + 2 C_{\mathrm{EM}}} = 2 C_{\mathrm{loc}}^2 (1 + C_{\mathrm{EM}}) T h^2,
    \end{aligned}
    \end{equation}
    
    Inserting \eqref{eq_Lrdmv} into the second line of \eqref{eq_defI6} and then substituting the resulting inequality into the right side of \eqref{eq_estglo}, we obtain the desired result \eqref{eq_converge2} with 
    \begin{equation}\label{eq_C1C2}
        C_1 := C_{\mathrm{st1}} \max\bbr{1,\, 2 + 4 C_{\mathrm{loc}}^2 (1 + C_{\mathrm{EM}})}, \quad C_2 := C_{\mathrm{st2}},
    \end{equation}
    where the constants $C_{\mathrm{st1}}$, $C_{\mathrm{st2}}$, $C_{\mathrm{loc}}$, and $C_{\mathrm{EM}}$ are given in \eqref{eq_Cst1}, \eqref{eq_estCh}, \eqref{eq_defCD}, and \Cref{assum_em}, respectively.

\end{proof}

We have the following discussions on the three error terms on the right-hand side of \eqref{eq_converge2}:

\begin{itemize}
    \item
    The first term $\mathbb{M}_{N}^2\rbr{v - \hat{v}}$ quantifies the error due to approximating the terminal condition $v(T, x) = g(x)$. 
    This term vanishes if the approximator $\hat{v}$ is the neural network $v_{\theta}$ constructed in \eqref{eq_netv}, which enforces the terminal condition exactly.
    
    \item The second term $L_{\mathrm{rdm}, \pi}(\hat{v})$, defined by \eqref{eq_defL2v}, represents the mean squared residual of $\hat{v}$ with respect to the formulation \eqref{eq_ER0}, weighted by the PDF $x \mapsto P_{X_{n}}\br{x}$. 
    This error can be effectively reduced by adversarially training the neural network $v_{\theta}$ according to \eqref{eq_minmax}. 
    
    \item The third term $T h^2$ corresponds to the truncation error arising from the random difference approximation~\eqref{eq_rdo}.
    This result indicates that our method achieves first-order accuracy in time; see \Cref{rmk_cr} for further discussion.
    
    \item Recalling \eqref{eq_defMtp}, the mean square error $\mathbb{M}_n^2\rbr{v - \hat{v}}$ is taken with respect to the PDF of $X_n^m$ generated by \eqref{eq_pathX}.
    Thus the result in \eqref{eq_converge2} indicates that the accuracy of $\hat{v}$ is ensured within the time-space region explored by the sample paths $\{X_n^m\}_{n=0}^{N-1}$. 
    This theoretical result is also supported by our numerical experiments; see the third column of \Cref{fig_HP9abcd}.
\end{itemize}

\begin{remark}[Convergence rate comparison]\label{rmk_cr}
As discussed in \cref{sec_con_mart}, the RDM formulation~\eqref{eq_rdo} is equivalent to the discretized martingale formulation~\eqref{eq0_mart_cond3}.
Therefore, \Cref{thm_convergence} also establishes a first-order convergence rate in time for the martingale methods in \cite{cai2024socmartnet,cai2023deepmartnet,cai2024martnet}.
This rate is higher than that of DeepBSDE-type methods \cite{weinan2017deep,han2018solving,raissi2018forwardbackward,hure2020deep,Zhou2021Actor,Zhang2022FBSDE}, which typically achieve at most order $1/2$ in time for general coefficients in \eqref{eq_pde}; see the theoretical analyses in \cite{han2018convergence,hure2020deep,germain2022Approximation,Andersson2023Convergence} and the numerical results in \cite{Zhang2022FBSDE}.
This discrepancy arises from different uses of the Euler-Maruyama scheme (EMS) in SDE discretization.
In martingale methods, the EMS is used in \eqref{eq_EulerApprox} to approximate the conditional law of $X_{t+h}$ given $X_t=x$; consequently, the error is governed by the weak convergence order of the EMS, namely $\mO(h)$ in sense of global truncation error \cite[Theorem 14.1.5]{Kloeden1992Numerical}.
In contrast, DeepBSDE-type methods use the EMS to approximate sample paths of the SDE in \eqref{eq_bsde}, and are therefore limited by the strong convergence rate $\mO(h^{1/2})$ of the EMS \cite[Theorem 10.2.2]{Kloeden1992Numerical}.
\end{remark}

\section{Derivation of variance of mini-batch estimation}\label{sec_variance}

In this section, we introduce the unbiased mini-batch estimators for the RDM loss in \eqref{eq_Lrdm} and the minimax loss in \cref{eq_minmax}. 
We then derive the variances of these mini-batch estimators, leading to \eqref{eq_varLrdm} and \eqref{eq_VarhatL2}. 
These results clarify the motivation for introducing the minimax formulation in \cref{eq_minmax}.

A mini-batch version of the RDM loss in \eqref{eq_Lrdm} can be defined as 
\begin{equation}\label{eq_Lrdm_mini}
    \hat{L}_{\mathrm{rdm}}(\hat{v}) := \frac{1}{M} \sum_{m=1}^{M} R_m^2,
\end{equation}
where $R_m^2$ is an unbiased mini-batch estimation of $\abs{\E{R(t, x, \xi; \hat{v})}}^2$ at $(t, x) = (t_m, x_m)$, given by
\begin{align}
    &R_m^2 := \br{\frac{1}{K} \sum_{k=1}^K R_{m, k}} \times \br{\frac{1}{K} \sum_{k=K+1}^{2K} R_{m, k}},\label{eq_defRnm} \\ 
    &R_{m, k} := R(t_{m}, x_{m}, \xi_{m, k}; \hat{v}). \label{eq_defRmk}
\end{align} 
Here, with a slight abuse of notation, $\{(t_{m}, x_{m})\}_{m=1}^M$ denotes i.i.d. samples drawn from the joint PDF $\tilde{p}$ on $[0, T-h] \times \R^d$, defined as
\begin{equation}\label{eq_def_bar_p}
    \tilde{p}(t, x) := \frac{p(t, x)}{T-h},\quad (t, x) \in [0, T-h] \times \R^d,
\end{equation}
where $p$ is the weight function specified in \eqref{eq_Lrdm}. 
The variables $\xi_{m, k}$, for all $m, k$, are i.i.d. samples of $\xi$.
Here, the subscript $m$ of $\xi_{m, k}$ indicates that we use different samples of $\xi$ for each individual $(t_{m}, x_{m})$. 
This strategy reduces the variance of the mini-batch estimator $\hat{L}_{\mathrm{rdm}}(\hat{v})$ compared to using a shared sample set of $\xi$ for all $(t_{m}, x_{m})$.

By the independence of $R_m$ across different $m$, we can compute the variance of $\hat{L}_{\mathrm{rdm}}(\hat{v})$ as follows:
\begin{equation}\label{eq_VarLrdm}
    \mathrm{Var}\rbr{\hat{L}_{\mathrm{rdm}}(\hat{v})} = \frac{1}{M^2} \sum_{m=1}^{M} \mathrm{Var}\rbr{R_m^2} = \frac{1}{M} \mathrm{Var}\rbr{R_1^2},
\end{equation}
where the second equality follows from the fact that all $R_m^2$ share the same distribution as the first sample $R_1^2$.
By the property of conditional variance, we can further decompose the variance of $R_1^2$ as follows:
\begin{equation}\label{eq_decompVar}
    \mathrm{Var}\rbr{R_1^2} = \mathrm{Var}\rbr{\E{R_1^2 \big \vert x_1}} + \E{\mathrm{Var}\rbr{R_1^2 \big \vert x_1}},
\end{equation}
where $\mathrm{Var}\rbr{R_1^2 \big \vert x_1}$ denotes the variance of $R_1^2$ conditioned on the sampled point $x_1$.
Using the definition of $R_1^2$ given by \eqref{eq_defRnm}, we have
\begin{align}
    \mathrm{Var}\rbr{\E{R_1^2 \big \vert x_1}} =\;& \mathrm{Var}\rbr{\frac{1}{K^2} \;\sum_{k=1}^K \sum_{l=K+1}^{2K} \E{R_{1, k} R_{1, l} \big \vert x_1}} \notag \\
    =\;& \mathrm{Var}\rbr{\E{R_{1, 1} R_{1, 2} \big \vert x_1}}, \label{eq_VarE}
\end{align} 
where the last equality follows from the fact that all $R_{1, k} R_{1, l}$ share the same distribution as $R_{1, 1} R_{1, 2}$.
Similarly,  
\begin{align}
    \E{\mathrm{Var}\rbr{R_1^2 \big \vert x_1}} =\;& \E{\mathrm{Var}\rbr{\left.\frac{1}{K^2} \sum_{k=1}^K \;\sum_{l=K+1}^{2K} R_{1, k} R_{1, l} \right\vert x_1}} \notag\\
    =\;& \E{\frac{1}{K^4} \sum_{k=1}^K \;\sum_{l=K+1}^{2K} \mathrm{Var}\rbr{R_{1, k} R_{1, l} \big\vert x_1}}, \notag \\
    =\;& \frac{1}{K^2} \E{\mathrm{Var}\rbr{R_{1, 1} R_{1, 2} \big\vert x_1}}, \label{eq_Evar}
\end{align}
where the second equality follows from the conditional independence of $R_{1, k} R_{1, l}$ for different $k,l$ given $x_1$.

Substituting \eqref{eq_VarE} and \eqref{eq_Evar} into \eqref{eq_decompVar}, and then inserting the result into \eqref{eq_VarLrdm}, we obtain
\begin{equation}\label{eq_VarLrdm_final}
    \mathrm{Var}\rbr{\hat{L}_{\mathrm{rdm}}(\hat{v})} = \frac{1}{M} \br{C_1 + \frac{C_2}{K^2}},
\end{equation}
where $C_1 := \mathrm{Var}\rbr{\E{R_{1, 1} R_{1, 2} \big \vert x_1}}$ and $C_2 := \E{\mathrm{Var}\rbr{R_{1, 1} R_{1, 2} \big\vert x_1}}$ are both constants independent of $M$ and $K$.
The result in \eqref{eq_VarLrdm_final} leads to \eqref{eq_varLrdm}, with a computational cost of $2 MK$ evaluations of the residual $R(t, x, \xi; \hat{v})$ in \eqref{eq_Lrdm}.


Next, we derive the variance of the mini-batch version of the minimax loss introduced in \eqref{eq_minmax}. 
Building on the notation used for \eqref{eq_Lrdm_mini}, we define the mini-batch estimator of $L(\hat{v}, \rho)$ as
\begin{equation}\label{eq_hatL}
    \abs{\hat{L}(\hat{v}, \rho)}^2 := \vbr{\Big}{\frac{1}{MK} \sum_{m=1}^{MK} \rho_m R_{m, 1}}^{\top} \vbr{\Big}{\frac{1}{MK} \sum_{m=MK+1}^{2MK} \rho_m R_{m, 1}},
\end{equation}
where $\rho_m := \rho(t_m, x_m)$ and $R_{m, 1}$ is defined in \eqref{eq_defRmk}, with the index $k$ of $\xi_{m, k}$ fixed to $1$ for all $m$. 
We emphasize that in \eqref{eq_hatL}, the index $m$ runs from $1$ to $2 MK$ rather than $2 M$. This ensures that the computational cost of evaluating $|\hat{L}(\hat{v}, \rho)|^2$ is comparable to that of \eqref{eq_Lrdm_mini}, both requiring $2 MK$ evaluations of the residual $R(t, x, \xi; \hat{v})$.
Using a similar analysis as in \eqref{eq_VarLrdm}, we can show that the variance of $\big|\hat{L}(\hat{v}, \rho)\big|^2$ is given by
\begin{equation}\label{eq_VarLrdm_minmax}
    \mathrm{Var}\rbr{\big|\hat{L}(\hat{v}, \rho)\big|^2} = \frac{1}{M K} \mathrm{Var}\rbr{\abs{\rho_1}^2 R_{1, 1}^2} = \mO\vbr{\Big}{\frac{1}{M K}}, 
\end{equation}
which is just \eqref{eq_VarhatL2}. 

\section{Cole-Hopf transformation for HJB equations}\label{sec_constPDE}

This section presents the derivation of the analytic solution in \eqref{eq_vtx_hjb}, using the Cole-Hopf transformation. 
The infimum in \eqref{eq_hjb_test} is attained at $\kappa = - c \sigma^{\top} \partial_x v$, and thus \eqref{eq_hjb_test} can be rewritten as
\begin{equation}\label{eq_pde_ch}
    \mathcal{R}v(t, x) = 0, \quad (t, x) \in [0, T) \times \R^d,
\end{equation}
where $\mathcal{R}$ is a differential operator defined by
\begin{equation}\label{eq_defRv}
    \mathcal{R}v := \vbr{\Big}{\partial_t + b^{\top} \partial_x + \frac{1}{2} \mathrm{Tr}[\sigma \sigma^{\top} \partial_{xx}]} v - \frac{c^2}{2} \big|\sigma^{\top} \partial_x v\big|^2. 
\end{equation}

For fixed $(s, x) \in [0, T) \times \R^d$, let $\X: [s, T] \times \Omega \to \R^d$ be the solution to the following SDE:
\begin{equation}\label{eq_sdeXtx}
    \X_t = x + \int_s^t b(r, \X_r) \di r + \int_s^t \sigma(r, \X_r) \di B_r, \quad t \in [s, T].  
\end{equation}
With $v$ being the solution to \eqref{eq_pde_ch}, define the process 
\begin{equation}\label{eq_defYt}
    Y_t := \exp\br{-c^2 v(t, \X_t)}, \quad t \in [s, T].
\end{equation}
Direct calculations show that
\begin{align*}
    \partial_{x_i} Y_t =\;& -c^2 Y_t \frac{\partial v}{\partial x_i}(t, \X_t), \\
    \partial_{x_j x_i}^2 Y_t =\;& -c^2 Y_t \frac{\partial^2 v}{\partial x_j \partial x_i}(t, \X_t) + c^4 Y_t \frac{\partial v}{\partial x_i} \frac{\partial v}{\partial x_j}(t, \X_t), \\
    \mathrm{Tr}\rbr{\sigma\sigma^{\top} \partial_{xx} Y_t} =\;& -c^2 Y_t \sum_{ijk} \sigma_{ij} \sigma_{kj} \frac{\partial^2 v}{\partial x_i \partial x_k}(t, \X_t) \\
    &+ c^4 Y_t \sum_{ijk} \sigma_{ij} \sigma_{kj} \frac{\partial v}{\partial x_i} \frac{\partial v}{\partial x_k}(t, \X_t).
\end{align*}
Applying the It\^o formula \cite[Theorem 2.3.3]{Zhang2017Backward} to the process $t \mapsto Y_t$, and substituting the above derivative terms into the resulting equation, we obtain
\begin{equation}\label{eq_diYt}
    Y_s = Y_T + \int_s^T c^2 Y_t \mathcal{R}v(t, \X_t) \di t + \int_s^T \bbr{\cdots} \di B_t,
\end{equation}
where we have used the definition of $\mathcal{R}v$ in \eqref{eq_defRv} to simplify the first integral term and $\bbr{\cdots}$ denotes the integrand of the It\^o integral, whose explicit form is omitted as it does not affect the remaining derivation.

Substituting \eqref{eq_pde_ch} into \eqref{eq_diYt}, the first integral term on the right-hand side vanishes. 
Taking expectations on both sides of the resulting equation yields
\begin{equation*}
    \E{Y_s} = \E{Y_T}, \quad s \in [0, T], 
\end{equation*}
where the It\^o integral vanishes by the martingale property under standard conditions on the integrand; see \cite[Theorem 2.2.7]{Zhang2017Backward}.
Substituting \eqref{eq_defYt} into the above relation and using $\X_s = x$ implied by \eqref{eq_sdeXtx}, we get $\exp\br{-c^2 v(s, x)} = \E{\exp\br{-c^2 g(\X_T)}}$, and hence
\begin{equation}\label{eq_vtx_hc}
\begin{aligned}
    v(s, x) &= -\frac{1}{c^2} \ln \E{\exp\br{-c^2 g(\X_T)}} \\
    &= -\frac{1}{c^2}\ln\br{\E{\exp\br{-c^2 v(T, X_T^0)} \big\vert X_s^0 = x}},
\end{aligned}
\end{equation}
where the second equality follows from the fact that for $t \in [s, T]$, the law of $\X_t$ defined by \eqref{eq_sdeXtx} coincides with the conditional law of $X_t^{0}$ given by \eqref{eq_xtu} under the zero control $u \equiv 0$ and the condition $X_s^{0}=x$.
The second line of \eqref{eq_vtx_hc} is just \eqref{eq_vtx_hjb}.

\section{Parameter setting for numerical experiments}\label{sec_param_setting}

\subsection{Parameter setting for \cref{sec_cde,sec_qlpde,sec_hjb,sec_ablation}}\label{sec_param_usual}

If not otherwise specified, we take $N = 100$ for all involved loss functions, and all loss functions are minimized by the RMSProp algorithm. 
The learning rates of \Cref{alg_pde,alg_amnet} are set to $\delta_0 = \delta_1 = 0.1$, $\delta_2 = 3 \times 10^{-3 - 2 i/I} \times d^{-r}$ for $i = 0, 1, \cdots, I-1$ with $r = 0.5$ for $d \leq 1,000$ and $r = 0.8$ for $d > 1,000$.
The inner iteration steps are $J = 2$ and $K = 1$.
The number $I$ of outer iterations is reported in each experiment.
For the experiments in \cref{sec_cde,sec_hjb,sec_ablation}, 
the time step size in \eqref{eq_rdo} and \eqref{eq_timepart} is set to $h = T/N$ with $N = 100$. 
For the experiments in \cref{sec_qlpde}, we set $N = 32$ to reduce memory usage.

The adversarial network $\rho_{\eta}$ is constructed according to \eqref{eq_defrhonet} with $r = 600$.
The parameter $c$ in \eqref{eq_defLamb} is set to $c = 99/(r-1)$ for the experiments in \cref{sec_cde,sec_hjb,sec_ablation}, and to $c = 9/(r-1)$ for the experiments in \cref{sec_qlpde}. 
The architecture of $v_{\theta}$ follows \eqref{eq_netv} to ensure satisfaction of the terminal condition. 
The neural networks $u_{\alpha}$ and $\phi_{\theta}$ in \eqref{eq_netv} are fully connected feedforward neural networks with $6$ hidden layers and ReLU activation functions.
For the experiments in \cref{sec_cde,sec_hjb,sec_ablation},  
the number of ReLU units in each hidden layer is set to $\max\{100, k d + 10\}$,  
where $k = 5$ for $d \leq 100$, $k = 2$ for $100 < d \leq 1000$, and $k = 1$ for $d > 1000$. 
For the experiments in \cref{sec_qlpde}, we set $W = 10000$ for both $d = 10^4$ and $d = 10^5$ to reduce memory usage. 
For HJB-2 in \cref{sec_hjb}, we also adopt MscaleDNNs \cite{Liu2020Multi} to construct the networks $u_{\alpha}$ and $v_{\theta}$; see \eqref{eq_mscaleDNN}.

The start points $X_0^{m}$ of the sample paths in \eqref{eq_pathX} and \eqref{eq_defhatXpil0_hjb} are randomly taken from the set $D_0$ with replacement, where the total number of sample paths retained in memory is set to $M = 10^4$, $5000$ and $2000$ for $d \leq 5000$, $5000 < d \leq 10^4$, and $d > 10^4$, respectively.
The mini-batch sizes in \eqref{eq_defA1A2} are set to $\abs{\mathbb{M}_1} = \abs{\mathbb{M}_2} = 128$, $64$, $32$ and $16$ for $d \leq 1000$, $1000 < d \leq 5000$, $5000 < d \leq 10^4$, and $d > 10^4$, respectively.
The index sets $\mathbb{N}_1$ and $\mathbb{N}_2$ are both fixed as $\bbr{0, 1, \cdots, N-1}$.
In \Cref{alg_pde,alg_amnet}, the interval $I_0$ for updating spatial samples is chosen as $M/(\abs{\mathbb{A}_1} + \abs{\mathbb{A}_2})$, and the percentage $r\%$ of sample paths updated is set to $20\%$.
The distribution $P_{\xi}$ is chosen to be the standard Gaussian $\mathrm{N}(0, I_q)$.

All numerical experiments are conducted using Python 3.12 and PyTorch 2.7.0.
The networks $u_{\alpha}$ and $v_{\theta}$ are trained using the Automatic Mixed Precision technique\footnote{\url{https://docs.pytorch.org/docs/2.7/amp.html}}.
The algorithms are accelerated using Distributed Data Parallel\footnote{\url{https://github.com/pytorch/tutorials/blob/main/intermediate_source/ddp_tutorial.rst}} on a compute node equipped with 8 NVIDIA A100-SXM4-80GB GPUs.

\subsection{Parameter setting for \cref{sec_compare}}\label{sec_param_comp}

\begin{algorithm}[t]
    \caption{
    Strong-form DRDM for the quasilinear parabolic equation \eqref{eq_pde}
    }\label{alg_pde_strf}
    \begin{algorithmic}[1]
        \Require 
        $v_{\theta}$: neural network parameterized by $\theta$;
        $\delta_1$: learning rate for $v_{\theta}$;
        $I$: maximum number of iterations;
        $I_0$ and $r\%$: interval for updating sample paths and the percentage of paths updated;
        $\ell$: number of random samples of $\xi$ to approximate the expectation in \eqref{eq_rdo};
        $M$: total number of sample paths retained in memory;
        $P_{\xi}$: distribution of $\xi$ satisfying \eqref{eq_condi_xi}. 
        
        \State Initialize $v_{\theta}$
        \State Generate the sample paths $\{X_n^{m}: 0 \leq n \leq N, 1 \leq m \leq M\}$ by \eqref{eq_pathX} with $\hat{v} = v_{\theta}$
        \For{$i = 1, 2, \cdots, I-1$}

        \State Sample the index subset $\mathbb{A} \subset \{0, 1, \cdots, N-1\} \times \{1, 2, \cdots, M\}$. 
        
        \State Sample $\{\xi_n^{m, j} : (n, m) \in \mathbb{A},\; 1 \leq j \leq \ell\}$ i.i.d. from $P_{\xi}$.    
        
        \State $\theta \leftarrow \theta - \delta_1 \nabla_{\theta} L_{\mathrm{rdm}}(v_{\theta}; \mathbb{A})$ with 
        \begin{equation*}
            L_{\mathrm{rdm}}(v; \mathbb{A}) := \frac{2 T h}{\ell \abs{\mathbb{A}}} \sum_{(n, m) \in \mathbb{A}} \,\sum_{j=1}^{\ell/2} R(t_n, X_n^m, \xi_n^{m, j}; v) R(t_n, X_n^m, \xi_n^{m, \ell/2 + j}; v)
        \end{equation*}
        
        \If{$i$ is divisible by $I_0$} 
        \State Update $r\%$ of the sample paths $\{X_n^m\}_{n=0}^N$ by \eqref{eq_pathX} with $\hat{v} = v_{\theta}$
        \EndIf
        
        \EndFor
        \Ensure $v_{\theta}$
    \end{algorithmic}
\end{algorithm}

We first provide the setting for the weak-form DRDM described in \Cref{alg_pde}.
As the benchmark problems \eqref{eq_ac} and \eqref{eq_sg} are time-independent, we should adapt \Cref{alg_pde} accordingly.
In particular,
the architecture of $v_{\theta}$ is time-independent and follows \cite[subsection 5.4.2]{Hu2025Bias}:
\begin{equation}
    v_{\theta}(x) = \vbr{\big}{1 - \abs{x}^2} \phi_{\theta}(x), \quad x \in \R^d,
\end{equation}
where the factor $1 - \abs{x}^2$ enforces the zero boundary condition; 
$\phi_{\theta}$ is a fully connected network with 4 hidden layers, 128 units per layer, and $\tanh$ activation.
The adversarial network $\rho_{\eta}$ is a time-independent variant of \eqref{eq_defrhonet}: 
\begin{equation*}
    \rho_{\eta}(x) = \sin \vbr{\big}{\Lambda \br{W x + b}} \in \R^r, \quad r = 128, \quad \eta = \br{W, b},
\end{equation*}
where $\Lambda(\cdot)$ is given by \eqref{eq_defLamb} with $c = 99 / (r - 1)$. 
In this time-independent setting, we use a step size of $h = 1/100$ for the random difference approximation in \eqref{eq_rdo} and simply set the number of time steps to $N = 1$.

The number $I$ of outer iterations in \Cref{alg_pde} is set to $5 \times 10^3$ for $d \leq 10^4$ and $12 \times 10^3$ for $d = 10^5$, with the inner iterations fixed at $J = 2$ and $K = 1$. 
For $d \leq 10^4$, this yields $10^4$ gradient descent steps in total, matching \cite[subsection 5.4.2]{Hu2025Bias}; 
for $d = 10^5$, we increase the number of outer iterations because DRDM does not reach satisfactory accuracy within $10^4$ gradient descent steps.
The networks $v_{\theta}$ and $\rho_{\eta}$ are trained with the Adam optimizer \cite{Kingma2017Adam} using a learning rate of $10^{-3 - 3 i / I}$ at the $i$-th outer iteration.
The number $M$ of sample paths retained in memory is set to $1000$.
The mini-batch size in \eqref{eq_defA1A2} is set to 100, i.e., set $\mathbb{A}_i = \{0, 1\} \times \mathbb{M}_i$ with $\abs{\mathbb{M}_i} = 50$ for $i = 1, 2$, to be consistent with \cite[subsection~5.4.2]{Hu2025Bias}. 
At the beginning of training, and whenever $\{X_0^m\}_{m=1}^M$ is updated, the samples $X_0^m$ are drawn in an i.i.d. manner by
\footnote{
In \cite[subsection 5.4.2]{Hu2025Bias}, it is stated that spatial samples are drawn ``uniformly'' from the unit ball $\mathbb{B}^d$.
According to the source code at \url{https://github.com/zheyuanhu01/SDGD_PINN}, ``uniformly'' means uniform in direction and radial distance as in \eqref{eq_unif_ball}, rather than uniform with respect to the Lebesgue measure on $\mathbb{B}^d$.
}
\begin{equation}\label{eq_unif_ball}
    X_0^m = u_m\, y_m / \abs{y_m}, \quad u_m \sim \mathrm{U}(0, 1), \quad y_m \sim \mathrm{N}(0, I_d).
\end{equation}
The distribution $P_{\xi}$ in \Cref{alg_pde} is chosen to be the standard Gaussian $\mathrm{N}(0, I_q)$.

The accuracy of each method is evaluated using the relative $L_2$ error defined as
\begin{equation*}
    \mathrm{Rel.~} L_2 \mathrm{~error} := \left(\frac{\sum_{i=1}^{k} \abs{v_{\theta}(x_i) - v(x_i)}^2}{\sum_{i=1}^{k} \abs{v(x_i)}^2}\right)^{1/2}, \quad k = 2 \times 10^4,
\end{equation*}
where $\{x_i\}_{i=1}^k \subset \mathbb{B}^d$ are the test points sampled from the same distribution as $X_0^m$ in \eqref{eq_unif_ball}, and fixed throughout training.
The algorithm is run on a single NVIDIA A100-SXM4-80GB GPU using FP32 precision, without Automatic Mixed Precision, to match the setup in \cite[subsection 5.4.2]{Hu2025Bias}.

We now present the configuration for the strong-form DRDM, whose procedure is outlined in \Cref{alg_pde_strf}.
For brevity, we only list the parameters that differ from the weak-form DRDM in \Cref{alg_pde}.
The number of iterations $I$ is set to $10{,}000$ for $d \leq 10^4$ and $24{,}000$ for $d = 10^5$.
This matches the total number of gradient descent steps used in \Cref{alg_pde} and in \cite[subsection 5.4.2]{Hu2025Bias}.
The mini-batch size $\abs{\mathbb{A}}$ is set to $100$, i.e., $\mathbb{A} = \{0, 1\} \times \mathbb{M}$ with $\mathbb{M} \subset \{1, 2, \cdots, M\}$ and $\abs{\mathbb{M}} = 100$; the number $\ell$ of samples of $\xi$ is set to $\ell = 128$, consistent with the RS-PINN setting in \cite[subsection 5.4.2]{Hu2025Bias}.
All other settings follow \Cref{alg_pde}.

\section{RDM for quadratic gradient operator}\label{sec_rdm_quad}

Now we introduce a random difference approximation for the quadratic gradient $\abs{\partial_x v}^2$.
To connect with the expansion in \cref{sec_basic}, we still consider a generic function $F: \R^q \to \R$ whose third-order derivatives are bounded.
Let $\xi$ be a $q$-dimensional random variable satisfying \eqref{eq_condi_xi}.
The Taylor expansion of $F$ at $z=0$ yields the pointwise counterpart of \eqref{eq_taylor}:
\begin{equation}\label{eq_taylor2}
    F(\sqrt{h} \xi) = \sum_{k=0}^2 \frac{h^{k/2}}{k!} \vbr{\big}{\xi^{\top} \partial_z}^k F(0) + \frac{h^{1.5}}{3!} \vbr{\big}{\xi^{\top} \partial_z}^3 F(c\sqrt{h} \xi),
\end{equation}
where $h > 0$ is a step size, and $c$ is a random variable taking values in $[0, 1]$ that depends on $\sqrt{h}\,\xi$.
To extract a gradient term from \eqref{eq_taylor2}, we multiply both sides by $\sqrt{h}\,\xi$ and take expectations, which yields
\begin{equation}\label{eq_taylor3}
\begin{aligned}
\E{\sqrt{h} \xi F(\sqrt{h} \xi)} =\;& \sum_{k=0}^2 \frac{h^{\frac{k+1}{2}}}{k!} \E{\xi \vbr{\big}{\xi^{\top} \partial_z}^k F(0)} \\
&+ \frac{h^2}{3!} \E{\xi \vbr{\big}{\xi^{\top} \partial_z}^3 F(c\sqrt{h} \xi)},
\end{aligned}
\end{equation}
where, by the moment conditions in \eqref{eq_condi_xi}, we have 
\begin{equation*}
\E{\xi \vbr{\big}{\xi^{\top} \partial_z}^k F(0)} = \left\{\begin{aligned}
    &0, && k = 0, 2,\\
    &\partial_z F(0), && k = 1,
\end{aligned}\right.
\end{equation*}
and $\E{\xi \vbr{\big}{\xi^{\top} \partial_z}^3 F(c\sqrt{h} \xi)} = \mathcal{O}(1)$ as $h \to 0$. 
Then the expansion \eqref{eq_taylor3} reduces to
\begin{equation*}
    \E{\sqrt{h} \xi F(\sqrt{h} \xi)} = h \partial_z F(0) + \mathcal{O}(h^2), 
\end{equation*}
which yields the RDM approximation 
\begin{equation}\label{eq_rdm_grad}
    \partial_z F(0) := \E{\frac{\xi}{\sqrt{h}} F(\sqrt{h} \xi)} + \mathcal{O}(h). 
\end{equation}

For the function $v: [0, T] \times \R^d \to \R$ and a fixed point $(t, x) \in [0, T] \times \R^d$, the gradient $\partial_x v(t, x)$ can be approximated by \eqref{eq_rdm_grad} with $F$ replaced by $z \mapsto v(t, x + z)$, i.e., 
\begin{equation}\label{eq_rdm_grad_v}
    \partial_x v(t, x) = \E{\frac{\xi}{\sqrt{h}} v(t, x + \sqrt{h} \xi)} + \mathcal{O}(h). 
\end{equation}
Moreover, using the debiasing technique in \cite{Hu2025Bias}, we can obtain an unbiased approximation to $\abs{\partial_x v}^2$ as
\begin{equation}\label{eq_rdm_grad2}
\begin{aligned}
    \abs{\partial_x v(t, x)}^2 =\;& \E{\frac{\xi^{\top} \bar{\xi}}{h} v(t, x + \sqrt{h} \xi) \, v(t, x + \sqrt{h} \bar{\xi})} + \mathcal{O}(h),
\end{aligned}
\end{equation}
where $\bar{\xi}$ is an independent copy of $\xi$.
We refer to the first terms on the right-hand sides of \eqref{eq_rdm_grad_v} and \eqref{eq_rdm_grad2} as the RDM approximations to $\partial_x v(t, x)$ and $\abs{\partial_x v(t, x)}^2$, respectively.

\begin{remark}[Connection to RS-PINNs]
Comparing \eqref{eq_rdm_grad_v} with the first identity in \eqref{eq_stein}, we observe that although RS-PINNs for $\partial_x v$ are constructed to represent the gradient of the Gaussian-smoothed network $v_{\theta}(x)$, they simultaneously provide a first-order approximation to the gradient $\partial_x \psi_{\theta}(x)$ of the original (unsmoothed) network.
This suggests that, within the RS-PINN framework \cite{He2023Learning,Hu2025Bias}, one may directly use $\psi_{\theta}$ itself as an approximator of $v$.
This modification incurs an $\mathcal{O}(h)$ truncation error for $\partial_x v$ but avoids the computational cost of evaluating the Gaussian-smoothing expectation for model inference.
\end{remark}

\bibliographystyle{elsarticle-num-names}
\bibliography{bibliography}

\begin{thebibliography}{55}
\expandafter\ifx\csname natexlab\endcsname\relax\def\natexlab#1{#1}\fi
\providecommand{\url}[1]{\texttt{#1}}
\providecommand{\href}[2]{#2}
\providecommand{\path}[1]{#1}
\providecommand{\DOIprefix}{doi:}
\providecommand{\ArXivprefix}{arXiv:}
\providecommand{\URLprefix}{URL: }
\providecommand{\Pubmedprefix}{pmid:}
\providecommand{\doi}[1]{\href{http://dx.doi.org/#1}{\path{#1}}}
\providecommand{\Pubmed}[1]{\href{pmid:#1}{\path{#1}}}
\providecommand{\bibinfo}[2]{#2}
\ifx\xfnm\relax \def\xfnm[#1]{\unskip,\space#1}\fi
\bibitem[{Al-Aradi et~al.(2022)Al-Aradi, Correia, Jardim, Naiff, and Saporito}]{Al2022Extensions}
\bibinfo{author}{A.~Al-Aradi}, \bibinfo{author}{A.~Correia}, \bibinfo{author}{G.~Jardim}, \bibinfo{author}{D.~d.~F. Naiff}, \bibinfo{author}{Y.~Saporito},
\newblock \bibinfo{title}{Extensions of the deep {G}alerkin method},
\newblock \bibinfo{journal}{Appl. Math. Comput.} \bibinfo{volume}{430} (\bibinfo{year}{2022}) \bibinfo{pages}{Paper No. 127287, 18}. \DOIprefix\doi{10.1016/j.amc.2022.127287}.
\bibitem[{Andersson et~al.(2023)Andersson, Andersson, and Oosterlee}]{Andersson2023Convergence}
\bibinfo{author}{K.~Andersson}, \bibinfo{author}{A.~Andersson}, \bibinfo{author}{C.~W. Oosterlee},
\newblock \bibinfo{title}{Convergence of a robust deep {FBSDE} method for stochastic control},
\newblock \bibinfo{journal}{SIAM J. Sci. Comput.} \bibinfo{volume}{45} (\bibinfo{year}{2023}) \bibinfo{pages}{a226--a255}. \DOIprefix\doi{10.1137/22M1478057}.
\bibitem[{Bachouch et~al.(2022)Bachouch, Hur\'{e}, Langren\'{e}, and Pham}]{Bachouch2022Deep}
\bibinfo{author}{A.~Bachouch}, \bibinfo{author}{C.~Hur\'{e}}, \bibinfo{author}{N.~Langren\'{e}}, \bibinfo{author}{H.~Pham},
\newblock \bibinfo{title}{Deep neural networks algorithms for stochastic control problems on finite horizon: numerical applications},
\newblock \bibinfo{journal}{Methodol. Comput. Appl. Probab.} \bibinfo{volume}{24} (\bibinfo{year}{2022}) \bibinfo{pages}{143--178}. \DOIprefix\doi{10.1007/s11009-019-09767-9}.
\bibitem[{Beck et~al.(2021)Beck, Becker, Cheridito, Jentzen, and Neufeld}]{beck2019deep}
\bibinfo{author}{C.~Beck}, \bibinfo{author}{S.~Becker}, \bibinfo{author}{P.~Cheridito}, \bibinfo{author}{A.~Jentzen}, \bibinfo{author}{A.~Neufeld},
\newblock \bibinfo{title}{Deep splitting method for parabolic {PDE}s},
\newblock \bibinfo{journal}{SIAM J. Sci. Comput.} \bibinfo{volume}{43} (\bibinfo{year}{2021}) \bibinfo{pages}{A3135--A3154}. \DOIprefix\doi{10.1137/19M1297919}.
\bibitem[{Cai et~al.(2024)Cai, Fang, Zhang, and Zhou}]{cai2024martnet}
\bibinfo{author}{W.~Cai}, \bibinfo{author}{S.~Fang}, \bibinfo{author}{W.~Zhang}, \bibinfo{author}{T.~Zhou}, \bibinfo{title}{Martingale deep learning for very high dimensional quasi-linear partial differential equations and stochastic optimal controls}, \bibinfo{year}{2024}. \href{http://arxiv.org/abs/2408.14395}{{\tt arXiv:2408.14395}}.
\bibitem[{Cai et~al.(2025)Cai, Fang, and Zhou}]{cai2024socmartnet}
\bibinfo{author}{W.~Cai}, \bibinfo{author}{S.~Fang}, \bibinfo{author}{T.~Zhou},
\newblock \bibinfo{title}{{SOC-MartNet}: A martingale neural network for the hamilton-jacobi-bellman equation without explicit $\inf_{u \in U} {H}$ in stochastic optimal controls},
\newblock \bibinfo{journal}{SIAM J. Sci. Comput.} \bibinfo{volume}{47} (\bibinfo{year}{2025}) \bibinfo{pages}{C795--C819}. \DOIprefix\doi{10.1137/24M1681033}.
\bibitem[{Cai et~al.(2023)Cai, He, and Margolis}]{cai2023deepmartnet}
\bibinfo{author}{W.~Cai}, \bibinfo{author}{A.~He}, \bibinfo{author}{D.~Margolis}, \bibinfo{title}{{DeepMartNet} -- a martingale based deep neural network learning method for {Dirichlet} {BVP} and eigenvalue problems of elliptic pdes}, \bibinfo{year}{2023}. \href{http://arxiv.org/abs/arXiv:2311.09456 [math.NA]}{{\tt arXiv:arXiv:2311.09456 [math.NA]}}.
\bibitem[{Chassagneux(2014)}]{Chassagneux2014Linear}
\bibinfo{author}{J.-F. Chassagneux},
\newblock \bibinfo{title}{Linear multistep schemes for {BSDE}s},
\newblock \bibinfo{journal}{SIAM J. Numer. Anal.} \bibinfo{volume}{52} (\bibinfo{year}{2014}) \bibinfo{pages}{2815--2836}. \DOIprefix\doi{10.1137/120902951}.
\bibitem[{Chassagneux and Crisan(2014)}]{Chassagneux2014Runge}
\bibinfo{author}{J.-F. Chassagneux}, \bibinfo{author}{D.~Crisan},
\newblock \bibinfo{title}{Runge-{K}utta schemes for backward stochastic differential equations},
\newblock \bibinfo{journal}{Ann. Appl. Probab.} \bibinfo{volume}{24} (\bibinfo{year}{2014}) \bibinfo{pages}{679--720}. \DOIprefix\doi{10.1214/13-AAP933}.
\bibitem[{Chen et~al.(2023)Chen, Huang, Wang, and Yang}]{Chen2023Friedrichs}
\bibinfo{author}{F.~Chen}, \bibinfo{author}{J.~Huang}, \bibinfo{author}{C.~Wang}, \bibinfo{author}{H.~Yang},
\newblock \bibinfo{title}{Friedrichs learning: Weak solutions of partial differential equations via deep learning},
\newblock \bibinfo{journal}{SIAM Journal on Scientific Computing} \bibinfo{volume}{45} (\bibinfo{year}{2023}) \bibinfo{pages}{A1271--A1299}. \DOIprefix\doi{10.1137/22M1488405}.
\bibitem[{E et~al.(2017)E, Han, and Jentzen}]{weinan2017deep}
\bibinfo{author}{W.~E}, \bibinfo{author}{J.~Han}, \bibinfo{author}{A.~Jentzen},
\newblock \bibinfo{title}{Deep learning-based numerical methods for high-dimensional parabolic partial differential equations and backward stochastic differential equations},
\newblock \bibinfo{journal}{Commun. Math. Stat.} \bibinfo{volume}{5} (\bibinfo{year}{2017}) \bibinfo{pages}{349--380}. \DOIprefix\doi{10.1007/s40304-017-0117-6}.
\bibitem[{E and Yu(2018)}]{E2018deep}
\bibinfo{author}{W.~E}, \bibinfo{author}{B.~Yu},
\newblock \bibinfo{title}{The deep {R}itz method: a deep learning-based numerical algorithm for solving variational problems},
\newblock \bibinfo{journal}{Commun. Math. Stat.} \bibinfo{volume}{6} (\bibinfo{year}{2018}) \bibinfo{pages}{1--12}. \DOIprefix\doi{10.1007/s40304-018-0127-z}.
\bibitem[{Gao et~al.(2023)Gao, Yan, and Zhou}]{Gao2023Failure}
\bibinfo{author}{Z.~Gao}, \bibinfo{author}{L.~Yan}, \bibinfo{author}{T.~Zhou},
\newblock \bibinfo{title}{Failure-informed adaptive sampling for {PINN}s},
\newblock \bibinfo{journal}{SIAM J. Sci. Comput.} \bibinfo{volume}{45} (\bibinfo{year}{2023}) \bibinfo{pages}{A1971--A1994}. \DOIprefix\doi{10.1137/22M1527763}.
\bibitem[{Germain et~al.(2022)Germain, Pham, and Warin}]{germain2022Approximation}
\bibinfo{author}{M.~Germain}, \bibinfo{author}{H.~Pham}, \bibinfo{author}{X.~Warin},
\newblock \bibinfo{title}{Approximation error analysis of some deep backward schemes for nonlinear {PDE}s},
\newblock \bibinfo{journal}{SIAM J. Sci. Comput.} \bibinfo{volume}{44} (\bibinfo{year}{2022}) \bibinfo{pages}{A28--A56}. \DOIprefix\doi{10.1137/20M1355355}.
\bibitem[{Guo et~al.(2022)Guo, Wu, Yu, and Zhou}]{Guo2022Monte}
\bibinfo{author}{L.~Guo}, \bibinfo{author}{H.~Wu}, \bibinfo{author}{X.~Yu}, \bibinfo{author}{T.~Zhou},
\newblock \bibinfo{title}{Monte {C}arlo f{PINN}s: deep learning method for forward and inverse problems involving high dimensional fractional partial differential equations},
\newblock \bibinfo{journal}{Comput. Methods Appl. Mech. Engrg.} \bibinfo{volume}{400} (\bibinfo{year}{2022}) \bibinfo{pages}{Paper No. 115523, 17}. \DOIprefix\doi{10.1016/j.cma.2022.115523}.
\bibitem[{Han et~al.(2018)Han, Jentzen, and E}]{han2018solving}
\bibinfo{author}{J.~Han}, \bibinfo{author}{A.~Jentzen}, \bibinfo{author}{W.~E},
\newblock \bibinfo{title}{Solving high-dimensional partial differential equations using deep learning},
\newblock \bibinfo{journal}{Proceedings of the National Academy of Sciences} \bibinfo{volume}{115} (\bibinfo{year}{2018}) \bibinfo{pages}{8505--8510}.
\bibitem[{Han and Long(2020)}]{han2018convergence}
\bibinfo{author}{J.~Han}, \bibinfo{author}{J.~Long},
\newblock \bibinfo{title}{Convergence of the deep {BSDE} method for coupled {FBSDE}s},
\newblock \bibinfo{journal}{Probab. Uncertain. Quant. Risk} \bibinfo{volume}{5} (\bibinfo{year}{2020}) \bibinfo{pages}{Paper No. 5, 33}. \DOIprefix\doi{10.1186/s41546-020-00047-w}.
\bibitem[{Han et~al.(2020{\natexlab{a}})Han, Lu, and Zhou}]{Han2020Solving}
\bibinfo{author}{J.~Han}, \bibinfo{author}{J.~Lu}, \bibinfo{author}{M.~Zhou},
\newblock \bibinfo{title}{Solving high-dimensional eigenvalue problems using deep neural networks: a diffusion {M}onte {C}arlo like approach},
\newblock \bibinfo{journal}{J. Comput. Phys.} \bibinfo{volume}{423} (\bibinfo{year}{2020}{\natexlab{a}}) \bibinfo{pages}{109792, 13}. \DOIprefix\doi{10.1016/j.jcp.2020.109792}.
\bibitem[{Han et~al.(2020{\natexlab{b}})Han, Nica, and Stinchcombe}]{Han2020derivative}
\bibinfo{author}{J.~Han}, \bibinfo{author}{M.~Nica}, \bibinfo{author}{A.~R. Stinchcombe},
\newblock \bibinfo{title}{A derivative-free method for solving elliptic partial differential equations with deep neural networks},
\newblock \bibinfo{journal}{J. Comput. Phys.} \bibinfo{volume}{419} (\bibinfo{year}{2020}{\natexlab{b}}) \bibinfo{pages}{18}. \DOIprefix\doi{10.1016/j.jcp.2020.109672}.
\bibitem[{He et~al.(2023)He, Li, Shi, Gao, Zhang, Bian, Wang, and Liu}]{He2023Learning}
\bibinfo{author}{D.~He}, \bibinfo{author}{S.~Li}, \bibinfo{author}{W.~Shi}, \bibinfo{author}{X.~Gao}, \bibinfo{author}{J.~Zhang}, \bibinfo{author}{J.~Bian}, \bibinfo{author}{L.~Wang}, \bibinfo{author}{T.-Y. Liu},
\newblock \bibinfo{title}{Learning physics-informed neural networks without stacked back-propagation},
\newblock in: \bibinfo{editor}{F.~Ruiz}, \bibinfo{editor}{J.~Dy}, \bibinfo{editor}{J.-W. van~de Meent} (Eds.), \bibinfo{booktitle}{Proceedings of The 26th International Conference on Artificial Intelligence and Statistics}, volume \bibinfo{volume}{206} of \textit{\bibinfo{series}{Proceedings of Machine Learning Research}}, \bibinfo{publisher}{PMLR}, \bibinfo{year}{2023}, pp. \bibinfo{pages}{3034--3047}.
\bibitem[{Hu et~al.(2024{\natexlab{a}})Hu, Shi, Karniadakis, and Kawaguchi}]{Hu2024Hutchinson}
\bibinfo{author}{Z.~Hu}, \bibinfo{author}{Z.~Shi}, \bibinfo{author}{G.~E. Karniadakis}, \bibinfo{author}{K.~Kawaguchi},
\newblock \bibinfo{title}{Hutchinson trace estimation for high-dimensional and high-order physics-informed neural networks},
\newblock \bibinfo{journal}{Comput. Methods Appl. Mech. Engrg.} \bibinfo{volume}{424} (\bibinfo{year}{2024}{\natexlab{a}}) \bibinfo{pages}{Paper No. 116883, 17}. \DOIprefix\doi{10.1016/j.cma.2024.116883}.
\bibitem[{Hu et~al.(2024{\natexlab{b}})Hu, Shukla, Karniadakis, and Kawaguchi}]{hu2024sdgd}
\bibinfo{author}{Z.~Hu}, \bibinfo{author}{K.~Shukla}, \bibinfo{author}{G.~E. Karniadakis}, \bibinfo{author}{K.~Kawaguchi},
\newblock \bibinfo{title}{Tackling the curse of dimensionality with physics-informed neural networks},
\newblock \bibinfo{journal}{Neural Networks} \bibinfo{volume}{176} (\bibinfo{year}{2024}{\natexlab{b}}) \bibinfo{pages}{106369}. \DOIprefix\doi{https://doi.org/10.1016/j.neunet.2024.106369}.
\bibitem[{Hu et~al.(2025)Hu, Yang, Wang, Karniadakis, and Kawaguchi}]{Hu2025Bias}
\bibinfo{author}{Z.~Hu}, \bibinfo{author}{Z.~Yang}, \bibinfo{author}{Y.~Wang}, \bibinfo{author}{G.~E. Karniadakis}, \bibinfo{author}{K.~Kawaguchi},
\newblock \bibinfo{title}{Bias-{V}ariance {T}rade-{O}ff in {P}hysics-{I}nformed {N}eural {N}etworks with {R}andomized {S}moothing for {H}igh-{D}imensional {PDE}s},
\newblock \bibinfo{journal}{SIAM J. Sci. Comput.} \bibinfo{volume}{47} (\bibinfo{year}{2025}) \bibinfo{pages}{C846--C872}. \URLprefix \url{https://doi.org/10.1137/23M1621356}. \DOIprefix\doi{10.1137/23M1621356}.
\bibitem[{Hur\'{e} et~al.(2021)Hur\'{e}, Pham, Bachouch, and Langren\'{e}}]{Hure2021Deep}
\bibinfo{author}{C.~Hur\'{e}}, \bibinfo{author}{H.~Pham}, \bibinfo{author}{A.~Bachouch}, \bibinfo{author}{N.~Langren\'{e}},
\newblock \bibinfo{title}{Deep neural networks algorithms for stochastic control problems on finite horizon: convergence analysis},
\newblock \bibinfo{journal}{SIAM J. Numer. Anal.} \bibinfo{volume}{59} (\bibinfo{year}{2021}) \bibinfo{pages}{525--557}. \DOIprefix\doi{10.1137/20M1316640}.
\bibitem[{Hur\'{e} et~al.(2020)Hur\'{e}, Pham, and Warin}]{hure2020deep}
\bibinfo{author}{C.~Hur\'{e}}, \bibinfo{author}{H.~Pham}, \bibinfo{author}{X.~Warin},
\newblock \bibinfo{title}{Deep backward schemes for high-dimensional nonlinear {PDE}s},
\newblock \bibinfo{journal}{Math. Comp.} \bibinfo{volume}{89} (\bibinfo{year}{2020}) \bibinfo{pages}{1547--1579}. \DOIprefix\doi{10.1090/mcom/3514}.
\bibitem[{Kingma and Ba(2017)}]{Kingma2017Adam}
\bibinfo{author}{D.~P. Kingma}, \bibinfo{author}{J.~Ba}, \bibinfo{title}{Adam: A method for stochastic optimization}, \bibinfo{year}{2017}. \href{http://arxiv.org/abs/1412.6980}{{\tt arXiv:1412.6980}}.
\bibitem[{Kloeden and Platen(1992)}]{Kloeden1992Numerical}
\bibinfo{author}{P.~E. Kloeden}, \bibinfo{author}{E.~Platen}, \bibinfo{title}{Numerical solution of stochastic differential equations}, volume~\bibinfo{volume}{23} of \textit{\bibinfo{series}{Applications of Mathematics (New York)}}, \bibinfo{publisher}{Springer-Verlag, Berlin}, \bibinfo{year}{1992}. \DOIprefix\doi{10.1007/978-3-662-12616-5}.
\bibitem[{Li et~al.(2024{\natexlab{a}})Li, Wang, Ye, He, and Wang}]{Li2024DOF}
\bibinfo{author}{R.~Li}, \bibinfo{author}{C.~Wang}, \bibinfo{author}{H.~Ye}, \bibinfo{author}{D.~He}, \bibinfo{author}{L.~Wang},
\newblock \bibinfo{title}{{DOF}: Accelerating high-order differential operators with forward propagation},
\newblock in: \bibinfo{booktitle}{ICLR 2024 Workshop on AI4DifferentialEquations In Science}, \bibinfo{year}{2024}{\natexlab{a}}. \URLprefix \url{https://openreview.net/forum?id=yQsLKpkRoS}.
\bibitem[{Li et~al.(2024{\natexlab{b}})Li, Ye, Jiang, Wen, Wang, Li, Li, He, Chen, Ren, and Wang}]{Li2024computational}
\bibinfo{author}{R.~Li}, \bibinfo{author}{H.~Ye}, \bibinfo{author}{D.~Jiang}, \bibinfo{author}{X.~Wen}, \bibinfo{author}{C.~Wang}, \bibinfo{author}{Z.~Li}, \bibinfo{author}{X.~Li}, \bibinfo{author}{D.~He}, \bibinfo{author}{J.~Chen}, \bibinfo{author}{W.~Ren}, \bibinfo{author}{L.~Wang},
\newblock \bibinfo{title}{A computational framework for neural network-based variational monte carlo with forward laplacian},
\newblock \bibinfo{journal}{Nat. Mach. Intell.} \bibinfo{volume}{6} (\bibinfo{year}{2024}{\natexlab{b}}) \bibinfo{pages}{209--219}. \DOIprefix\doi{10.1038/s42256-024-00794-x}.
\bibitem[{Li et~al.(2024{\natexlab{c}})Li, Verma, and Ruthotto}]{Li2024Neural}
\bibinfo{author}{X.~Li}, \bibinfo{author}{D.~Verma}, \bibinfo{author}{L.~Ruthotto},
\newblock \bibinfo{title}{A neural network approach for stochastic optimal control},
\newblock \bibinfo{journal}{SIAM J. Sci. Comput.} \bibinfo{volume}{46} (\bibinfo{year}{2024}{\natexlab{c}}) \bibinfo{pages}{C535--C556}. \DOIprefix\doi{10.1137/23M155832X}.
\bibitem[{Liu et~al.(2020)Liu, Cai, and Xu}]{Liu2020Multi}
\bibinfo{author}{Z.~Liu}, \bibinfo{author}{W.~Cai}, \bibinfo{author}{Z.-Q.~J. Xu},
\newblock \bibinfo{title}{Multi-scale deep neural network ({M}scale{DNN}) for solving {P}oisson-{B}oltzmann equation in complex domains},
\newblock \bibinfo{journal}{Commun. Comput. Phys.} \bibinfo{volume}{28} (\bibinfo{year}{2020}) \bibinfo{pages}{1970--2001}. \DOIprefix\doi{10.4208/cicp.oa-2020-0179}.
\bibitem[{N\"usken and Richter(2023)}]{Nusken2023Interpolating}
\bibinfo{author}{N.~N\"usken}, \bibinfo{author}{L.~Richter},
\newblock \bibinfo{title}{Interpolating between {BSDE}s and {PINN}s: deep learning for elliptic and parabolic boundary value problems},
\newblock \bibinfo{journal}{J. Mach. Learn.} \bibinfo{volume}{2} (\bibinfo{year}{2023}) \bibinfo{pages}{31--64}.
\bibitem[{\O~ksendal(2003)}]{Oksendal2003Stochastic}
\bibinfo{author}{B.~\O~ksendal}, \bibinfo{title}{Stochastic differential equations}, Universitext, \bibinfo{edition}{sixth} ed., \bibinfo{publisher}{Springer-Verlag, Berlin}, \bibinfo{year}{2003}. \DOIprefix\doi{10.1007/978-3-642-14394-6}, \bibinfo{note}{an introduction with applications}.
\bibitem[{Pavliotis(2014)}]{pavliotis2014stochastic}
\bibinfo{author}{G.~A. Pavliotis}, \bibinfo{title}{Stochastic Processes and Applications: Diffusion Processes, the Fokker-Planck and Langevin Equations}, Texts in Applied Mathematics, \bibinfo{publisher}{Springer New York}, \bibinfo{year}{2014}.
\bibitem[{Raissi(2018)}]{raissi2018forwardbackward}
\bibinfo{author}{M.~Raissi}, \bibinfo{title}{Forward-backward stochastic neural networks: Deep learning of high-dimensional partial differential equations}, \bibinfo{year}{2018}. \href{http://arxiv.org/abs/arXiv:1804.07010 [stat.ML]}{{\tt arXiv:arXiv:1804.07010 [stat.ML]}}.
\bibitem[{Raissi et~al.(2019)Raissi, Perdikaris, and Karniadakis}]{Raissi2019Physics}
\bibinfo{author}{M.~Raissi}, \bibinfo{author}{P.~Perdikaris}, \bibinfo{author}{G.~E. Karniadakis},
\newblock \bibinfo{title}{Physics-informed neural networks: a deep learning framework for solving forward and inverse problems involving nonlinear partial differential equations},
\newblock \bibinfo{journal}{J. Comput. Phys.} \bibinfo{volume}{378} (\bibinfo{year}{2019}) \bibinfo{pages}{686--707}. \DOIprefix\doi{10.1016/j.jcp.2018.10.045}.
\bibitem[{Shi et~al.(2024)Shi, Hu, Lin, and Kawaguchi}]{shi2024stochastic}
\bibinfo{author}{Z.~Shi}, \bibinfo{author}{Z.~Hu}, \bibinfo{author}{M.~Lin}, \bibinfo{author}{K.~Kawaguchi},
\newblock \bibinfo{title}{Stochastic taylor derivative estimator: Efficient amortization for arbitrary differential operators},
\newblock in: \bibinfo{booktitle}{The Thirty-eighth Annual Conference on Neural Information Processing Systems}, \bibinfo{year}{2024}.
\bibitem[{Sirignano and Spiliopoulos(2018)}]{Sirignano2018DGM}
\bibinfo{author}{J.~Sirignano}, \bibinfo{author}{K.~Spiliopoulos},
\newblock \bibinfo{title}{D{GM}: a deep learning algorithm for solving partial differential equations},
\newblock \bibinfo{journal}{J. Comput. Phys.} \bibinfo{volume}{375} (\bibinfo{year}{2018}) \bibinfo{pages}{1339--1364}. \DOIprefix\doi{10.1016/j.jcp.2018.08.029}.
\bibitem[{Stein(1981)}]{Stein1981Estimation}
\bibinfo{author}{C.~M. Stein},
\newblock \bibinfo{title}{Estimation of the mean of a multivariate normal distribution},
\newblock \bibinfo{journal}{Ann. Statist.} \bibinfo{volume}{9} (\bibinfo{year}{1981}) \bibinfo{pages}{1135--1151}.
\bibitem[{Taylor(2023)}]{Taylor2023PDEIII}
\bibinfo{author}{M.~E. Taylor}, \bibinfo{title}{Partial differential equations {III}. {Nonlinear} equations}, volume \bibinfo{volume}{117} of \textit{\bibinfo{series}{Appl. Math. Sci.}}, \bibinfo{edition}{3rd corrected and expanded edition} ed., \bibinfo{publisher}{Cham: Springer}, \bibinfo{year}{2023}. \DOIprefix\doi{10.1007/978-3-031-33928-8}.
\bibitem[{Wang et~al.(2022{\natexlab{a}})Wang, Li, He, and Wang}]{wang2022is}
\bibinfo{author}{C.~Wang}, \bibinfo{author}{S.~Li}, \bibinfo{author}{D.~He}, \bibinfo{author}{L.~Wang},
\newblock \bibinfo{title}{Is \$l{\textasciicircum}2\$ physics informed loss always suitable for training physics informed neural network?},
\newblock in: \bibinfo{editor}{A.~H. Oh}, \bibinfo{editor}{A.~Agarwal}, \bibinfo{editor}{D.~Belgrave}, \bibinfo{editor}{K.~Cho} (Eds.), \bibinfo{booktitle}{Advances in Neural Information Processing Systems}, \bibinfo{year}{2022}{\natexlab{a}}.
\bibitem[{Wang et~al.(2022{\natexlab{b}})Wang, Zhao, and Zhou}]{Wang2022Sinc}
\bibinfo{author}{X.~Wang}, \bibinfo{author}{W.~Zhao}, \bibinfo{author}{T.~Zhou},
\newblock \bibinfo{title}{Sinc-{$\theta$} schemes for backward stochastic differential equations},
\newblock \bibinfo{journal}{SIAM J. Numer. Anal.} \bibinfo{volume}{60} (\bibinfo{year}{2022}{\natexlab{b}}) \bibinfo{pages}{1799--1823}. \DOIprefix\doi{10.1137/21M1444679}.
\bibitem[{Xu and Zhang(2025)}]{zhang2025shotgun}
\bibinfo{author}{W.~Xu}, \bibinfo{author}{W.~Zhang},
\newblock \bibinfo{title}{A deep shotgun method for solving high-dimensional parabolic partial differential equations},
\newblock \bibinfo{journal}{J. Sci. Comput.} \bibinfo{volume}{104} (\bibinfo{year}{2025}) \bibinfo{pages}{69}. \DOIprefix\doi{10.1007/s10915-025-02983-1}.
\bibitem[{Yang et~al.(2020)Yang, Zhao, and Zhou}]{Yang2020unified}
\bibinfo{author}{J.~Yang}, \bibinfo{author}{W.~Zhao}, \bibinfo{author}{T.~Zhou},
\newblock \bibinfo{title}{A unified probabilistic discretization scheme for {FBSDE}s: stability, consistency, and convergence analysis},
\newblock \bibinfo{journal}{SIAM J. Numer. Anal.} \bibinfo{volume}{58} (\bibinfo{year}{2020}) \bibinfo{pages}{2351--2375}. \DOIprefix\doi{10.1137/19M1260177}.
\bibitem[{Yong and Zhou(1999)}]{Yong1999Stochastic}
\bibinfo{author}{J.~Yong}, \bibinfo{author}{X.~Y. Zhou}, \bibinfo{title}{Stochastic controls}, volume~\bibinfo{volume}{43} of \textit{\bibinfo{series}{Applications of Mathematics (New York)}}, \bibinfo{publisher}{Springer-Verlag, New York}, \bibinfo{year}{1999}. \DOIprefix\doi{10.1007/978-1-4612-1466-3}, \bibinfo{note}{hamiltonian systems and HJB equations}.
\bibitem[{Yuan and Mao(2008)}]{Yuan2008note}
\bibinfo{author}{C.~Yuan}, \bibinfo{author}{X.~Mao},
\newblock \bibinfo{title}{A note on the rate of convergence of the {Euler}-{Maruyama} method for stochastic differential equations},
\newblock \bibinfo{journal}{Stochastic Anal. Appl.} \bibinfo{volume}{26} (\bibinfo{year}{2008}) \bibinfo{pages}{325--333}. \DOIprefix\doi{10.1080/07362990701857251}.
\bibitem[{Zang et~al.(2020)Zang, Bao, Ye, and Zhou}]{Zang2020Weak}
\bibinfo{author}{Y.~Zang}, \bibinfo{author}{G.~Bao}, \bibinfo{author}{X.~Ye}, \bibinfo{author}{H.~Zhou},
\newblock \bibinfo{title}{Weak adversarial networks for high-dimensional partial differential equations},
\newblock \bibinfo{journal}{J. Comput. Phys.} \bibinfo{volume}{411} (\bibinfo{year}{2020}) \bibinfo{pages}{109409, 14}. \DOIprefix\doi{10.1016/j.jcp.2020.109409}.
\bibitem[{Zhang(2004)}]{Zhang2004numerical}
\bibinfo{author}{J.~Zhang},
\newblock \bibinfo{title}{A numerical scheme for {BSDE}s},
\newblock \bibinfo{journal}{Ann. Appl. Probab.} \bibinfo{volume}{14} (\bibinfo{year}{2004}) \bibinfo{pages}{459--488}. \DOIprefix\doi{10.1214/aoap/1075828058}.
\bibitem[{Zhang(2017)}]{Zhang2017Backward}
\bibinfo{author}{J.~Zhang}, \bibinfo{title}{Backward stochastic differential equations}, volume~\bibinfo{volume}{86} of \textit{\bibinfo{series}{Probability Theory and Stochastic Modelling}}, \bibinfo{publisher}{Springer, New York}, \bibinfo{year}{2017}. \DOIprefix\doi{10.1007/978-1-4939-7256-2}, \bibinfo{note}{from linear to fully nonlinear theory}.
\bibitem[{Zhang and Cai(2022)}]{Zhang2022FBSDE}
\bibinfo{author}{W.~Zhang}, \bibinfo{author}{W.~Cai},
\newblock \bibinfo{title}{F{BSDE} based neural network algorithms for high-dimensional quasilinear parabolic {PDE}s},
\newblock \bibinfo{journal}{J. Comput. Phys.} \bibinfo{volume}{470} (\bibinfo{year}{2022}) \bibinfo{pages}{Paper No. 111557, 14}. \DOIprefix\doi{10.1016/j.jcp.2022.111557}.
\bibitem[{Zhao et~al.(2006)Zhao, Chen, and Peng}]{zhao2006new}
\bibinfo{author}{W.~Zhao}, \bibinfo{author}{L.~Chen}, \bibinfo{author}{S.~Peng},
\newblock \bibinfo{title}{A new kind of accurate numerical method for backward stochastic differential equations},
\newblock \bibinfo{journal}{SIAM J. Sci. Comput.} \bibinfo{volume}{28} (\bibinfo{year}{2006}) \bibinfo{pages}{1563--1581}. \DOIprefix\doi{10.1137/05063341X}.
\bibitem[{Zhao et~al.(2014)Zhao, Fu, and Zhou}]{zhao2014new}
\bibinfo{author}{W.~Zhao}, \bibinfo{author}{Y.~Fu}, \bibinfo{author}{T.~Zhou},
\newblock \bibinfo{title}{New kinds of high-order multistep schemes for coupled forward backward stochastic differential equations},
\newblock \bibinfo{journal}{SIAM J. Sci. Comput.} \bibinfo{volume}{36} (\bibinfo{year}{2014}) \bibinfo{pages}{A1731--A1751}. \DOIprefix\doi{10.1137/130941274}.
\bibitem[{Zhao et~al.(2012)Zhao, Li, and Zhang}]{Zhao2012A}
\bibinfo{author}{W.~Zhao}, \bibinfo{author}{Y.~Li}, \bibinfo{author}{G.~Zhang},
\newblock \bibinfo{title}{A generalized {$\theta$}-scheme for solving backward stochastic differential equations},
\newblock \bibinfo{journal}{Discrete Contin. Dyn. Syst. Ser. B} \bibinfo{volume}{17} (\bibinfo{year}{2012}) \bibinfo{pages}{1585--1603}. \DOIprefix\doi{10.3934/dcdsb.2012.17.1585}.
\bibitem[{Zhao et~al.(2010)Zhao, Zhang, and Ju}]{zhao2010stable}
\bibinfo{author}{W.~Zhao}, \bibinfo{author}{G.~Zhang}, \bibinfo{author}{L.~Ju},
\newblock \bibinfo{title}{A stable multistep scheme for solving backward stochastic differential equations},
\newblock \bibinfo{journal}{SIAM J. Numer. Anal.} \bibinfo{volume}{48} (\bibinfo{year}{2010}) \bibinfo{pages}{1369--1394}. \DOIprefix\doi{10.1137/09076979X}.
\bibitem[{Zhou et~al.(2021)Zhou, Han, and Lu}]{Zhou2021Actor}
\bibinfo{author}{M.~Zhou}, \bibinfo{author}{J.~Han}, \bibinfo{author}{J.~Lu},
\newblock \bibinfo{title}{Actor-critic method for high dimensional static {H}amilton-{J}acobi-{B}ellman partial differential equations based on neural networks},
\newblock \bibinfo{journal}{SIAM J. Sci. Comput.} \bibinfo{volume}{43} (\bibinfo{year}{2021}) \bibinfo{pages}{A4043--A4066}. \DOIprefix\doi{10.1137/21M1402303}.

\end{thebibliography}

\end{document}